\documentclass[11pt]{amsart}
\usepackage{bbm}
\usepackage{amssymb, amsthm, amsmath, amscd}
\usepackage{enumerate}
\usepackage[usenames,dvipsnames]{color}
\usepackage{hyperref}

\theoremstyle{definition}
\newtheorem{example}{Example}[section]      
\newtheorem{definition}[example]{Definition}
\newtheorem{theorem}[example]{Theorem}
\newtheorem{remark}[example]{Remark}
\newtheorem{proposition}[example]{Proposition}
\newtheorem{lemma}[example]{Lemma}
\newtheorem{corollary}[example]{Corollary}
\newtheorem{notation}[example]{Notation}
\numberwithin{equation}{section}

\usepackage{amsfonts}
\usepackage{mathrsfs}
\usepackage{textcomp}
\usepackage[all]{xy}

\title[The weak tracial Rokhlin property]{\bf Crossed products by compact group actions with the weak tracial Rokhlin property}

\author{Xiaochun Fang}
\address{School of Mathematical Sciences,
	Key Laboratory of Intelligent Computing and Applications(Ministry of Education), 
	Tongji University, Shanghai 200092, CHINA}
\email{xfang@tongji.edu.cn}

\author{Haotian Tian}
\address{School of Mathematical Sciences,
	Tongji University, Shanghai 200092, CHINA}
\email{2011284@tongji.edu.cn}
\thanks{Corresponding author: Haotian Tian}


\begin{document}
	\begin{abstract}
		In this paper, we introduce compact group actions with the weak tracial Rokhlin property. This concept generalizes both the finite group actions with the weak tracial Rokhlin property and the compact group actions with the tracial Rokhlin property on classifiable C*-algebras (in the sense of the Elliott program). Under this framework, we prove that simplicity, pure infiniteness, tracial $\mathcal{Z}$-stability and the combination of nuclearity and $\mathcal{Z}$-stability can be transferred from the original algebra to the crossed product. We also show that the radius of comparison of the fixed point algebra does not exceed that of the original algebra. Furthermore, we discuss the relationship between our definition and natural generalization of the finite group case in non-Elliott program settings. Finally, we provide a nontrivial example of a compact group action with the weak tracial Rokhlin property with comparison: an action of $(S_2)^\mathbb{N}$ on the Jiang-Su algebra $\mathcal{Z}$. Since $\mathcal{Z}$ contains no nontrivial projections, this action does not possess the tracial Rokhlin property.
	\end{abstract}
	
	\maketitle
	
	\tableofcontents
	
	\section{Introduction}
	The Rokhlin property for the case of a single automorphism was originally introduced for von Neumann algebras by Connes in \cite{Con1975}. Later, the Rokhlin property for finite group actions on C*-algebras first appeared in the work of Herman and Jones in \cite{HJ1982} and \cite{HJ1983}. This property is useful to understand the structure of the crossed product of C*-algebras and properties passing from the original algebra to the crossed product \cite{OP2012}. However, the finite group acitions with the Rokhlin property are rare. Phillips, in \cite{Phi2011}, introduced the tracial Rokhlin property for finite group actions on unital simple C*-algebras. The tracial Rokhlin property is generic in many cases, and also can be used to study properties passing from the original algebra to the crossed product. Weak versions of the tracial Rokhlin property in which one uses orthogonal positive contractions instead of orthogonal projections were studied for actions on unital simple C*-algebras with few projections \cite{GHS2021,HO2013,MS2012b,Phi2012,Wan2013,Wan2018} (see Definition \ref{WTRF}). As an example, the flip action on the Jiang-Su algebra $\mathcal{Z}\cong \mathcal{Z}\otimes\mathcal{Z}$ has the weak tracial Rokhlin property but it does not have the tracial Rokhlin property \cite{HO2013}. For the non-unital case, Santiago and Gardella studied the Rokhlin property for finite group actions on non-unital simple C*-algebras in \cite{San2015} and \cite{GS2016}. Forough and Golestani studied the (weak) tracial Rokhlin property for finite group actions on non-unital simple C*-algebras in \cite{FG2020}.
	
	In \cite{HW2007}, Hirshberg and Winter also introduced the Rokhlin property for second-countable compact group actions on unital C*-algebras. Since then, crossed products by compact group actions with the Rokhlin property have been studied by several authors. In particular, permanence properties are proved in \cite{HW2007}, \cite{Gar2017} and \cite{Gar2019}. The same as finite groups, Rokhlin actions of compact groups are rare, especially when the group is connected. More recently, Mohammadkarimi and Phillips studied the tracial Rokhlin property with comparison for compact group actions and proved that the crossed product of a unital separable simple infinite dimensional C*-algebra with tracial rank zero by an action of a second-countable compact group with the tracial Rokhlin property with comparison has again tracial rank zero in \cite{MP2021} and some other permanence properties. Moreover, they gave some examples of compact group actions with the tracial Rokhlin property with comparison. The authors have studied many permanence properties in \cite{TF2025} including stable rank one, real rank zero, $\beta$-comparison, Winter's $n$-comparison, $m$-almost divisibility and weakly ($m$,$n$)-divisibility.
	
	In this paper, We introduce compact group actions with the weak tracial Rokhlin property. This concept generalizes both the finite group actions with the weak tracial Rokhlin property and the compact group actions with the tracial Rokhlin property on classifiable C*-algebras (in the sense of the Elliott program). As in the case of finite groups, it is impossible for a compact group action on a C*-algebra containing no non-trivial projections (such as the Jiang-Su algebra $\mathcal{Z}$) to have the tracial Rokhlin property. Therefore, this generalization has practical research value. Under this concept, we prove that simplicity, pure infiniteness, tracial $\mathcal{Z}$-stability, and the combination of amenability and $\mathcal{Z}$-stability can be transferred from the original algebra to the crossed product algebra. We also prove that the radius of comparison of the fixed point algebra does not exceed that of the original algebra. It is worth mentioning that during the research process, we did not directly study the crossed product algebra, but instead employed the so-called ``averaging process'' introduced by Gardella in \cite{Gar2014}. That is, we first study how structural properties transfer from the original algebra to the fixed point algebra, and then utilize the fact that the fixed point algebra and the crossed product algebra are Morita equivalent (and hence stably isomorphic in the separable case) to obtain the structural properties of the crossed product algebra. Simultaneously, we discuss the relationship between our definition and the natural generalization of the finite group case in a non-Elliott program sense. Finally, we provide a non-trivial example of a compact group action with the weak tracial Rokhlin property and finite radius of comparison. This is an action of $(S_2)^\mathbb{N}$ on the Jiang-Su algebra $\mathcal{Z}$. Since the Jiang-Su algebra $\mathcal{Z}$ contains no non-trivial projections, this action does not have the tracial Rokhlin property.
	
	
	The paper is organized as follows. Section \ref{sec8.0} contains some preliminaries about central sequence algebras, Cuntz subequivalence and some other lemmas that we will use in this paper. Seciton \ref{sec8.1} contains the definition of the weak tracial Rokhlin property with comparison for compact group actions, and its basic properties. In particular, we give a central sequence formulation of the definition, and give an averaging process which is the key technical tool for proofs of permanence properties.
	
	In Section \ref{sec8.2}, we prove that the fixed point algebra and the crossed product of an infinite dimensional simple separable unital C*-algebra by an action of a compact group with the weak tracial Rokhlin property with comparison are again simple. The rest of our permanence properties are in Section \ref{sec8.3}.
	
	Section \ref{sec8.4} discusses the relation between the weak tracial Rokhlin property with comparison and the weak tracial Rokhlin property. In particular, under some reasonable conditions on the algebra, we prove that they are equivalent.
	
	In Section \ref{sec8.5}, we construct an example of action which has the weak tracial Rokhlin property with comparison. It is an action of $(S_2)^\mathbb{N}$ on the Jiang-Su algebra $\mathcal{Z}$. This example does not have the tracial Rokhlin property. In fact, through this example and other properties, we present a general method for constructing actions with the weak tracial Rokhlin property, not limited to constructing this specific example.
	
	\section{Preliminaries and Definitions}
	\label{sec8.0}
	In this section, we recall some definitions and known facts about central sequence algebras, Cuntz subequivalence, completely positive order zero maps and other notations.
	
	\subsection{Central sequence algebras}
	\begin{definition}
		Let $A$ be a unital C*-algebra. We use $l^\infty(\mathbb{N},A)$ to denote the set of all bounded sequences in $A$ with the supremum norm which is a unital C*-algebra with the unit as the constant sequence $1$. Let
		\[c_0(\mathbb{N},A)=\{(a_n)_{n\in\mathbb{N}}\in l^\infty(\mathbb{N},A):\lim_{n\to\infty}\|a_n\|=0\}.\]
		It is obvious that $c_0(\mathbb{N},A)$ is a closed two-side ideal in $l^\infty(\mathbb{N},A)$, and we use the notation $A_\infty$ to denote the quotient $l^\infty(\mathbb{N},A)/c_0(\mathbb{N},A)$. Denoted by $\pi_A\colon l^\infty(\mathbb{N},A)\to A_\infty$ the quotient map. Define $\iota\colon A\to l^\infty(\mathbb{N},A)$ by $\iota(a)=(a,a,a,\dots)$, the constant sequence, for all $a\in A$. Identify $A$ with $\pi_A\circ\iota(A)$. Denoted by $A_\infty\cap A'$ the relative commutant of $A$ inside of $A_\infty$.
		
		let $G$ be a topological group and let $\alpha\colon G\to \mathrm{Aut}(A)$ be an action of $G$ on $A$, then it induced actions of $G$ on $l^\infty(\mathbb{N},A)$ and on $A_\infty$, denoted by $\alpha^\infty$ and $\alpha_\infty$. Since for any $g\in G$,
		\[(\alpha_\infty)_g(A_\infty\cap A')\subseteq A_\infty\cap A',\]
		so we also use $\alpha_\infty$ to denote the restricted action on $A_\infty\cap A'$. These actions are not necessarily continuous when $G$ is not discrete. Therefore, we set
		\[l^\infty_\alpha(\mathbb{N},A)=\{a\in l^\infty(\mathbb{N},A):g\mapsto\alpha^\infty_g(a)\ \mathrm{is}\ \mathrm{continuous}\},\]
		and $A_{\infty,\alpha}=\pi_A(l^\infty_\alpha(\mathbb{N},A))$. Then, $A_{\infty,\alpha}$ is invariant under $\alpha_\infty$ and the action $\alpha_\infty$ is continuous.
	\end{definition}
	
	\subsection{Cuntz comparison}
	\begin{definition}
		Let $A$ be a C*-algebra, $a\in A_+$ and $\varepsilon>0$. Then we denote $f(a)$ by $(a-\varepsilon)_+$, where $f(t)=max\{0,t-\varepsilon\}$ is continuous from $[0,\infty)$ to $[0,\infty)$.
	\end{definition}
	
	The following definitions related to Cuntz comparison are from \cite{Cun1978}, for more information, you can refer to \cite{GP2024} and \cite{APT2018}.
	
	\begin{definition}
		Let $A$ be a C*-algebra. Let $a,b\in (A\otimes K)_+$.
		\begin{enumerate}
			\item We say that $a$ is Cuntz subequivalent to $b$ (written $a\precsim_A b$), if there is a sequence $(r_n)_{n=1}^\infty$ in $A\otimes K$ such that $\lim\limits_{n\to\infty}\|r_n^*br_n-a\|=0$.
			
			\item We say that $a$ is Cuntz equivalent to $b$ (written $a\sim_A b$), if $a\precsim_A b$ and $b\precsim_A a$. This is an equivalence relation, we use $\langle a\rangle_A$ to denote the equivalence class of $a$. With the addition operation $\langle a\rangle_A +\langle b\rangle_A=\langle a\oplus b\rangle_A$ and the order operation $\langle a\rangle_A \leq\langle b\rangle_A$ if $a\precsim_A b$, $\mathrm{Cu}(A)=(A\otimes K)_+/\sim_A$ is an ordered semigroup which we called Cuntz semigroup. $\mathrm{W}(A)=M_\infty(A)_+/\sim_A$ is also an ordered semigroup with the same operation and order as above.
		\end{enumerate}
		If $B$ is a hereditary C*-subalgebra of $A$, and $a,b\in B_+$, then it is easy to check that $a\precsim_A b\iff a\precsim_B b$.
	\end{definition}
	
	\begin{remark}
		Adopting the usual identifications $A \subset M_n(A) \subset M_\infty(A) \subset A \otimes K$. If $a, b \in A_+$ and $a \precsim_A b$, then we assert that there exists a sequence $(r_n)_{n=1}^\infty \subset A$ such that $\lim_{n \to \infty} \| r_n^* b r_n - a \| = 0$. Indeed, choose a sequence $(v_n)_{n=1}^\infty$ taking values in $A \otimes K$ satisfying $\lim_{n \to \infty} \| v_n^* b v_n - a \| = 0$. Let $(e_{j,k})_{j,k \in \mathbb{N}}$ be the standard matrix units of $K$, and set $r_n = (e_{1,1} \otimes 1)v_n (e_{1,1} \otimes 1)$.
		
		For the same reason, if $a,b\in M_n(A)_+$ for some $n\in\mathbb{N}$, then the sequence $(r_n)_{n=1}^\infty$ can be taken in $M_n(A)$; replacing $M_n(A)$ with $M_\infty(A)$, the same conclusion holds. Therefore, if $a$ and $b$ belong to any one of $A_+, M_n(A)_+, M_\infty(A)_+, (A\otimes K)_+$ (they do not need to belong to the same algebra), we write $a\precsim_A b$ (or $a\sim_A b$) to mean that this equivalence relation holds in $A\otimes K$. This is equivalent to the relation holding in the smallest containing subalgebra (i.e., one of $A, M_n(A), M_\infty(A), A\otimes K$) that contains both $a$ and $b$.
	\end{remark}
	
	Below we give some well-known facts about Cuntz comparison. (1) is contained in \cite[Proposition 2.6]{KR2000} (and also appears in a slightly different form in the earlier \cite[Proposition 2.4]{Ror1992}). (2) is \cite[Lemma 1.7]{Phi2014}. (3) is straightforward (and is also \cite[Lemma 2.5(i)]{KR2000}). (4) is \cite[Corollary 1.6]{Phi2014}. (5) is \cite[Lemma 2.8(ii)]{KR2000}. (6) is \cite[Lemma 2.8(iii)]{KR2000}. (7) is \cite[Lemma 2.9]{KR2000}. (8) is \cite[Lemma 2.4(2)]{KR2002}. For (9), it is straightforward to verify that $a$ and $f(a)$ generate the same hereditary subalgebra in $A$. (10) is contained in the discussion following \cite[Definition 2.3]{KR2000}.
	\begin{lemma}\label{CCP}
		Let $A$ be a C*-algebra.
		\begin{enumerate}[{(}1{)}]
			\item Let $a,b\in A_+$. Then the following are equivalent:
			\begin{enumerate}
				\item $a\precsim_A b$;
				\item $(a-\varepsilon)_+\precsim_A b$ for all $\varepsilon>0$;
				\item For all $\varepsilon>0$, there exists $\delta>0$ such that $(a-\varepsilon)_+\precsim_A (b-\delta)_+$;
				\item For all $\varepsilon>0$, there exist $\delta>0$ and $x\in A$ such that $(a-\varepsilon)_+=x^*(b-\delta)_+x$.
			\end{enumerate}
			\item Let $a,b\in A_+$. If $a\leq b$ and $\lambda\geq 0$, then $(a-\lambda)_+\precsim_A (b-\lambda)_+$.
			\item Let $a\in A_+$ and $\varepsilon_1,\varepsilon_2>0$. Then $((a-\varepsilon_1)_+-\varepsilon_2)_+=(a-(\varepsilon_1+\varepsilon_2))_+$.
			\item Let $\varepsilon>0$, $\lambda\geq 0$, and $a,b\in A_+$. If $\|a-b\|<\varepsilon$, then $(a-\varepsilon-\lambda)_+\precsim_A(b-\lambda)_+$.
			\item Let $a,b\in A_+$, then $a+b\precsim_A a\oplus b$.
			\item Let $a,b\in A_+$ with $ab=0$, then $a+b\sim_A a\oplus b$.
			\item Let $a_1,a_2,b_1,b_2\in A_+$. If $a_1\precsim_A a_2$ and $b_1\precsim_A b_2$, then $a_1\oplus b_1\precsim a_2\oplus b_2$.
			\item Let $a,b\in A_+$ and $\delta>0$. If $a\precsim_A (b-\delta)_+$, then for each $r>1$ there exists $y\in A$ with $a=y^*by$ and $\|y\|\leq r\delta^{-\frac{1}{2}}\|a\|^\frac{1}{2}$.
			\item Let $a\in A_+$ and let $f\colon[0,\|a\|]\to[0,\infty)$ be a continuous function satisfying $f(0)=0$ and $f(\lambda)>0$ for $\lambda>0$. Then $f(a)\sim_A a$.
			\item Let $c\in A$ and $\varepsilon\geq0$. Then $(c^*c-\varepsilon)_+\sim_A (cc^*-\varepsilon)_+$.
		\end{enumerate}
	\end{lemma}
	
	Lemma \ref{CCP} does not provide all the tools concerning Cuntz comparison needed for this paper. Below we supplement some equally well-known instrumental lemmas.
	\begin{lemma}\cite[Lemma 2.4]{Phi2014}\label{OCP}
		Let $A$ be a simple C*-algebra which is not of type I. Let $a\in A$ be a nonzero positive element, and $l$ a positive integer. Then there exist $l$ mutually orthogonal nonzero positive elements $b_1,\ldots,b_l\in A$ such that $b_1\sim_A\cdots \sim_A b_l$, and $b_1+\cdots + b_l\in\overline{aAa}$.
	\end{lemma}
	
	\begin{lemma}\cite[Lemma 1.5]{Phi2014}\label{CPA}
		Let $A$ be a C*-algebra, let $a,b\in A_+$, and let $\alpha,\beta\geq0$. Then
		\[((a+b)-(\alpha+\beta))_+\precsim_A(a-\alpha)_++(b-\beta)_+\precsim_A(a-\alpha)_+\oplus(b-\beta)_+\]
	\end{lemma}
	
	The following lemma is known as Kirchberg's Slice Lemma, see \cite[Lemma 4.1.9]{RS2002}.
	\begin{lemma}\label{KS}
		Let $A$ and $B$ be C*-algebras, and let $D$ be a nonzero hereditary subalgebra of the minimal tensor product $A\otimes B$. Then there exists a nonzero positive element $z\in A\otimes B$ such that $z^*z$ is of the form $a\otimes b$ for some $a\in A$, $b\in B$, and $zz^*\in D$.
	\end{lemma}
	
	Below we introduce the strict comparison of positive elements and the radius of comparison. The following definition of tracial state is from \cite{Lin2001}, and the definition of 2-quasitrace is from \cite{BH1982}.
	\begin{definition}
		Let $A$ be a C*-algebra. A linear functional $\tau\colon A\to \mathbb{C}$ is called a \textit{trace} on $A$ if $\tau$ is positive and satisfies $\tau(xy)=\tau(yx)$ for all $x,y\in A$. If furthermore $\|\tau\|=1$, then $\tau$ is called a \textit{tracial state}. The set of all tracial states of $A$ is denoted by $\mathrm{T}(A)$.
		
		In fact, for any approximate unit $(e_\lambda)_{\lambda\in \Lambda}$ of $A$, $\|\tau\|=\lim_\lambda \tau(e_\lambda)$. In particular, if $A$ has a unit, then $\|\tau\|=\tau(1_A)$.
	\end{definition}
	
	\begin{definition}
		Let $A$ be a C*-algebra. A functional $\tau\colon A\to\mathbb{C}$ is called a \textit{quasitrace} on $A$ if it satisfies:
		\begin{enumerate}[{(}1{)}]
			\item $\tau(x^*x)=\tau(xx^*)\geq 0$ for all $x\in A$;
			\item $\tau$ is linear on commutative subalgebras of $A$;
			\item If $x=a+ib$ where $a,b$ are self-adjoint, then $\tau(x)=\tau(a)+i\tau(b)$.
		\end{enumerate}
		$\tau$ is called a \textit{$n$-quasitrace} ($n\geq 2$) if it can be extended to a quasitrace on $M_n(A)$. $\tau$ is called \textit{normalized} if $\|\tau\|=1$.
	\end{definition}
	
	From \cite{BH1982}, it follows that for a unital C*-algebra, a 2-quasitrace $\tau$ automatically extends to an $n$-quasitrace for all $n\geq 2$. Denote by $\mathrm{QT}(A)$ the set of all normalised 2-quasitraces on $A$. It is proved in \cite{Han1981} that if $A$ is a stably finite C*-algebra, then $\mathrm{QT}(A)\neq \varnothing$.
	\begin{remark}
		\begin{enumerate}[{(}1{)}]
			\item Let $A$ be a unital C*-algebra. For all $a\in M_\infty(A)_+$ and $\tau\in \mathrm{QT}(A)$, define $d_\tau\colon M_\infty(A)\to [0,\infty)$ by $d_\tau(a)=\lim_{n\to\infty}\tau(a^\frac{1}{n})$.
			\item For $a\in(A\otimes K)_+$, $\langle a\rangle\in \mathrm{Cu}(A)$ and $\langle a\rangle\in\mathrm{W}(A)$, we use the same notation for the corresponding functional. From part of the proof of \cite[Proposition 4.2]{ERS2011}, it follows that $d_\tau$ is well-defined on $\mathrm{Cu}(A)$ and $\mathrm{W}(A)$.
			\item By \cite[Theorem 2.32]{APT2009}, $d_\tau$ defines a state on $\mathrm{W}(A)$.
			\item According to \cite[Theorem I.2.2]{BH1982}, for any $\tau\in \mathrm{QT}(A)$, $d_\tau$ is a lower semicontinuous function on $A$. Consequently, for any $a\in A_+\setminus \{0\}$ we have
			\[d_\tau(a)=\lim_{\varepsilon\searrow 0}d_\tau((a-\varepsilon)_+).\]
		\end{enumerate}
	\end{remark}
	
	\begin{definition}
		Let $A$ be a unital simple stably finite C*-algebra. We say that $A$ has the \textit{strict comparison of positive elements property} if whenever positive elements $a,b$ in $M_\infty(A)$ satisfy
		\[d_\tau(a)<d_\tau(b),\quad \forall \tau\in \mathrm{QT}(A),\]
		then $a\precsim b$.
	\end{definition}
	
	By replacing normalised 2-quasitraces with densely defined 2-quasitraces, the strict comparison property for positive elements can also be defined for non-unital C*-algebras (see \cite{Lin2023}). Recall that an ordered semigroup $S$ is called \textit{almost unperforated} if for any $x,y\in S$ and $n\in\mathbb{N}$ satisfying $(n+1)x\leq ny$, we have $x\leq y$. By \cite[Lemma 6.1]{TT2015}, for a simple C*-algebra, having strict comparison of positive elements is equivalent to its Cuntz semigroup being almost unperforated.
	
	\begin{definition}\cite[Definition 6.1]{Tom2006}
		Let $A$ be a unital stably finite C*-algebra. Let $r\in [0,\infty)$. We say that $A$ has \textit{$r$-comparison} if whenever positive elements $a,b$ in $M_\infty(A)$ satisfy $d_\tau(a)+r<d_\tau(b)$ for all $\tau\in \mathrm{QT}(A)$, then $a\precsim b$. Define the \textit{radius of comparison} of $A$ as
		\[\mathrm{rc}(A)\colon=\inf(\{r\in[0,\infty)\colon A\text{ has } r\text{-comparison}\}).\]
		If the infimum does not exist, we set $\mathrm{rc}(A)=\infty$.
	\end{definition}
	
	The above definition is equivalent to defining strict comparison and the radius of comparison using $\mathrm{W}(A)$; using $\mathrm{Cu}(A)$ yields the same result (see \cite[Proposition 6.12]{Phi2014}). If $A$ is simple and has finite radius of comparison, then the above infimum is attained, meaning $A$ has $\mathrm{rc}(A)$-comparison (see \cite[Proposition 6.3]{Tom2006} and \cite[Lemma 1.23]{Phi2014}).
	
	\subsection{Completely positive order zero map}
	We briefly review some core knowledge of completely positive order zero maps (for more details and conclusions, see \cite{WZ2009}).
	
	Let $A$ be a C*-algebra, and let $a,b\in A$. We say $a$ and $b$ are \emph{orthogonal}, denoted $a\perp b$, if $ab=ba=a^*b=ab^*=0$. If $a,b\in A$ are self-adjoint, then $a\perp b$ if and only if $ab=0$.
	
	\begin{definition}\cite[Definition 2.3]{WZ2009}
		Let $A$ and $B$ be C*-algebras and let $\varphi\colon A\to B$ be a completely positive map. We say $\varphi$ has \emph{order zero} if for every $a,b\in A$, we have $\varphi(a)\perp\varphi(b)$ whenever $a\perp b$.
	\end{definition}
	
	\begin{remark}
		It can be directly verified that homomorphisms between C*-algebras have order zero, and the composition of two order zero maps is again an order zero map.
	\end{remark}
	
	\begin{theorem}\cite[Theorem 2.2, Corollary 4.8]{GT2022}\label{OZAH1}
		Let $A$ and $B$ be C*-algebras and let $\varphi\colon A\to B$ be a completely positive map. Then $\varphi$ has order zero if and only if for any $a, b, c\in A$, we have
		\[\varphi(ab)\varphi(c)=\varphi(a)\varphi(bc).\]
		
		In particular, if $A$ has a unit, then $\varphi$ has order zero if and only if for any $a,b\in A$, we have
		\[\varphi(ab)\varphi(1)=\varphi(a)\varphi(b).\]
	\end{theorem}
	
	In the case where $A$ has a unit, the above theorem can be easily derived from the decomposition theorem for completely positive order zero maps; for the general case, one may refer to the relevant proof in \cite{GT2022}.
	
	\begin{theorem}\cite[Theorem 3.3]{WZ2009}\label{OZAH}
		Let $A$ and $B$ be C*-algebras and let $\varphi\colon A\to B$ be a completely positive order zero map. Let $C\colon = C^*(\varphi(A))\subset B$. Then there exists a positive element $h\in M(C)\cap C'$ satisfying $\|h\|=\|\varphi\|$ and a $*$-homomorphism
		\[\pi_\varphi\colon A\to M(C)\cap \{h\}'\]
		such that for $a\in A$ we have
		\[\pi_\varphi(a)h=\varphi(a).\]
		If $A$ is unital, then $h=\varphi(1_A)\in C$.
	\end{theorem}
	
	The following are some corollaries of Theorem \ref{OZAH}, giving important properties of order zero maps. The most important is Corollary \ref{OZCP}, which shows that order zero maps preserve Cuntz subequivalence. This provides crucial support for proving the transfer of many regularity properties later.
	\begin{corollary}\cite[Corollary 4.1]{WZ2009}
		Let $A$ and $B$ be C*-algebras and let $\varphi\colon A\to B$ be a completely positive contractive order zero map. Then, the map defined by $\rho_\varphi(\mathrm{id}_{(0,1]}\otimes a)\colon =\varphi(a)$ induces a $*$-homomorphism $\rho_\varphi\colon C_0((0,1])\otimes A\to B$.
		
		Conversely, any $*$-homomorphism $\rho\colon C_0((0,1])\otimes A\to B$ induces a completely positive contractive order zero map $\varphi_\rho\colon A\to B$ via $\varphi_\rho(a)\colon =\rho(\mathrm{id}_{(0,1]}\otimes a)$.
		
		This gives a canonical bijection between the space of completely positive contractive order zero maps from $A$ to $B$ and the space of homomorphisms from $C_0((0,1])\otimes A$ to $B$.
	\end{corollary}
	
	As in \cite{Win2010}, Theorem \ref{OZAH} allows us to perform the functional calculus for positive functions on completely positive contractive order zero maps.
	\begin{corollary}\cite[Corollary 4.2]{WZ2009}\label{OZF}
		Let $\varphi\colon A\to B$ be a completely positive contractive order zero map, and let $f\in C_0((0,1])$ be a positive function. Let $C$, $h$, and $\pi_\varphi$ be as in Theorem \ref{OZAH}. Then, the map 	\[f(\varphi)\colon A\to C\subset B\] 
		defined by
		\[f(\varphi)(a)\colon =f(h)\pi_\varphi(a)\quad \forall a\in A,\]
		is a well-defined completely positive order zero map. If $\|f\|\leq 1$, then $f(\varphi)$ is also contractive.
	\end{corollary}
	
	\begin{corollary}\cite[Corollary 4.3]{WZ2009}
		Let $A, B, C$ and $D$ be C*-algebras and let $\varphi\colon A\to B$ and $\psi\colon C\to D$ be completely positive contractive order zero maps. Then the induced completely positive contractive map
		\[\varphi\otimes_\mu \psi\colon A\otimes_\mu C\to B\otimes_\mu D\]
		has order zero, where $\otimes_\mu$ denotes the minimal or maximal tensor product. In particular, for any $k\in\mathbb{N}$,
		\[\varphi^{(k)}\colon M_k(A)\to M_k(B)\]
		has order zero.
	\end{corollary}
	
	\begin{corollary}\cite[Corollary 4.4]{WZ2009}
		Let $A$ and $B$ be C*-algebras, $\varphi\colon A\to B$ a completely positive order zero map, and $\tau$ a trace on $B$. Then, the composition $\tau\circ\varphi$ is a trace.
		
		The same conclusion holds when replacing ``trace'' with ``2-quasitrace''.
	\end{corollary}
	
	\begin{corollary}\cite[Corollary 4.5]{WZ2009}\label{OZCP}
		Let $A$ and $B$ be C*-algebras and let $\varphi\colon A\to B$ be a completely positive order zero map.
		
		Then, $\varphi$ induces a morphism of ordered semigroups
		\[\mathrm{W}(\varphi)\colon \mathrm{W}(A)\to \mathrm{W}(B)\]
		between the Cuntz semigroups via
		\[\mathrm{W}(\varphi)(\langle a\rangle)=\langle \varphi^{(k)}(a)\rangle\quad \forall a\in M_k(A)_+.\]
	\end{corollary}
	
	\begin{remark}
		The same conclusion holds when replacing $\mathrm{W}(A)$ with $\mathrm{Cu}(A)$ in the above corollary.
	\end{remark}
	
	\subsection{Maps under group actions}
	\begin{definition}
		Let $G$ be a locally compact group. We denote by $\mathrm{Lt}\colon G\to \mathrm{Aut}(C_0(G))$ the action induced by the left translation of $G$ on itself.
	\end{definition}
	
	\begin{definition}\cite[Definition 1.3]{HP2015}
		Let $G$ be a compact group, and let $A$, $B$ be C*-algebras. Let $\alpha\colon G\to \mathrm{Aut}(A)$ and $\gamma\colon G\to \mathrm{Aut}(B)$ be actions of $G$ on $A$ and $B$, respectively. Let $F\subseteq A$, $S\subseteq B$, and $\varepsilon>0$. A completely positive contractive map $\varphi\colon A\to B$ is said to be an $(F, S, \varepsilon)$-approximately equivariant central multiplicative map if:
		\begin{enumerate}[{(}1{)}]
			\item $\|\varphi(xy)-\varphi(x)\varphi(y)\|<\varepsilon$ for all $x, y\in F$.
			\item $\|\varphi(x)a-a\varphi(x)\|<\varepsilon$ for all $x\in F$ and all $a\in S$.
			\item $\sup_{g\in G}\|\varphi(\alpha_g(x))-\gamma_g(\varphi(x))\|<\varepsilon$ for all $x\in F$.\label{d2.15.3}
		\end{enumerate}
	\end{definition}
	
	By ignoring the group action and omitting condition (\ref{d2.15.3}), we obtain the definition of a general $(F, S, \varepsilon)$-approximately central multiplicative map.
	
	The following version of approximate centrality and approximate multiplicativity is convenient to use.
	\begin{definition}\cite[Definition 1.4]{MP2021}
		Let $A$ and $B$ be C*-algebras, and let $F\subseteq A$. A completely positive contractive map $\varphi\colon A\to B$ is said to be an $(n, F, \varepsilon)$-approximately multiplicative map if for any $m\in\{1, 2, \ldots, n\}$ and any $x_1, x_2, \ldots, x_m\in F$, we always have
		\[\|\varphi(x_1x_2\cdots x_m)-\varphi(x_1)\varphi(x_2)\cdots\varphi(x_m)\|<\varepsilon.\]
		If $S\subseteq B$ is also given, then $\varphi$ is said to be an $(n, F, S, \varepsilon)$-approximately central multiplicative map if it also satisfies $\|\varphi(x)a-a\varphi(x)\|<\varepsilon$ for all $x\in F$ and all $a\in S$.
	\end{definition}
	
	Next, we introduce some concepts slightly weaker than approximate multiplicativity.
	\begin{definition}
		Let $A$ be a unital C*-algebra, $B$ a C*-algebra, and $\varphi\colon A\to B$ a completely positive map. Let $\varepsilon>0$.
		\begin{enumerate}[{(}1{)}]
			\item We say $\varphi$ is $\varepsilon$-order zero if for $a, b\in A$ with $a\perp b$, we have
			\[\|\varphi(a)\varphi(b)\|<\varepsilon.\]
			\item We say $\varphi$ is strongly $\varepsilon$-order zero if for $a, b\in A$, we have
			\[\|\varphi(1)\varphi(ab)-\varphi(a)\varphi(b)\|<\varepsilon.\]
		\end{enumerate}
	\end{definition}
	
	Both definitions above are generalizations of order zero maps. However, since orthogonality is not complete under limits, while the strongly $\varepsilon$-order zero map avoids the use of orthogonal elements. Therefore, in Definition \ref{WTRC} we adopt the strongly $\varepsilon$-order zero map instead of the $\varepsilon$-order zero map, in order to have completeness when subsequently using sequence algebras.
	
	Replacing approximate multiplicativity with strongly approximate order zero yields the following definition.
	\begin{definition}
		Let $G$ be a compact group, let $A$ be a unital C*-algebra, and $B$ a C*-algebra. Let $\alpha\colon G\to \mathrm{Aut}(A)$ and $\gamma\colon G\to \mathrm{Aut}(B)$ be actions of $G$ on $A$ and $B$, respectively. Let $F\subseteq A$, $S\subseteq B$, and $\varepsilon>0$. A completely positive contractive map $\varphi\colon A\to B$ is said to be an $(F, S, \varepsilon)$-approximately central equivariant strongly order zero map if:
		\begin{enumerate}[{(}1{)}]
			\item $\|\varphi(1)\varphi(xy)-\varphi(x)\varphi(y)\|<\varepsilon$ for all $x, y\in F$.
			\item $\|\varphi(x)a-a\varphi(x)\|<\varepsilon$ for all $x\in F$ and all $a\in S$.
			\item $\sup_{g\in G}\|\varphi(\alpha_g(x))-\gamma_g(\varphi(x))\|<\varepsilon$ for all $x\in F$.
		\end{enumerate}
	\end{definition}

	\section{The weak tracial Rokhlin property with comparison}
	\label{sec8.1}
	
	This section establishes a definitional framework for the weak tracial Rokhlin property with comparison for actions of second-countable compact groups. As theoretical groundwork, we first systematically outline the core features of two existing theoretical frameworks: the weak tracial Rokhlin property for finite group actions and the tracial Rokhlin property with comparison for actions of second-countable compact groups.
	
	It is particularly noteworthy that even when restricted to finite groups, the new definition proposed herein constitutes a strictly stronger version than the traditional one. To elucidate this theoretical distinction, in Section \ref{sec8.4}, by introducing several auxiliary definitions and deeply analyzing their interrelationships, we reveal the key fact that a simple direct translation of the weak tracial Rokhlin property for finite group actions would fail to yield the desired theoretical results.
	
	\begin{definition}\cite[Definition 3.2]{AGP2021}\label{WTRF}
		Let $G$ be a finite group, let $A$ be a unital simple infinite-dimensional C*-algebra, and let $\alpha\colon G\to \mathrm{Aut}(A)$ be an action. We say that $\alpha$ has the weak tracial Rokhlin property if for every finite set $F\subseteq A$, every $\varepsilon>0$, and every $x\in A_+$ with $\|x\|=1$, there exist mutually orthogonal positive contractions $(d_g)_{g\in G}\in A$ with $d=\sum_{g\in G}d_g$ such that the following hold:
		\begin{enumerate}[{(}1{)}]
			\item $\|d_ga-ad_g\|<\varepsilon$ for all $a\in F$ and all $g\in G$.
			\item $\|\alpha_g(d_h)-d_{gh}\|<\varepsilon$ for all $g,h\in G$.
			\item $1-d\precsim_A x$.
			\item $\|dxd\|>1-\varepsilon$.
		\end{enumerate}
	\end{definition}
	
	The following lemma states that $d$ can be taken in $A^\alpha$, but with slightly weakened conditions.
	\begin{lemma}\cite[Lemma 2.14]{Asa2024}\label{WTRFE}
		Let $G$ be a finite group, let $A$ be a unital simple infinite-dimensional C*-algebra, and let $\alpha\colon G\to \mathrm{Aut}(A)$ be an action. Then $\alpha$ has the weak tracial Rokhlin property if and only if for every finite set $F\subseteq A$, every $\varepsilon>0$, and every $x\in A_+$ with $\|x\|=1$, there exists a positive element $d\in A^\alpha$ with $\|d\|=1$ and positive contractions $(d_g)_{g\in G}\in A$ such that the following hold:
		\begin{enumerate}[{(}1{)}]
			\item $d=\sum_{g\in G}d_g$.
			\item $\|d_ga-ad_g\|<\varepsilon$ for all $a\in F$ and all $g\in G$.
			\item $\|d_gd_h\|<\varepsilon$ for all $g,h\in G$ with $g\neq h$.
			\item $\|\alpha_g(d_h)-d_{gh}\|<\varepsilon$ for all $g,h\in G$.
			\item $(1-d-\varepsilon)_+\precsim_A x$.
			\item $\|dxd\|>1-\varepsilon$.
		\end{enumerate}
	\end{lemma}
	
	\begin{definition}\cite[Definition 2.4]{MP2021}\label{TRC}
		Let $G$ be a second-countable compact group, let $A$ be a unital simple infinite-dimensional C*-algebra, and let $\alpha\colon G\to \mathrm{Aut}(A)$ be an action. We say that the action $\alpha$ has the tracial Rokhlin property with comparison if for every $\varepsilon>0$, every finite set $F\subseteq A$, every finite set $S\subseteq C(G)$, every $x\in A_+$ with $\|x\|=1$, and every $y\in(A^\alpha)_+\setminus\{0\}$, there exist a projection $p\in A^\alpha$ and a unital completely positive map $\psi\colon C(G)\to pAp$ such that the following hold:
		\begin{enumerate}[{(}1{)}]
			\item $\psi$ is an $(S,F,\varepsilon)$-approximately central equivariant multiplicative map.
			\item $1-p\precsim_A x$.
			\item $1-p\precsim_{A^\alpha} y$.
			\item $1-p\precsim_{A^\alpha} p$.
			\item $\|pxp\|>1-\varepsilon$.
		\end{enumerate}
	\end{definition}
	
	The approximate multiplicativity requirement in Definition \ref{TRC} implies that $A$ must contain sufficiently many projections. However, for C*-algebras without non-trivial projections (such as the Jiang-Su algebra $\mathcal{Z}$), this means no such group action can exist on them. Strongly order zero maps impose no such restriction and can also approximately preserve orthogonality. Synthesizing the above two definitions, we now give the definition of the weak tracial Rokhlin property with comparison for compact group actions.
	
	\begin{definition}\label{WTRC}
		Let $G$ be a second-countable compact group, let $A$ be a unital simple infinite-dimensional C*-algebra, and let $\alpha\colon G\to \mathrm{Aut}(A)$ be an action. We say that the action $\alpha$ has the weak tracial Rokhlin property with comparison if for every $\varepsilon>0$, every finite set $F\subseteq A$, every finite set $S\subseteq C(G)$, every $x\in A_+$ with $\|x\|=1$, and every $y\in(A^\alpha)_+\setminus\{0\}$, there exist a positive contraction $d\in A^\alpha$ and a contractive completely positive map $\psi\colon C(G)\to \overline{dAd}$ such that the following hold:
		\begin{enumerate}[{(}1{)}]
			\item $\psi$ is an $(S,F,\varepsilon)$-approximately central equivariant strongly order zero map.
			\item $\psi(1)=d$.
			\item $(1-d-\varepsilon)_+\precsim_A x$.\label{d8.3.3}
			\item $(1-d-\varepsilon)_+\precsim_{A^\alpha} y$.\label{d8.3.4}
			\item $(1-d-\varepsilon)_+\precsim_{A^\alpha} d$.\label{d8.3.5}
			\item $\|dxd\|>1-\varepsilon$.\label{d8.3.6}
		\end{enumerate}
	\end{definition}
	
	Definition \ref{WTRC} requires $d$ to be $\alpha$-invariant, which is not necessary in Definition \ref{WTRF}. However, due to Lemma \ref{WTRFE}, we can require this condition. Like Definition \ref{TRC}, Definition \ref{WTRC} includes two conditions that have no analogues in Definition \ref{WTRF}. Condition (\ref{d8.3.4}) is automatic in the finite group case (see Proposition \ref{CIFA}). Condition (\ref{d8.3.5}) is automatic if $\mathrm{rc}(A)<1$ (see Proposition \ref{RCL1C}) or if $A$ is purely infinite (see Proposition \ref{PIC}). Even for finite groups, we do not know if it is automatic. Therefore, the term ``weak tracial Rokhlin property'' is not used here.
	
	The main technical result of this chapter is Theorem \ref{MTF}. This theorem, for given tolerance $\varepsilon>0$ and given compact subsets of $A$ and $A^\alpha$, constructs an $\alpha$-invariant positive contraction $d$ (which is ``sufficiently large'' in the sense described by conditions (\ref{d8.3.3}), (\ref{d8.3.4}), (\ref{d8.3.5}) of Definition \ref{WTRC}), and simultaneously constructs a completely positive contractive approximate order zero map from $A$ into $\overline{dA^\alpha d}$.
	
	For ease of comparison, we reformulate Definition \ref{WTRC} for finite groups to resemble the statement of Lemma \ref{WTRFE}.
	
	\begin{lemma}
		Let $G$ be a finite group, let $A$ be a unital simple infinite-dimensional C*-algebra, and let $\alpha\colon G\to \mathrm{Aut}(A)$ be an action. Then $\alpha$ has the weak tracial Rokhlin property with comparison if and only if for every finite set $F\subseteq A$, every $\varepsilon>0$, every $x\in A_+$ with $\|x\|=1$, and every $y\in(A^\alpha)_+\setminus\{0\}$, there exist a positive contraction $d\in A^\alpha$ and positive contractions $(d_g)_{g\in G}\in A$ such that the following hold:
		\begin{enumerate}[{(}1{)}]
			\item $d=\sum_{g\in G}d_g$.
			\item $\|d_ga-ad_g\|<\varepsilon$ for all $a\in F$ and all $g\in G$.
			\item $\|d_gd_h\|<\varepsilon$ for all $g,h\in G$ with $g\neq h$.
			\item $\|\alpha_g(d_h)-d_{gh}\|<\varepsilon$ for all $g,h\in G$.
			\item $(1-d-\varepsilon)_+\precsim_A x$.
			\item $(1-d-\varepsilon)_+\precsim_{A^\alpha} y$.
			\item $(1-d-\varepsilon)_+\precsim_{A^\alpha} d$.
			\item $\|dxd\|>1-\varepsilon$.
		\end{enumerate}
		\begin{proof}
			Set $S=\{\chi_{\{g\}}\colon g\in G\}\subseteq C(G)$. It is straightforward to verify that Definition \ref{WTRC} is equivalent to the same statement.
			
			For sufficiency: Let $F\subseteq A$, $\varepsilon>0$, $x\in A_+$ with $\|x\|=1$, and $y\in(A^\alpha)_+\setminus\{0\}$. Let $d$ and $(d_g)_{g\in G}$ be as stated in the lemma. We can define a contractive completely positive map $\psi\colon C(G)\to \overline{dAd}$ by $\psi(f)=\sum_{g\in G}f(g)d_g$ for $f\in C(G)$. Clearly $\psi(1)=d$. Then we have
			\[\max_{g\in G}\|(\psi\circ{\textrm{Lt}}_g)(\chi_{\{h\}})-(\alpha_g\circ\psi)(\chi_{\{h\}})\|=\max_{g\in G}\|d_{gh}-\alpha_g(d_h)\|<\varepsilon.\]
			And for $g\in G$ and $a\in F$, we have
			\[\|\psi(\chi_{\{g\}})a-a\psi(\chi_{\{g\}})\|=\|d_ga-ad_g\|<\varepsilon.\]
			Similarly, for $g,h\in G$ with $g\neq h$, we have
			\[\|\psi(1)\psi(\chi_{\{g\}}\chi_{\{h\}})-\psi(\chi_{\{g\}})\psi(\chi_{\{h\}})\|=\|d_gd_h\|<\varepsilon.\]
			Conditions (\ref{d8.3.3}), (\ref{d8.3.4}), (\ref{d8.3.5}), and (\ref{d8.3.6}) from Definition \ref{WTRC} follow immediately.
			
			For necessity: Assume $\alpha$ has the weak tracial Rokhlin property with comparison. Let $F\subseteq A$ be finite, $\varepsilon>0$, $x\in A_+$ with $\|x\|=1$, and $y\in(A^\alpha)_+\setminus\{0\}$. Without loss of generality, we may assume $\|a\|\leq 1$ for all $a\in F$.
			
			Applying Definition \ref{WTRC} with the given $\varepsilon, F, S, x$, and $y$, we obtain a positive contraction $d\in A^\alpha$ and a contractive completely positive map $\psi\colon C(G)\to \overline{dAd}$ such that the following hold:
			\begin{enumerate}[{(}1{)}]
				\setcounter{enumi}{8}
				\item $\psi(1)=d$.
				\item $\|\psi(\chi_{\{g\}})\psi(\chi_{\{h\}})\|<\varepsilon$ for all $g,h\in G$ with $g\neq h$.
				\item $\|\psi(\chi_{\{g\}})a-a\psi(\chi_{\{g\}})\|<\varepsilon$ for all $g\in G$ and $a\in F$.
				\item $\|(\alpha_g\circ\psi)(\chi_{\{h\}})-\psi(\chi_{\{gh\}})\|<\varepsilon$ for all $g,h\in G$.
				\item $(1-d-\varepsilon)_+\precsim_A x$.
				\item $(1-d-\varepsilon)_+\precsim_{A^\alpha} y$.
				\item $(1-d-\varepsilon)_+\precsim_{A^\alpha} d$.
				\item $\|dxd\|>1-\varepsilon$.
			\end{enumerate}
			
			Set $d_g=\psi(\chi_{\{g\}})$. Since $\sum_{g\in G}d_g=d\in A^\alpha$, it is straightforward to verify that $d$ and $(d_g)_{g\in G}$ satisfy the conditions of the lemma.
		\end{proof}
	\end{lemma}
	
	Definition \ref{WTRC} remains valid under the following modification: if we assume $F$ and $S$ are compact sets instead of finite sets, and require $\psi\colon C(G)\to \overline{dAd}$ to be exactly equivariant.
	\begin{lemma}\label{WTRE}
		Let $G$ be a second-countable compact group, let $A$ be a unital simple infinite-dimensional C*-algebra, and let $\alpha\colon G\to \mathrm{Aut}(A)$ be an action. Then $\alpha$ has the weak tracial Rokhlin property with comparison if and only if for every $\varepsilon>0$, every compact set $F\subseteq A$, every compact set $S\subseteq C(G)$, every $x\in A_+$ with $\|x\|=1$, and every $y\in(A^\alpha)_+\setminus\{0\}$, there exist a positive contraction $d\in A^\alpha$ and a contractive completely positive equivariant map $\psi\colon C(G)\to \overline{dAd}$ such that the following hold:
		\begin{enumerate}[{(}1{)}]
			\item $\psi$ is an $(S,F,\varepsilon)$-approximately central strongly order zero map.
			\item $\psi(1)=d$.
			\item $(1-d-\varepsilon)_+\precsim_A x$.
			\item $(1-d-\varepsilon)_+\precsim_{A^\alpha} y$.
			\item $(1-d-\varepsilon)_+\precsim_{A^\alpha} d$.
			\item $\|dxd\|>1-\varepsilon$.\label{l8.3.6}
		\end{enumerate}
		\begin{proof}
			The proof is essentially the same as \cite[Lemma 2.9]{MP2021}, so we omit it.
		\end{proof}
	\end{lemma}
	
	The following lemma shows that the element $x$ in condition (\ref{l8.3.6}) of Lemma \ref{WTRE} can be replaced by another positive contraction. This will be used in subsequent proofs, so we state it here.
	\begin{lemma}\label{WTRD}
		Let $G$ be a second-countable compact group, let $A$ be a unital simple infinite-dimensional C*-algebra, and let $\alpha\colon G\to \mathrm{Aut}(A)$ be an action. Then $\alpha$ has the weak tracial Rokhlin property with comparison if and only if for every $\varepsilon>0$, every finite set $F\subseteq A$, every finite set $S\subseteq C(G)$, every $x\in A_+$, every $z\in A_+$ with $\|z\|=1$, and every $y\in(A^\alpha)_+\setminus\{0\}$, there exist a positive contraction $d\in A^\alpha$ and a contractive completely positive equivariant map $\psi\colon C(G)\to \overline{dAd}$ such that the following hold:
		\begin{enumerate}[{(}1{)}]
			\item $\psi$ is an $(S,F,\varepsilon)$-approximately central strongly order zero map.
			\item $\psi(1)=d$.
			\item $(1-d-\varepsilon)_+\precsim_A x$.
			\item $(1-d-\varepsilon)_+\precsim_{A^\alpha} y$.
			\item $(1-d-\varepsilon)_+\precsim_{A^\alpha} d$.
			\item $\|dzd\|>1-\varepsilon$.\label{l8.4.6}
		\end{enumerate}
		\begin{proof}
			The proof is essentially the same as \cite[Lemma 3.9]{FG2020}, so we omit it.
		\end{proof}
	\end{lemma}
	
	When $A$ is finite, condition (\ref{d8.3.6}) of Definition \ref{WTRC}, condition (\ref{l8.3.6}) of Lemma \ref{WTRE}, and condition (\ref{l8.4.6}) of Lemma \ref{WTRD} are redundant. To prove this, we need a technique from \cite[Lemma 4.5]{Phi2014}.
	
	\begin{lemma}\cite[Lemma 2.6]{AGP2021}\label{CBG}
		Let $A$ be a unital C*-algebra, let $a,g\in A_+$ with $0\leq a,g \leq 1$, and let $\varepsilon_1,\varepsilon_2\geq 0$. Then
		\[(a-(\varepsilon_1+\varepsilon_2))_+\precsim_A[(1-g)a(1-g)-\varepsilon_1]_+\oplus(g-\frac{\varepsilon_2}{2})_+.\]
	\end{lemma}
	
	\begin{lemma}\label{NL}
		Let $A$ be a unital simple infinite-dimensional finite C*-algebra. Let $x\in A_+$ with $\|x\|=1$. Then for every $\varepsilon>0$, there exists $y\in(\overline{xAx})_+\setminus\{0\}$ such that whenever $g\in A_+$ satisfies $0\leq g\leq 1$, we have $\|(1-g)x(1-g)\|>1-\varepsilon$.
		\begin{proof}
			The proof is essentially the same as \cite[Lemma 2.9]{Phi2014}; we simply replace \cite[Lemma 1.8]{Phi2014} with Lemma \ref{CBG}.
		\end{proof}
	\end{lemma}
	
	\begin{lemma}
		Let $G$ be a second-countable compact group, let $A$ be a unital simple infinite-dimensional C*-algebra, and let $\alpha\colon G\to \mathrm{Aut}(A)$ be an action. Then the action $\alpha$ has the weak tracial Rokhlin property with comparison if and only if for every $\varepsilon>0$, every finite set $F\subseteq A$, every finite set $S\subseteq C(G)$, every $x\in A_+$ with $\|x\|=1$, and every $y\in(A^\alpha)_+\setminus\{0\}$, there exist a positive contraction $d\in A^\alpha$ and a contractive completely positive map $\psi\colon C(G)\to \overline{dAd}$ such that the following hold:
		\begin{enumerate}[{(}1{)}]
			\item $\psi$ is an $(S,F,\varepsilon)$-approximately central equivariant strongly order zero map.
			\item $\psi(1)=d$.
			\item $(1-d-\varepsilon)_+\precsim_A x$.
			\item $(1-d-\varepsilon)_+\precsim_{A^\alpha} y$.
			\item $(1-d-\varepsilon)_+\precsim_{A^\alpha} d$.
		\end{enumerate}
		\begin{proof}
			The proof is essentially the same as \cite[Lemma 4.5]{Phi2014}; we simply replace \cite[Lemma 2.9]{Phi2014} with Lemma \ref{NL}.
		\end{proof}
	\end{lemma}
	
	When $A$ is separable, we characterize the weak tracial Rokhlin property with comparison in terms of central sequences. We begin with a definition.
	
	\begin{definition}
		Let $A$ be a C*-algebra and let $\alpha\colon G\to\mathrm{Aut}(A)$ be an action of a topological group. Let $f,g\in (A_\infty,\alpha)_+$. We say $f$ is $\alpha$-pointwise Cuntz subequivalent to $g$ and write $f\precsim_{\alpha-p.w.-A} g$, if there exist a representing sequence $(f_n)_{n\in\mathbb{N}}\in l_\alpha^\infty(\mathbb{N},A)$ for $f$ and a representing sequence $(g_n)_{n\in\mathbb{N}}\in l_\alpha^\infty(\mathbb{N},A)$ for $g$ such that each $f_n, g_n$ is positive and for every $\varepsilon>0$ and all sufficiently large $n\in\mathbb{N}$, we have $(f_n-\varepsilon)_+\precsim_A g_n$. Without the group action, we say $f$ is pointwise Cuntz subequivalent to $g$ and write $f\precsim_{p.w.-A} g$.
	\end{definition}
	
	Next, we show that it suffices to consider a single lift instead of all.
	\begin{lemma}
		Let $f,g\in (A_\infty,\alpha)_+$. Then $f$ is $\alpha$-pointwise Cuntz subequivalent to $g$ if and only if there exist a representing sequence $(f_n)_{n\in\mathbb{N}}\in l_\alpha^\infty(\mathbb{N},A)$ for $f$ and a representing sequence $(g_n)_{n\in\mathbb{N}}\in l_\alpha^\infty(\mathbb{N},A)$ for $g$ such that each $f_n, g_n$ is positive and for every $\varepsilon>0$ and all sufficiently large $n\in\mathbb{N}$, we have $(f_n-\varepsilon)_+\precsim_A g_n$.
		\begin{proof}
			We only prove sufficiency. Suppose there exist a representing sequence $(f_n)_{n\in\mathbb{N}}\in l_\alpha^\infty(\mathbb{N},A)$ for $f$ and a representing sequence $(g_n)_{n\in\mathbb{N}}\in l_\alpha^\infty(\mathbb{N},A)$ for $g$ such that each $f_n, g_n$ is positive and for every $\varepsilon>0$ and all sufficiently large $n\in\mathbb{N}$, we have $(f_n-\frac{\varepsilon}{4})_+\precsim_A g_n$. If $(c_n)_{n\in\mathbb{N}}\in l_\alpha^\infty(\mathbb{N},A)$ is a representing sequence for $f$, and $(d_n)_{n\in\mathbb{N}}\in l_\alpha^\infty(\mathbb{N},A)$ is a representing sequence for $g$ such that each $c_n, d_n$ is positive, we need to show that for all sufficiently large $n\in\mathbb{N}$, we have
			\[(c_n-\varepsilon)_+\precsim_A d_n.\]
			Since $(f_n-\frac{\varepsilon}{4})_+\precsim_A g_n$ for all sufficiently large $n\in\mathbb{N}$, by part (1) of Lemma \ref{CCP}, there exists $\delta>0$ such that for all sufficiently large $n\in\mathbb{N}$, we have
			\[(f_n-\frac{\varepsilon}{2})_+\precsim_A (g_n-\delta)_+.\]
			We have $\lim_{n\to\infty}\|c_n-f_n\|=0$ and $\lim_{n\to\infty}\|d_n-g_n\|=0$, hence $\lim_{n\to\infty}\|c_n-(f_n-\frac{\varepsilon}{2})_+\|< \varepsilon$. Therefore, by part (4) of Lemma \ref{CCP}, for all sufficiently large $n\in\mathbb{N}$, we have
			\[(c_n-\varepsilon)_+\precsim_A (f_n-\frac{\varepsilon}{2})_+\ \mbox{and}\ (g_n-\delta)_+\precsim_A d_n.\]
			Thus for all sufficiently large $n\in\mathbb{N}$, we have
			\[(c_n-\varepsilon)_+\precsim_A(f_n-\frac{\varepsilon}{2})_+\precsim_A(g_n-\delta)_+\precsim_A d_n.\]
		\end{proof}
	\end{lemma}
	
	\begin{lemma}\label{WTRCS}
		Let $G$ be a second-countable compact group, let $A$ be a unital simple separable infinite-dimensional C*-algebra, and let $\alpha\colon G\to \mathrm{Aut}(A)$ be an action. Then $\alpha$ has the weak tracial Rokhlin property with comparison if and only if for every $x\in A_+$ with $\|x\|=1$ and every $y\in(A^\alpha)_+\setminus\{0\}$, there exist a positive contraction $d\in (A_{\infty,\alpha}\cap A')^{\alpha_\infty}$ and an equivariant completely positive contractive order zero map
		\[\varphi\colon C(G)\to\overline{d(A_{\infty,\alpha}\cap A')d}\]
		such that the following hold:
		\begin{enumerate}[{(}1{)}]
			\item $\varphi(1)=d$.\label{l8.9.1}
			\item $1-d\precsim_{\alpha-p.w.-A} x$.\label{l8.9.2}
			\item $1-d\precsim_{p.w.-A^\alpha} y$.\label{l8.9.3}
			\item $1-d\precsim_{p.w.-A^\alpha} d$.\label{l8.9.4}
			\item $\|dxd\|=1$.\label{l8.9.5}
		\end{enumerate}
		\begin{proof}
			Proof of sufficiency. Let $\varepsilon>0$, finite set $F\subseteq A$, finite set $S\subseteq C(G)$, $x\in A_+$ with $\|x\|=1$, and $y\in(A^\alpha)_+\setminus\{0\}$. Choose the given $d$ and $\varphi$. Choose a positive contraction $h=(h_n)_{n\in\mathbb{N}}\in l^\infty_\alpha(\mathbb{N},A)$ such that $\pi_A(h)=d$. By averaging over $G$, we may assume $h$ is $\alpha^\infty$-invariant. Then $\pi_A$ defines a surjection from $\overline{h l^\infty_\alpha(\mathbb{N},A) h}$ to $\overline{d A_{\infty,\alpha}d}$. By the Choi-Effros lifting Theorem, we can choose a lift $\theta=(\theta_n)_{n\in\mathbb{N}}\colon C(G)\to \overline{h l^\infty_\alpha(\mathbb{N},A) h}$ for $\varphi$ consisting of completely positive contractive maps $\theta_n\colon C(G)\to \overline{h_n A h_n}$. We claim the following properties hold:
			\begin{enumerate}[{(}1{)}]
				\setcounter{enumi}{5}
				\item $\theta_n(1)=h_n$ for all $n\in\mathbb{N}$.\label{l8.9.6}
				\item $\lim\limits_{n\to\infty}\|\theta_n(1)\theta_n(f_1f_2)-\theta_n(f_1)\theta_n(f_2)\|=0$ for all $f_1,f_2\in C(G)$.\label{l8.9.7}
				\item $\lim\limits_{n\to\infty}\|\theta_n(f)a-a\theta_n(f)\|=0$ for all $f\in C(G)$ and $a\in F$.\label{l8.9.8}
				\item $\lim\limits_{n\to\infty}\sup_{g\in G}\|(\theta_n\circ \mathrm{Lt}_g)(f)-(\alpha_g\circ\theta_n)(f)\|=0$ for all $f\in C(G)$.\label{l8.9.9}
			\end{enumerate}
			Properties (\ref{l8.9.6}), (\ref{l8.9.7}), (\ref{l8.9.8}) follow from the facts that $\pi_A\circ\theta=\varphi$ is a completely positive contractive order zero map, $\varphi(1)=d$, and its range is contained in $A_{\infty,\alpha}\cap A'$. For (\ref{l8.9.9}): Let $f\in C(G)$ and define $\rho_n\colon G\to[0,\infty)$ by
			\[\rho_n(g)=\sup_{m\geq n}\|(\theta_m\circ \mathrm{Lt}_g)(f)-(\alpha_g\circ\theta_m)(f)\|\]
			for $g\in G$. Since $\pi_A\circ\theta=\varphi$ is equivariant, we have $\lim_{n\to\infty}\rho_n(g)=0$ for all $g\in G$, and $\rho_m(g)\leq\rho_n(g)$ whenever $g\in G$ and $m\geq n$. Let $\sigma\colon l_\alpha^\infty(\mathbb{N},A)\to l_\alpha^\infty(\mathbb{N},A)$ be the map defined by $\sigma((a_n)_{n\in\mathbb{N}})=(a_2,a_3,\dots)$. Now
			\[\rho_n(g)=\|\sigma^{n-1}((\theta\circ \mathrm{Lt}_g)(f)-(\alpha_g^\infty\circ\theta)(f))\|,\]
			so $\rho_n$ is continuous. Therefore, (\ref{l8.9.9}) follows by Dini's theorem.
			
			Let $\mu$ be the normalized Haar measure on $G$. For $n\in\mathbb{N}$ and $f\in C(G)$, define
			\[\gamma_n(f)=\int_G(\alpha_g^{-1}\circ\theta_n\circ \mathrm{Lt}_g)(f)d\mu(g).\]
			Then $\gamma_n$ is an equivariant completely positive contractive map from $C(G)$ to $A$ with $\gamma_n(1)=h_n$, and by (\ref{l8.9.9}), we have
			\begin{enumerate}[{(}1{)}]
				\setcounter{enumi}{9}
				\item $\lim_{n\to\infty}\|\gamma_n(f)-\theta_n(f)\|=0$ for all $f\in C(G)$.\label{l8.9.10}
			\end{enumerate}
			By (\ref{l8.9.7}), (\ref{l8.9.8}), and (\ref{l8.9.10}), we can choose $n_0\in\mathbb{N}$ such that for all $n\geq n_0$, we have $\|\gamma_n(1)\gamma_n(f_1f_2)-\gamma_n(f_1)\gamma_n(f_2)\|<\varepsilon$ for all $f_1,f_2\in S$ and $\|\gamma_n(f)a-a\gamma_n(f)\|<\varepsilon$ for all $a\in F$ and $f\in S$.
			
			Since $1-d\precsim_{\alpha-p.w.-A} x$, we can choose $n_1\in\mathbb{N}$ such that for all $n\geq n_1$, we have $(1-h_n-\varepsilon)_+\precsim_Ax$. Similarly, since $1-d\precsim_{p.w.-A^\alpha} y$, we can choose $n_2\in\mathbb{N}$ such that for all $n\geq n_2$, we have $(1-h_n-\varepsilon)_+\precsim_{A^\alpha}y$. Similarly, since $1-d\precsim_{p.w.-A^\alpha} d$, we can choose $n_3\in\mathbb{N}$ such that for all $n\geq n_3$, we have $(1-h_n-\varepsilon)_+\precsim_{A^\alpha}h_n$. Also, we can choose $n_4\in\mathbb{N}$ such that for all $n\geq n_4$, we have $\|h_nxh_n\|>1-\varepsilon$. Set $N=\max(n_0,n_1,n_2,n_3,n_4)$, choose $\psi=\gamma_N$ and $d=h_N$, then all conditions are satisfied.
			
			Now assume the weak tracial Rokhlin property with comparison holds. Let $x\in A_+$ with $\|x\|=1$ and $y\in(A^\alpha)_+\setminus\{0\}$. We may assume $\|y\|=1$. Choose dense sequences
			\[f_1,f_2,\dots\in C(G),\quad a_1,a_2,\dots\in A.\]
			Since the weak tracial Rokhlin property with comparison holds, by Lemma \ref{WTRE}, for $S_n=\{f_1,f_2,\dots,f_n\}$, $F_n=\{a_1,a_2,\dots,a_n\}$, $\frac{1}{n}$, $x$, $y$, we obtain positive contractions $h_n\in A^\alpha$ and equivariant completely positive contractive maps $\psi_n\colon C(G)\to \overline{h_nAh_n}$ such that the following hold:
			\begin{enumerate}[{(}1{)}]
				\setcounter{enumi}{10}
				\item $\psi_n$ is an ($S_n$,$F_n$,$\frac{1}{n}$)-approximately central strongly order zero map.
				\item $\psi_n(1)=h_n$.
				\item $(1-h_n-\frac{1}{n})_+\precsim_A x$.
				\item $(1-h_n-\frac{1}{n})_+\precsim_{A^\alpha} y$.
				\item $(1-h_n-\frac{1}{n})_+\precsim_{A^\alpha} h_n$.
				\item $\|h_nxh_n\|>1-\frac{1}{n}$. 
			\end{enumerate}
			
			Define an $\alpha^\infty$-invariant positive contraction $h\in l^\infty_\alpha(\mathbb{N},A)$ by $h=(h_n)_{n\in\mathbb{N}}$. Define an equivariant completely positive contractive map $\psi\colon C(G)\to l^\infty(\mathbb{N},A)$ by $\psi(f)=(\psi_n(f))_{n\in\mathbb{N}}$ for $f\in C(G)$. Clearly $\psi(1)=h$. By equivariance, the range of $\psi$ is contained in $l^\infty_\alpha(\mathbb{N},A)$. Set $d=\pi_A(h)$ and $\varphi=\pi_A\circ\psi$. It is easy to see that $\varphi$ is an equivariant completely positive contractive order zero map from $C(G)$ to $\overline{d(A_{\infty,\alpha}\cap A')d}$ with $\varphi(1)=d$.
			
			Conditions (\ref{l8.9.2}), (\ref{l8.9.3}), and (\ref{l8.9.4}) follow by taking $n\in\mathbb{N}$ sufficiently large such that $\frac{1}{n}<\varepsilon$.
			
			For (\ref{l8.9.5}): Since $1>\|h_nxh_n\|>1-\frac{1}{n}$, we have $\lim_{n\to\infty}\|h_nxh_n\|=1$. Therefore $\|dxd\|=1$.
		\end{proof}
	\end{lemma}
	
	We now introduce an ``averaging process'' introduced by Gardella in \cite{Gar2014}, which allows us to transfer certain properties from the original algebra to the fixed point algebra.
	
	Let $\alpha\colon G\to \mathrm{Aut}(A)$ be a continuous action of a compact group on a unital C*-algebra. Identify $C(G)\otimes A$ with $C(G,A)$ (sending $f\otimes a$ to the function $h\mapsto f(h)a$), and let $\gamma\colon G\to \mathrm{Aut}(C(G,A))$ denote the diagonal action, i.e., $\gamma_g(a)(h)=\alpha_g(a(g^{-1}h))$ for all $g,h\in G$ and $a\in C(G,A)$. Define an averaging process $\varphi\colon C(G,A)\to C(G,A)$ by, for all $a\in C(G,A)$ and all $g\in G$,
	\[\varphi(a)(g)=\alpha_g(a(1)).\]
	It is easy to verify that $\gamma_g(\varphi(a))=\varphi(a)$ for all $g\in G$ and all $a\in C(G,A)$.
	\begin{lemma}\cite[Lemma 4.2]{Gar2014}\label{POU1}
		Let $A$ be a unital C*-algebra, let $G$ be a compact group, and let $\alpha\colon G\to \mathrm{Aut}(A)$ be a continuous action. Let $\varphi\colon C(G,A)\to C(G,A)$ be the averaging process defined above, and let $\gamma\colon G\to \mathrm{Aut}(C(G,A))$ be the diagonal action. Given $\varepsilon>0$ and a compact set $F\subseteq C(G,A)$, there exist $\delta>0$, a finite set $K\subseteq G$, and continuous functions $f_k\in C(G)$ for $k\in K$ such that the following hold:
		\begin{enumerate}[{(}1{)}]
			\item If $g,h\in G$ satisfy $d(g,h)<\delta$, then $\|\gamma_g(a)-\gamma_h(a)\|<\varepsilon$ for all $a\in \cup_{g\in G}\gamma_g(F)$.\label{l8.10.1}
			\item $0\leq f_k\leq 1$ for all $k\in K$.\label{l8.10.2}
			\item The family $(f_k)_{k\in K}$ is a partition of unity on $G$.\label{l8.10.3}
			\item For $k_1, k_2\in K$, if $f_{k_1}f_{k_2}\neq 0$, then $d(k_1,k_2)<\delta$.\label{l8.10.4}
			\item For each $g\in G$ and each $a\in \cup_{g\in G}\gamma_g(F)$, we have
			\[\|\varphi(a)(g)-\sum_{k\in K} f_k(g)\alpha_k(a(1))\|<\varepsilon.\label{l8.10.5}\]
		\end{enumerate}
	\end{lemma}
	
	The following proof references \cite[Proposition 4.3]{Gar2014}.
	\begin{proposition}\label{POU}
		Let $G$ be a second-countable compact group, let $A$ be a unital simple infinite-dimensional C*-algebra, and let $\alpha\colon G\to \mathrm{Aut}(A)$ be an action with the weak tracial Rokhlin property with comparison. For every $\varepsilon>0$, every compact set $F\subseteq A$, every $x\in A_+$ with $\|x\|=1$, and every $y\in(A^\alpha)_+\setminus\{0\}$, there exist $\delta>0$, a finite set $K\subseteq G$, a partition of unity $(f_k)_{k\in K}\subseteq C(G)$, a positive contraction $d\in A^\alpha$, and a completely positive contractive map $\psi\colon C(G)\to\overline{dAd}$ such that the following hold:
		\begin{enumerate}[{(}1{)}]
			\item $\psi(1)=d$.\label{p8.1.1}
			\item If $g,g'\in G$ satisfy $d(g,g')<\delta$, then $\|\alpha_g(a)-\alpha_{g'}(a)\|<\frac{\varepsilon}{|K|^2}$ for all $a\in F$.\label{p8.1.2}
			\item $0\leq f_k\leq 1$ for all $k\in K$.\label{p8.1.3}
			\item If $k,k'\in K$ satisfy $f_kf_{k'}\neq 0$, then $d(k,k')<\delta$.\label{p8.1.4}
			\item For all $g\in G$ and all $a\in F$, we have
			\[\|\alpha_g(\sum_{k\in K}\psi(f_k)^\frac{1}{2}\alpha_k(a)\psi(f_k)^\frac{1}{2})-\sum_{k\in K}\psi(f_k)^\frac{1}{2}\alpha_k(a)\psi(f_k)^\frac{1}{2}\|<\varepsilon.\label{p8.1.5}\]
			\item For all $a\in F$ and all $k\in K$, we have
			\[\|a\psi(f_k)-\psi(f_k)a\|<\frac{\varepsilon}{|K|^2}\ \mbox{and}\ \|a\psi(f_k)^\frac{1}{2}-\psi(f_k)^\frac{1}{2}a\|<\frac{\varepsilon}{|K|^2}.\label{p8.1.6}\]
			\item If $k,k'\in K$ satisfy $f_kf_{k'}=0$, then
			\[\|\psi(f_k)^\frac{1}{2}\psi(f_{k'})^\frac{1}{2}\|<\frac{\varepsilon}{|K|^2}.\label{p8.1.7}\]
			\item $(1-d-\varepsilon)_+\precsim_A x$.\label{p8.1.8}
			\item $(1-d-\varepsilon)_+\precsim_{A^\alpha} y$.\label{p8.1.9}
			\item $(1-d-\varepsilon)_+\precsim_{A^\alpha} d$.\label{p8.1.10}
			\item $\|dxd\|>1-\varepsilon$.\label{p8.1.11}
		\end{enumerate}
		\begin{proof}
			Without loss of generality, we can assume that $F$ is $\alpha$-invariant. For $F$ and $\frac{\varepsilon}{4}$, by Lemma \ref{POU1}, we can obtain $\delta>0$, a finite subset $K\subseteq G$, and a partition of unity $(f_k)_{k\in K}\subseteq C(G)$ such that conditions (\ref{l8.10.1}) - (\ref{l8.10.5}) in Lemma \ref{POU1} hold. Taking the continuous function $f(t)=t^\frac{1}{2}$, $t\in[0,1]$, and applying \cite[Lemma 1.5]{AP2020} to $\frac{\varepsilon}{|K|^2}$, we obtain $\delta_1>0$. Let $\delta_2=\min(\frac{\varepsilon}{|K|^2},\frac{\varepsilon^4}{|K|^8},\delta_1)$. Let $S=\{f_k\colon k\in K\}$, for the finite subset $S\subseteq C(G)$, the compact subset $F\subseteq A$, $\frac{1}{m}$, $x$ and $y$, by Lemma \ref{WTRE}, we can obtain $d_m\in A^\alpha$ and a completely positive contractive equivariant $(S,F,\frac{1}{m})$-approximately central strongly order zero map $\psi_m\colon C(G)\to \overline{d_mAd_m}$. Identify $C(G)\otimes A$ with $C(G,A)$, and for $m\in\mathbb{N}$, define a contractive linear map $\varphi_m\colon C(G,A)\to A$ by $\varphi_m(f\otimes a)=\psi_m(f)a$ for all $f\in C(G)$ and all $a\in A$. Clearly, $\varphi_m$ is equivariant, where we endow $C(G,A)$ with the diagonal action $\gamma$ of $G$. First, verify condition (\ref{p8.1.5}).
			
			Given $a\in F$ and given $g\in G$, in the first step using the approximate centrality of $\psi_m$, in the second step using the equivariance of $\varphi_m$, and in the last step applying condition (\ref{l8.10.5}) from the conclusion of Lemma \ref{POU1} and noting that for the averaging process $\gamma_g(\varphi(a))=\varphi(a)$, we have
			\begin{align*}
				&{}\limsup_{m\to\infty}\|\alpha_g(\sum_{k\in K}\psi_m(f_k)^\frac{1}{2}\alpha_k(a)\psi_m(f_k)^\frac{1}{2})-\sum_{k\in K}\psi_m(f_k)^\frac{1}{2}\alpha_k(a)\psi_m(f_k)^\frac{1}{2}\|\\
				=&{}\limsup_{m\to\infty}\|\alpha_g(\varphi_m(\sum_{k\in K}f_k\otimes\alpha_k(a)))-\varphi_m(\sum_{k\in K}f_k\otimes\alpha_k(a))\|\\
				=&{}\limsup_{m\to\infty}\|\varphi_m(\gamma_g(\sum_{k\in K}f_k\otimes\alpha_k(a))-\sum_{k\in K}f_k\otimes\alpha_k(a))\|\\
				\leq&{}\limsup_{m\to\infty}\|\gamma_g(\sum_{k\in K}f_k\otimes\alpha_k(a))-\sum_{k\in K}f_k\otimes\alpha_k(a)\|<\frac{\varepsilon}{2}.
			\end{align*}
			Choose a sufficiently large $m>\frac{1}{\delta_2}$ such that
			\[\|\alpha_g(\sum_{k\in K}\psi_m(f_k)^\frac{1}{2}\alpha_k(a)\psi_m(f_k)^\frac{1}{2})-\sum_{k\in K}\psi_m(f_k)^\frac{1}{2}\alpha_k(a)\psi_m(f_k)^\frac{1}{2}\|<\varepsilon.\]
			Then for such $m$, condition (\ref{p8.1.1}) is obvious. Due to the continuity of $\alpha$, condition (\ref{p8.1.2}) is also obvious. Conditions (\ref{p8.1.3}) and (\ref{p8.1.4}) come from Lemma \ref{POU1}. Conditions (\ref{p8.1.8}) - (\ref{p8.1.11}) come from Lemma \ref{WTRE}. 
			
			For condition (\ref{p8.1.6}), for all $a\in F$ and all $k\in K$, we have
			\[\|a\psi(f_k)-\psi(f_k)a\|<\delta_2=\min(\frac{\varepsilon}{|K|^2},\frac{\varepsilon^4}{|K|^8},\
			\delta_1).\]
			Therefore, by \cite[Lemma 1.5]{AP2020}, for all $a\in F$ and all $k\in K$, we have
			\[\|a\psi(f_k)^\frac{1}{2}-\psi(f_k)^\frac{1}{2}a\|<\frac{\varepsilon}{|K|^2}.\]
			For condition (\ref{p8.1.7}), whenever $k,k'\in K$ satisfy $f_kf_{k'}=0$, we have
			\[\|\psi(f_k)\psi(f_{k'})\|<\delta_2\leq\frac{\varepsilon}{|K|^2}.\]
			Hence,
			\begin{align*}
				\|\psi(f_k)^\frac{1}{2}\psi(f_{k'})^\frac{1}{2}\|^4=&{}\|\psi(f_k)^\frac{1}{2}\psi(f_{k'})\psi(f_k)\psi(f_{k'})\psi(f_k)^\frac{1}{2}\|\\
				\leq&{}\|\psi(f_k)\|\|\psi(f_{k'})\|\|\psi(f_k)\psi(f_{k'})\|\\
				\leq&{}\|\psi(f_k)\psi(f_{k'})\|<\delta_2.
			\end{align*}
			Thus $\|\psi(f_k)^\frac{1}{2}\psi(f_{k'})^\frac{1}{2}\|<\delta_2^\frac{1}{4}\leq \frac{\varepsilon}{|K|^2}$.
		\end{proof}
	\end{proposition}
	
	We now present the most important theorem of this section, which is a key tool for transferring properties from the original algebra to the fixed point algebra. The proof references \cite[Theorem 3.1]{Gar2019}.
	\begin{theorem}\label{MTF}
		Let $G$ be a second-countable compact group, let $A$ be a unital simple separable infinite-dimensional C*-algebra, and let $\alpha\colon G\to \mathrm{Aut}(A)$ be an action with the weak tracial Rokhlin property with comparison. For every $\varepsilon>0$, every compact set $F_1\subseteq A$, every compact set $F_2\subseteq A^\alpha$, every $x\in A_+$ with $\|x\|=1$, and every $y\in(A^\alpha)_+\setminus\{0\}$, there exist a positive contraction $d\in A^\alpha$ and a completely positive contractive map $\varphi\colon A\to \overline{dA^\alpha d}$ such that the following hold:
		\begin{enumerate}[{(}1{)}]
			\item $\varphi(1)=d$.\label{t8.1.1}
			\item $\|\varphi(a)\varphi(b)-d\varphi(ab)\|<\varepsilon$ for all $a,b\in F_1\cup F_2$.\label{t8.1.2}
			\item $\|da-ad\|<\varepsilon$ for all $a\in F_1\cup F_2$.\label{t8.1.3}
			\item $\|\varphi(a)-d^\frac{1}{2}ad^\frac{1}{2}\|<\varepsilon$ for all $a\in F_2$.\label{t8.1.4}
			\item $\|\varphi(a)\|\geq\|a\|-\varepsilon$ for all $a\in F_1\cup F_2$.\label{t8.1.5}
			\item $(1-d-\varepsilon)_+\precsim_A x$.\label{t8.1.6}
			\item $(1-d-\varepsilon)_+\precsim_{A^\alpha} y$.\label{t8.1.7}
			\item $(1-d-\varepsilon)_+\precsim_{A^\alpha} d$.\label{t8.1.8}
			\item $\|dxd\|>1-\varepsilon$.\label{t8.1.9}
		\end{enumerate}
		\begin{proof}
			Let $\varepsilon_0=\frac{\varepsilon}{7}$, take the continuous function $f(t)=t^\frac{1}{2}$, $t\in[0,1]$, and apply \cite[Lemma 1.5]{AP2020} to $\varepsilon_0$ to obtain $\delta_1>0$. Let $\delta_2=\min(\varepsilon_0,\
			\delta_1)$. Without loss of generality, we may assume that $\|a\|\leq1$ for all $a\in F_1\cup F_2$. For the compact set $F=F_1\cup F_2$ (without loss of generality, assume it is $\alpha$-invariant), $\delta_2$, $x\in A_+$ with $\|x\|=1$, and $y\in(A^\alpha)_+\setminus\{0\}$, using Proposition \ref{POU}, we can obtain a positive number $\delta>0$, a finite subset $K\subseteq G$, a partition of unity $(f_k)_{k\in K}\subseteq C(G)$, a positive contraction $d\in A^\alpha$, and a completely positive contractive map $\psi\colon C(G)\to\overline{dAd}$ such that the following hold:
			\begin{enumerate}[{(}1{)}]
				\setcounter{enumi}{9}
				\item $\psi(1)=d$.\label{t8.1.10}
				\item If $g,g'\in G$ and $d(g,g')<\delta$, then for all $a\in F$, $\|\alpha_g(a)-\alpha_{g'}(a)\|<\frac{\delta_2}{|K|^2}$.\label{t8.1.11}
				\item Whenever $k,k'\in K$ satisfy $f_kf_{k'}\neq 0$, then $d(k,k')<\delta$.\label{t8.1.12}
				\item For all $g\in G$ and all $a\in F$, we have
				\[\|\alpha_g(\sum_{k\in K}\psi(f_k)^\frac{1}{2}\alpha_k(a)\psi(f_k)^\frac{1}{2})-\sum_{k\in K}\psi(f_k)^\frac{1}{2}\alpha_k(a)\psi(f_k)^\frac{1}{2}\|<\delta_2.\label{t8.1.13}\]
				\item For all $a\in F$ and all $k\in K$, we have
				\[\|a\psi(f_k)-\psi(f_k)a\|<\frac{\delta_2}{|K|^2}\ \mbox{and} \ \|a\psi(f_k)^\frac{1}{2}-\psi(f_k)^\frac{1}{2}a\|<\frac{\delta_2}{|K|^2}.\label{t8.1.14}\]
				\item Whenever $k,k'\in K$ satisfy $f_kf_{k'}=0$, we have
				\[\|\psi(f_k)^\frac{1}{2}\psi(f_{k'})^\frac{1}{2}\|<\frac{\delta_2}{|K|^2}.\label{t8.1.15}\]
				\item $(1-d-\delta_2)_+\precsim_A x$.\label{t8.1.16}
				\item $(1-d-\delta_2)_+\precsim_{A^\alpha} y$.\label{t8.1.17}
				\item $(1-d-\delta_2)_+\precsim_{A^\alpha} d$.\label{t8.1.18}
				\item $\|dxd\|>1-\delta_2$.\label{t8.1.19}
			\end{enumerate}
				
			Denote the canonical conditional expectation by $E\colon A\to A^\alpha$. Define a completely positive contractive map $\varphi\colon A\to \overline{dA^\alpha d}$ by, for all $a\in F$,
			\[\varphi(a)=E(\sum_{k\in K}\psi(f_k)^\frac{1}{2}\alpha_k(a)\psi(f_k)^\frac{1}{2}).\]
			We will show that $\varphi$ satisfies all the properties. Clearly $\varphi(1)=d$, which is condition (\ref{t8.1.1}). Since $(1-d-\varepsilon)_+\leq(1-d-\delta_2)_+$, conditions (\ref{t8.1.6}),  (\ref{t8.1.7}),  (\ref{t8.1.8}), and (\ref{t8.1.9}) immediately hold.
				
			For condition (\ref{t8.1.2}), given $a,b\in F$, using (\ref{t8.1.13}) in the second step, (\ref{t8.1.11}),  (\ref{t8.1.12}),  (\ref{t8.1.14}),  (\ref{t8.1.15}) in the third step, (\ref{t8.1.10}) and the fact that $(f_k)_{k\in K}\subseteq C(G)$ is a partition of unity in the fourth step, and (\ref{t8.1.13}) in the fifth step, we have	
			\begin{align*}
				\varphi(a)\varphi(b)&{}=E(\sum_{k'\in K}\psi(f_{k'})^\frac{1}{2}\alpha_{k'}(a)\psi(f_{k'})^\frac{1}{2})E(\sum_{k\in K}\psi(f_k)^\frac{1}{2}\alpha_k(b)\psi(f_k)^\frac{1}{2})\\
				&{}\approx_{2\delta_2}\sum_{k'\in K}\sum_{k\in K}\psi(f_{k'})^\frac{1}{2}\alpha_{k'}(a)\psi(f_{k'})^\frac{1}{2}\psi(f_k)^\frac{1}{2}\alpha_k(b)\psi(f_k)^\frac{1}{2}\\
				&{}\approx_{4\delta_2}\sum_{k'\in K}\sum_{k\in K}\psi(f_{k'})\psi(f_k)^\frac{1}{2}\alpha_k(ab)\psi(f_k)^\frac{1}{2}\\
				&{}=d\sum_{k\in K}\psi(f_k)^\frac{1}{2}\alpha_k(ab)\psi(f_k)^\frac{1}{2}\\
				&{}\approx_{\delta_2}dE(\sum_{k\in K}\psi(f_k)^\frac{1}{2}\alpha_k(ab)\psi(f_k)^\frac{1}{2})\\
				&{}=d\varphi(ab).
			\end{align*}
			Hence $\|\varphi(a)\varphi(b)-d\varphi(ab)\|<7\delta_2\leq\varepsilon$, i.e., condition (\ref{t8.1.2}) holds.
				
			To prove condition (\ref{t8.1.3}), let $a\in F$. Using (\ref{t8.1.10}) and the fact that $(f_k)_{k\in K}\subseteq C(G)$ is a partition of unity in the first step, and (\ref{t8.1.14}) in the third step, we have
			\begin{align*}
				\|da-ad\|&{}=\|\sum_{k\in K}\psi(f_k)a-a\sum_{k\in K}\psi(f_k)\|\\
				&{}\leq\sum_{k\in K}\|\psi(f_k)a-a\psi(f_k)\|\\
				&{}<|K|\frac{\delta_2}{|K|^2}<\delta_2=\min(\varepsilon_0,\
				\delta_1)<\varepsilon.
			\end{align*}
			Then by \cite[Lemma 1.5]{AP2020}, we have
			\begin{equation}
				\|d^\frac{1}{2}a-ad^\frac{1}{2}\|<\varepsilon_0. \label{eq8.1.1} 
			\end{equation}
				
			To verify condition (\ref{t8.1.4}), let $a\in F_2\subseteq F$. Using (\ref{t8.1.14}) in the second step, (\ref{t8.1.10}) and the fact that $(f_k)_{k\in K}\subseteq C(G)$ is a partition of unity in the fourth step, and (\ref{eq8.1.1}) in the fifth step, we have
			\begin{align*}
				\varphi(a)&{}=E(\sum_{k\in K}\psi(f_k)^\frac{1}{2}\alpha_k(a)\psi(f_k)^\frac{1}{2})\\
				&{}\approx_{\delta_2}E(\sum_{k\in K}\psi(f_k)\alpha_k(a))\\
				&{}=E(\sum_{k\in K}\psi(f_k))a\\
				&{}=da\approx_{\varepsilon_0}d^\frac{1}{2}ad^\frac{1}{2}.
			\end{align*}
			Therefore $\|\varphi(a)-d^\frac{1}{2}ad^\frac{1}{2}\|<2\varepsilon_0<\varepsilon$.
				
			For condition (\ref{t8.1.5}), for any $n\in\mathbb{N}$, assume $\{S_n\}$ is a sequence of finite subsets in $A$ such that $\cup_{n\in\mathbb{N}}S_n$ is dense in $A$ and $F\subseteq S_n$ for all $n\in\mathbb{N}$. For $\frac{1}{n}$, $S_n$, $F_2$, $x$, and $y$, there exist a positive contraction $d_n\in A^\alpha$ and a completely positive contractive map $\varphi_n\colon A\to \overline{d_nA^\alpha d_n}$ such that the following hold:
			\begin{enumerate}[{(}1{)}]
				\setcounter{enumi}{19}
				\item $\varphi_n(1)=d_n$.
					
				\item For all $a,b\in S_n$, $\|\varphi_n(a)\psi_n(b)-d_n\varphi_n(ab)\|<\frac{1}{n}$.
					
				\item For all $a\in S_n$, $\|d_na-ad_n\|<\frac{1}{n}$.
					
				\item For all $a\in F_2$, $\|\varphi_n(a)-d_n^\frac{1}{2}ad_n^\frac{1}{2}\|<\frac{1}{n}$.
					
				\item $(1-d_n-\frac{1}{n})_+\precsim_Ax$.
					
				\item $(1-d_n-\frac{1}{n})_+\precsim_{A^\alpha} y$.
					
				\item $(1-d_n-\frac{1}{n})_+\precsim_{A^\alpha} d_n$.
					
				\item $\|d_nxd_n\|>1-\frac{1}{n}$.
			\end{enumerate}
				
			Let $\pi_{A^\alpha}\colon l^\infty(\mathbb{N},A^\alpha)\to (A^\alpha)_\infty$ be the quotient map. Denote $d=(d_n)_{n=1}^\infty\in l^\infty(\mathbb{N},A^\alpha)$ and $r=\pi_{A^\alpha}(d)$. Define a completely positive contractive order zero map $\Phi\colon A\to \overline{r(A^\alpha)_\infty r}$ by $\Phi(a)=\pi_{A^\alpha}(\{\varphi_1(a),\varphi_2(a),\dots\})$ for all $a\in A$. We have the following hold:
			\begin{enumerate}[{(}1{)}]
				\setcounter{enumi}{27}
				\item $\Phi(1)=r$.
					
				\item For all $a\in A$, $ar=ra$.
					
				\item For all $a\in F_2$, $\Phi(a)=r^\frac{1}{2}ar^\frac{1}{2}$.
					
				\item $1-r\precsim_{\alpha-p.w.-A}x$.
					
				\item $1-r\precsim_{p.w.-A^\alpha}y$.
				
				\item $1-r\precsim_{p.w.-A^\alpha}r$.
					
				\item $\|rxr\|=1$.\label{t8.1.34}
			\end{enumerate}
			From $\|rxr\|=1$ in (\ref{t8.1.34}), we can deduce $\|r\|=1$. Indeed, since $\|rxr\|=1$, we have $1=\|rxr\|\leq \|r\|\|x\|\|r\|= \|r\|^2\leq 1$, so $\|r\|=1$.
				
			Choose a finite set $T\subseteq F$ such that for every $a\in F$ there exists $b\in T$ with $\|b-a\|<\frac{\varepsilon}{3}$. From (\ref{t8.1.34}), we have $r\neq0$. Therefore, the order zero map $\Phi$ is nonzero. Since $A$ is simple, by \cite[Proposition 5.3]{FL2020}, $\|\Phi(a)\|=\|\Phi\|\|a\|$. Because $\|\Phi\|=\|\Phi(1)\|=\|r\|=1$, we conclude that this map is isometric. Hence, for sufficiently large $n$ and all $b\in T$, we have \[\|\varphi_n(b)\|>\|b\|-\frac{\varepsilon}{3}.\]
			Let $a\in F$ and choose $b\in T$ such that $\|b-a\|<\frac{\varepsilon}{3}$. Then for sufficiently large $n$, we have
			\begin{align*}
				\|\varphi_n(a)\|\geq&{}\|\varphi_n(b)\|-\|b-a\|\\
				>&{}\|b\|-\frac{\varepsilon}{3}-\frac{\varepsilon}{3}\\
				\geq&{}\|a\|-\|b-a\|-\frac{2\varepsilon}{3}\\
				>&{}\|a\|-\varepsilon.
			\end{align*}
			Now, for sufficiently large $n$, let $\varphi=\varphi_n$ and $d=d_n$, then conditions (\ref{t8.1.1}) - (\ref{t8.1.9}) hold.
		\end{proof}
	\end{theorem}
	
	\begin{remark}\label{MTFR}
		The naive weak tracial Rokhlin property (Definition \ref{NWTR}) is obtained from Definition \ref{WTRC} by omitting conditions (\ref{d8.3.4}) and (\ref{d8.3.5}). Similarly, the weak tracial Rokhlin property (Definition \ref{WTR}) is obtained from Definition \ref{WTRC} by omitting condition (\ref{d8.3.4}). We can obtain corresponding versions of Theorem \ref{MTF} by omitting the corresponding conditions.
	\end{remark}
	
	\section{Simplicity of the fixed point algebra and the crossed product}
	\label{sec8.2}
	In this section, we prove that for an action of a second-countable compact group with the weak tracial Rokhlin property with comparison, the fixed point algebra and the crossed product algebra formed by a unital simple separable infinite-dimensional C*-algebra remain simple.
	
	For the case of finite groups, it follows from \cite[Proposition 3.2]{FG2020} that the weak tracial Rokhlin property implies pointwise outerness. Therefore, by \cite[Theorem 3.1]{Kis1981}, the crossed product is simple. However, without the discreteness of the group, pointwise outerness of the action does not imply simplicity of the crossed product.
	
	The core concept of saturation for compact group actions allows us to prove the simplicity of the crossed product. First, we use Theorem \ref{MTF} to prove the simplicity of the fixed point algebra. Then we will prove that the weak tracial Rokhlin property with comparison implies saturation, and consequently, Theorem \ref{CI} gives the simplicity of the crossed product.
	
	\begin{theorem}\label{FAS}
		Let $A$ be a unital simple separable infinite-dimensional C*-algebra. Let $\alpha\colon G\to \mathrm{Aut}(A)$ be an action of a second-countable compact group with the weak tracial Rokhlin property with comparison. Then the fixed point algebra $A^\alpha$ is simple.
		\begin{proof}
			Let $I\subseteq A^\alpha$ be a nonzero ideal. We will show that $I$ contains an invertible element.
			
			Since $A$ is unital and simple, we have $AIA=A$. Therefore, there exist $n\in\mathbb{N}$, $a_1,a_2,\dots,a_n\in A$, and $x_1,x_2,\dots,x_n\in I$ such that
			$\sum_{i=1}^na_i^*x_ia_i=1$.
			
			Take a nonzero positive element $z\in I$. Let $F=\{x_i:i=1,2,\dots,n\}$. Suppose $\{S_n\}$ is a sequence of finite sets in $A$ such that $\cup_{n\in\mathbb{N}}S_n$ is dense in $A$ and $F\subseteq S_n$ for all $n\in\mathbb{N}$. For $\frac{1}{n}$, $z$, $S_n$, and $F$, there exist a positive contraction $d_n\in A^\alpha$ and a completely positive contractive map $\varphi_n\colon A\to \overline{d_nA^\alpha d_n}$ such that the following hold:
			\begin{enumerate}[{(}1{)}]
				\item $\varphi_n(1)=d_n$.
				\item $\|\varphi_n(a)\varphi_n(b)-d_n\varphi_n(ab)\|<\frac{1}{n}$ for all $a,b\in S_n$.
				\item $\|d_na-ad_n\|<\frac{1}{n}$ for all $a\in S_n$.
				\item $\|\varphi_n(a)-d_n^{\frac{1}{2}}ad_n^{\frac{1}{2}}\|<\frac{1}{n}$ for all $a\in F$.
				\item $(1-d_n-\frac{1}{n})_+\precsim_{A^\alpha} z$.
			\end{enumerate}
			Let $\pi_{A^\alpha}\colon l^\infty(\mathbb{N},A^\alpha)\to (A^\alpha)_\infty$ be the quotient map. Denote $d=(d_n)_{n=1}^\infty\in l^\infty(\mathbb{N},A^\alpha)$ and $r=\pi_{A^\alpha}(d)$. Define a completely positive contractive order zero map $\varphi\colon A\to \overline{r(A^\alpha)_\infty r}$ by $\varphi(a)=\pi_{A^\alpha}(\{\varphi_1(a),\varphi_2(a),\dots\})$ for all $a\in A$. We have the following hold:
			\begin{enumerate}[{(}1{)}]
				\setcounter{enumi}{5}
				\item $\varphi(1)=r$.
				\item $ar=ra$ for all $a\in A$.
				\item $\varphi(a)=r^\frac{1}{2}ar^\frac{1}{2}$ for all $a\in F$.
				\item $1-r\precsim_{p.w.-A^\alpha}z$.		 
			\end{enumerate}
			Therefore, by Theorem \ref{OZAH} and Corollary \ref{OZF}, we have
			\begin{align*}
				r=&{}\varphi(\sum_{i=1}^na_i^*x_ia_i)
				=r\pi_{\varphi}(\sum_{i=1}^na_i^*x_ia_i)\\
				=&{}\sum_{i=1}^n\pi_{\varphi}(a_i^*)r\pi_{\varphi}(x_i)\pi_{\varphi}(a_i)\\
				=&{}\sum_{i=1}^n\pi_{\varphi}(a_i^*)\varphi(x_i)\pi_{\varphi}(a_i)\\
				=&{}\sum_{i=1}^n\pi_{\varphi}(a_i^*)r^\frac{1}{2}x_ir^\frac{1}{2}\pi_{\varphi}(a_i)\\
				=&{}\sum_{i=1}^n\sqrt{\varphi}(a_i^*)x_i\sqrt{\varphi}(a_i).
			\end{align*}
			For $i=1,\dots,n$, we can lift $\sqrt{\varphi}(a_i)$ and $r$ to obtain $b_i\in A^\alpha$ and $d_{n_0}\in A^\alpha$ such that $(1-d_{n_0}-\frac{1}{2})_+\precsim_{A^\alpha} z$ and
			\[\|d_{n_0}-\sum_{i=1}^nb_i^*x_ib_i\|<\frac{1}{2}.\]
			Note that
			\begin{align*}
				&{}\|1-(1-d_{n_0}-\frac{1}{2})_+-\sum_{i=1}^nb_i^*x_ib_i\|\\
				=&{}\|1-(1-d_{n_0})+(1-d_{n_0})-(1-d_{n_0}-\frac{1}{2})_+-\sum_{i=1}^nb_i^*x_ib_i\|\\
				\leq&{}\|1-(1-d_{n_0})-\sum_{i=1}^nb_i^*x_ib_i\|+\|(1-d_{n_0})-(1-d_{n_0}-\frac{1}{2})_+\|\\
				<&{}\frac{1}{2}+\frac{1}{2}=1.
			\end{align*}
			Therefore, $(1-d_{n_0}-\frac{1}{2})_++\sum_{i=1}^nb_i^*x_ib_i$ is invertible and belongs to $I$.
		\end{proof}
	\end{theorem}
	
	\begin{theorem}
		Let $A$ be a unital simple separable infinite-dimensional C*-algebra. Let $\alpha\colon G\to \mathrm{Aut}(A)$ be an action of a second-countable compact group with the weak tracial Rokhlin property with comparison. Then the fixed point algebra $A^\alpha$ is infinite-dimensional and not of type I.
		\begin{proof}	
			We first prove infinite-dimensionality. We will show that for any $n\in \mathbb{N}$, there exist nonzero mutually orthogonal positive elements $x_1,\dots,x_n\in A^\alpha$. Since $A$ is simple and infinite-dimensional, there exist $y_1,\dots,y_n\in A_+$ with $\|y_k\|=1$ for $k=1,\dots,n$. For given $\frac{1}{3}$, $n$, choose $0<\delta<\frac{1}{3}$ satisfying the condition of \cite[Lemma 2.5.12]{Lin2001}.
			
			By Theorem \ref{MTF}, for $\delta<\frac{1}{3}$, $n=2$, $x=y=1$, $F_1=\{y_1,\dots,y_n\}$, and $F_2=\varnothing$. There exist a positive contraction $d\in A^\alpha$ and a completely positive contractive map $\varphi\colon A\to \overline{dA^\alpha d}$ such that the following hold:
			\begin{enumerate}[{(}1{)}]
				\item $\|d\varphi(y_jy_k)-\varphi(y_j)\varphi(y_k)\|<\delta$ for $j,k=1,2,\dots,n$.\label{t8.3.1}
				\item $\|\varphi(y_k)\|>1-\delta$ for $k=1,2,\dots,n$.
			\end{enumerate}
			From (\ref{t8.3.1}), we have $\|\varphi(y_j)\varphi(y_k)\|<\delta$ when $j\neq k$. By \cite[Lemma 2.5.12]{Lin2001}, there exist $x_1,\dots,x_n\in (A^\alpha)_+$ such that $\|x_k-\varphi(y_k)\|<\frac{1}{3}$ for $k=1,\dots,n$ and $x_jx_k=0$ for $j\neq k$. Moreover, for $k=1,\dots,n$,
			\[\|x_k\|>\|\varphi(y_k)\|-\frac{1}{3}>1-\delta-\frac{1}{3}\geq\frac{1}{3},\]
			hence $x_k\neq0$. Therefore, $A^\alpha$ is infinite-dimensional.
			
			Since $A^\alpha$ is unital, simple, and infinite-dimensional, it is not of type I.
		\end{proof}
	\end{theorem}
	
	We recall the definition of the strong Arveson spectrum for actions of compact groups on C*-algebras and related concepts. For more information, see \cite{GLP1994}.
	
	\begin{definition}\cite{GLP1994}\label{SAS}
		Let $A$ be a C*-algebra, let $\alpha\colon G\to \mathrm{Aut}(A)$ be an action of a compact group on $A$, and let $\pi\colon G\to \mathrm{U}(\mathcal{H}_\pi)$ be a unitary representation of $G$ on a Hilbert space. We define
		\[A_2(\pi)=\{x\in B(\mathcal{H}_\pi)\otimes A\colon (\mathrm{id}_{B(\mathcal{H}_\pi)}\otimes\alpha_g)(x)=x(\pi(g)\otimes1_A)\ \mbox{for all} g\in G\}.\]
		Let $\widehat{G}$ be a set containing exactly one representative from each unitary equivalence class of irreducible representations of $G$. Define the strong Arveson spectrum as
		\[\widetilde{\mathrm{Sp}}(\alpha)=\{\pi\in\widehat{G}\colon \overline{A_2(\pi)^*A_2(\pi)}=(B(\mathcal{H}_\pi)\otimes A)^{\mathrm{Ad}(\pi)\otimes\alpha}\}.\]
		Denote by $\mathcal{H}^\alpha(A)$ the family of all nonzero $\alpha$-invariant hereditary C*-subalgebras $B$ of $A$. Define the strong Connes spectrum as
		\[\widetilde{\Gamma}(\alpha)=\cap\{\widetilde{\mathrm{Sp}}(\alpha|B)\colon B\in\mathcal{H}^\alpha(A)\}.\]
	\end{definition}
	
	\begin{lemma}\cite[Lemma 3.5]{MP2021}\label{SS}
		In the setting of Definition \ref{SAS}, for any finite-dimensional unitary representation $\pi$ of $G$, the set $\overline{A_2(\pi)^*A_2(\pi)}$ is a closed two-sided ideal of $(B(\mathcal{H}_\pi)\otimes A)^{\mathrm{Ad}(\pi)\otimes\alpha}$.
	\end{lemma}
	
	Recall (\cite[Definition 7.1.4]{Phi2006}) that an action of a compact group $G$ on a C*-algebra $A$ is saturated if a suitable completion of $A$ is a Morita equivalence bimodule from $A^\alpha$ to $A\rtimes_\alpha G$. The following theorem gives us some essentially known results on saturation, the strong Arveson spectrum, and simplicity.
	
	\begin{theorem}\cite[Theorem 3.6]{MP2021}\label{CI}
		Let $A$ be a C*-algebra, let $\alpha\colon G\to \mathrm{Aut}(A)$ be an action of a compact group on $A$. The following are equivalent:
		\begin{enumerate}[{(}1{)}]
			\item $A^\alpha$ is simple and $\alpha$ is saturated.
			\item $A^\alpha$ is simple and $\widetilde{\mathrm{Sp}}(\alpha)=\widehat{G}$.
			\item $A\rtimes_\alpha G$ is simple.
			\item $A$ has no nontrivial $G$-invariant ideals and $\widetilde{\Gamma}(\alpha)=\widehat{G}$.
		\end{enumerate}
		When these conditions hold, $A^\alpha$ is isomorphic to a full hereditary subalgebra of $A\rtimes_\alpha G$ and is strongly Morita equivalent to $A\rtimes_\alpha G$.
	\end{theorem}
	
	Below we give an equivariant version of Lemma \ref{OZF}, that is, an equivariant completely positive order zero map remains an equivariant completely positive order zero map after continuous functional calculus by a nonnegative function.
	\begin{lemma}\label{OZFE}
		Let $G$ be a locally compact group, and let $A$ be a unital C*-algebra, $B$ a C*-algebra. Let $\alpha\colon G\to \mathrm{Aut}(A)$ and $\gamma\colon G\to \mathrm{Aut}(B)$ be actions of $G$ on $A$ and $B$, respectively. Let $\varphi\colon A\to B$ be a completely positive contractive order zero map, and let $f\in C_0((0,1])$ be a positive function. Let $C, h$, and $\pi_\varphi$ be as in Theorem \ref{OZAH}. Then, the map defined by
		\[f(\varphi)(a)\colon=f(h)\pi_\varphi(a)\quad \forall a\in A,\]
		\[f(\varphi)\colon A\to C\subset B\]
		is a well-defined equivariant completely positive order zero map. If the norm of $f$ is at most one, then $f(\varphi)$ is also contractive.
		\begin{proof}
			Here we only verify the equivariance. Since $\varphi$ is equivariant, for any $g\in G$ we have $\varphi(\alpha_g(1))=\beta_g(\varphi(1))$, i.e., $\varphi(1)=\beta_g(\varphi(1))$. By the arbitrariness of $g$, it follows that $\varphi(1)$ is $\beta$-invariant. Then for any $x\in A$, $g\in G$, by Lemma \ref{OZF} we have 
			\begin{align*}
				\pi_\varphi(\alpha_g(x))h=&{}\varphi(\alpha_g(x))
				=\gamma_g(\varphi(x))
				=\gamma_g(\pi_\varphi(x)h)
				=\gamma_g(\pi_\varphi(x))h.
			\end{align*} 
			Hence $(\pi_\varphi(\alpha_g(x))-\gamma_g(\pi_\varphi(x)))h=0$. From the continuous functional calculus for positive elements, it is easy to see that $(\pi_\varphi(\alpha_g(x))-\gamma_g(\pi_\varphi(x)))f(h)=0$, and $f(h)$ is also $\beta$-invariant. 
			
			Therefore, we have $\pi_\varphi(\alpha_g(x))f(h)=\gamma_g(\pi_\varphi(x))f(h)$, i.e., $f(\varphi)(\alpha_g(x))=\gamma_g(f(\varphi)(x))$. 
		\end{proof}
	\end{lemma}
	\begin{proposition}\label{AS}
		Let $A$ be a unital simple separable infinite-dimensional C*-algebra. Let $\alpha\colon G\to \mathrm{Aut}(A)$ be an action of a second-countable compact group with the weak tracial Rokhlin property with comparison. Then $\alpha$ is saturated.
		\begin{proof}
			Let $\pi\colon G\to \mathrm{U}(\mathcal{H}_\pi)$ belong to $\hat{G}$. By Theorem \ref{CI}, we will show that $\overline{A_2(\pi)^*A_2(\pi)}=(B(\mathcal{H}_\pi)\otimes A)^{\mathrm{Ad}(\pi)\otimes\alpha}$. By Lemma \ref{SS}, $\overline{A_2(\pi)^*A_2(\pi)}$ is a closed two-sided ideal of $(B(\mathcal{H}_\pi)\otimes A)^{\mathrm{Ad}(\pi)\otimes\alpha}$. Therefore, it suffices to prove that $\overline{A_2(\pi)^*A_2(\pi)}$ contains an invertible element.
			
			It follows from the Peter–Weyl theorem that all irreducible unitary representations of a compact group are finite-dimensional. Let $d_\pi=\mathrm{dim}(\mathcal{H}_\pi)$ and identify $\mathcal{H}_\pi$ with $\mathbb{C}^{d_\pi}$. Let $\delta=\min(\frac{1}{4d_\pi^3},\frac{1}{4})$. Endow $C(G)$ with the action $\mathrm{Lt}$. Denote by $\pi^*$ the function $g\mapsto\pi(g)^*$. A simple calculation shows that $\pi^*$ belongs to $C(G)_2(\pi)$. Indeed, fix $g\in G$, for $h\in G$, we have
			\[((\mathrm{id}_{B(\mathcal{H}_\pi)}\otimes \mathrm{Lt}_g)(\pi^*))(h)=\pi(g^{-1}h)^*\]
			and
			\[(\pi^*(\pi(g)\otimes1_{C(G)}))(h)=\pi(h^{-1})\pi(g)=\pi(g^{-1}h)^*.\]
			Thus, $(\mathrm{id}_{B(\mathcal{H}_\pi)}\otimes \mathrm{Lt}_g)(\pi^*)=\pi^*(\pi(g)\otimes1_{C(G)})$, so $\pi^*\in C(G)_2(\pi)$. For $j,k=1,2,\dots,d_\pi$, let $\pi_{j,k}\in C(G)$ be the function $g\mapsto \pi(g)_{j,k}$. Applying Lemma \ref{WTRCS}, we obtain a positive contraction $r\in (A_{\infty,\alpha}\cap A')^{\alpha_\infty}$ and an equivariant completely positive contractive order zero map $\varphi\colon C(G)\to \overline{r(A_{\infty,\alpha}\cap A')r}$ satisfying $\varphi(1)=r$ and $1-r\precsim_{p.w.-A^\alpha}r$. Note that $r$ is $\alpha_\infty$-invariant. By Lemma \ref{OZFE}, it is easy to see that we obtain an equivariant completely positive contractive order zero map $\sqrt{\varphi}$. Then
			\[\mathrm{id}_{B(\mathcal{H}_\pi)}\otimes\sqrt{\varphi}\colon B(\mathcal{H}_\pi)\otimes C(G)\to \overline{(1_{B(\mathcal{H}_\pi)}\otimes r)(B(\mathcal{H}_\pi)\otimes (A_{\infty,\alpha}\cap A'))(1_{B(\mathcal{H}_\pi)}\otimes r)}\]
			is an equivariant completely positive contractive order zero map with $(\mathrm{id}_{B(\mathcal{H}_\pi)}\otimes\sqrt{\varphi})(1)=1_{B(\mathcal{H}_\pi)}\otimes r^\frac{1}{2}$.
			
			Since $C(G)$ is amenable and separable, by the Choi–Effros lifting theorem, for $\delta$, $F=\{1_A\}$, and $S=\{\pi_{j,k},(\pi_{j,k})^*\colon j,k=1,2,\dots,d_\pi\}$, and $x=y=1$, we obtain a positive contraction $g\in A^\alpha$ such that $(1-g-\delta)_+\precsim_{A^\alpha} g$ and an $(S,F,\delta)$-approximately equivariant central strongly order zero completely positive contractive map $\psi\colon C(G)\to \overline {gAg}$ with $\psi(1)=g^\frac{1}{2}$ such that for all $f_1,f_2\in S$, $\|g^\frac{1}{2}\psi(f_1f_2)-\psi(f_1)\psi(f_2)\|<\delta$. As in Lemma \ref{WTRE}, we will assume that $\psi$ is equivariant. Then
			\[\mathrm{id}_{B(\mathcal{H}_\pi)}\otimes\psi\colon B(\mathcal{H}_\pi)\otimes C(G)\to \overline{(1_{B(\mathcal{H}_\pi)}\otimes g)(B(\mathcal{H}_\pi)\otimes A)(1_{B(\mathcal{H}_\pi)}\otimes g)}.\]
			is an equivariant completely positive contractive map, and $(\mathrm{id}_{B(\mathcal{H}_\pi)}\otimes\psi)(1)=1_{B(\mathcal{H}_\pi)}\otimes g^\frac{1}{2}$. 
			
			Let $(e_{j,k})_{j,k=1,2,\dots,d_\pi}$ be the standard matrix units of $B(\mathcal{H}_\pi)$. Let
			\[c=(\mathrm{id}_{B(\mathcal{H}_\pi)}\otimes\psi)(\pi^*)=\sum_{j,k=1}^{d_\pi}e_{j,k}\otimes\psi((\pi_{k,j})^*)\in B(\mathcal{H}_\pi)\otimes \overline{gAg}.\]
			Since $\mathrm{id}_{B(\mathcal{H}_\pi)}\otimes\psi$ is equivariant and $\pi^*\in C(G)_2(\pi)$, a simple calculation shows that $c\in A_2(\pi)$. Indeed, for $g\in G$, since $\pi^*=\sum_{j,k=1}^{d_\pi}e_{j,k}\otimes (\pi_{k,j})^*\in C(G)_2(\pi)$, we have
			\begin{align*}
				(\mathrm{id}_{B(\mathcal{H}_\pi)}\otimes\alpha_g)(c)=&{}(\mathrm{id}_{B(\mathcal{H}_\pi)}\otimes\alpha_g)((\mathrm{id}_{B(\mathcal{H}_\pi)}\otimes\psi)(\pi^*))\\
				=&{}(\mathrm{id}_{B(\mathcal{H}_\pi)}\otimes\psi)((\mathrm{id}_{B(\mathcal{H}_\pi)}\otimes \mathrm{Lt}_g)(\pi^*))\\
				=&{}(\mathrm{id}_{B(\mathcal{H}_\pi)}\otimes\psi)(\pi^*(\pi(g)\otimes 1_{C(G)}))\\
				=&{}(\mathrm{id}_{B(\mathcal{H}_\pi)}\otimes\psi)((\sum_{j,k=1}^{d_\pi}e_{j,k}\otimes (\pi_{k,j})^*)(\pi(g)\otimes 1_{C(G)}))\\
				=&{}(\mathrm{id}_{B(\mathcal{H}_\pi)}\otimes\psi)(\sum_{j,k=1}^{d_\pi}(e_{j,k}\pi(g))\otimes (\pi_{k,j})^*))\\
				=&{}\sum_{j,k=1}^{d_\pi}\mathrm{id}_{B(\mathcal{H}_\pi)}(e_{j,k}\pi(g))\otimes\psi((\pi_{k,j})^*))\\
				=&{}\sum_{j,k=1}^{d_\pi}(e_{j,k}\pi(g))\otimes\psi((\pi_{k,j})^*))\\
				=&{}(\mathrm{id}_{B(\mathcal{H}_\pi)}\otimes\psi)(\sum_{j,k=1}^{d_\pi}e_{j,k}\otimes (\pi_{k,j})^*)(\pi(g)\otimes 1_A)\\
				=&{}(\mathrm{id}_{B(\mathcal{H}_\pi)}\otimes\psi)(\pi^*)(\pi(g)\otimes 1_A)\\
				=&{}c(\pi(g)\otimes 1_A).
			\end{align*}
			Therefore, $c\in A_2(\pi)$, and $c^*c\in\overline{A_2(\pi)^*A_2(\pi)}$. Furthermore, since $\pi\pi^*=1$, we have
			\[\|c^*c-1_{B(\mathcal{H}_\pi)}\otimes g\|\leq\sum_{j,k,l=1}^{d_\pi}\|\psi(\pi_{j,k})\psi((\pi_{l,k})^*)-g^\frac{1}{2}\psi(\pi_{j,k}(\pi_{l,k})^*)\|<d_\pi^3\delta.\]
			Using $(1-g-\delta)_+\precsim_{A^\alpha} g$, by (1) of Lemma \ref{CCP}, we obtain $\delta_2>0$ such that $(1-g-2\delta)_+\precsim_{A^\alpha}(g-\delta_2)_+$. Then, by (8) of Lemma \ref{CCP}, we have $v\in A^\alpha$ such that $(1-g-2\delta)_+=v^*gv$. Clearly, $1_{B(\mathcal{H}_\pi)}\otimes v\in (B(\mathcal{H}_\pi)\otimes A)^{\mathrm{Ad}(\pi)\otimes\alpha}$. Since $\overline{A_2(\pi)^*A_2(\pi)}$ is a two-sided ideal of $(B(\mathcal{H}_\pi)\otimes A)^{\mathrm{Ad}(\pi)\otimes\alpha}$, it follows that
			\[(1_{B(\mathcal{H}_\pi)}\otimes v)^*c^*c(1_{B(\mathcal{H}_\pi)}\otimes v)\in\overline{A_2(\pi)^*A_2(\pi)}.\]
			Moreover, we have
			\[\|(1_{B(\mathcal{H}_\pi)}\otimes v)^*c^*c(1_{B(\mathcal{H}_\pi)}\otimes v)-1_{B(\mathcal{H}_\pi)}\otimes (1-g-2\delta)_+\|<d_\pi^3\delta.\]
			Therefore,
			\[z=c^*c+(1_{B(\mathcal{H}_\pi)}\otimes v)^*c^*c(1_{B(\mathcal{H}_\pi)}\otimes v)\]
			satisfies
			\[\|z-1_{B(\mathcal{H}_\pi)}\otimes1\|<2(d_\pi^3\delta+\delta)\leq1.\]
			Hence, $\overline{A_2(\pi)^*A_2(\pi)}$ contains the invertible element $z$.
		\end{proof}
	\end{proposition}
	
	\begin{theorem}
		Let $A$ be a unital simple separable infinite-dimensional C*-algebra. Let $\alpha\colon G\to \mathrm{Aut}(A)$ be an action of a second-countable compact group with the weak tracial Rokhlin property with comparison. Then the crossed product $A\rtimes_\alpha G$ is simple.
		\begin{proof}
			By Theorem \ref{FAS}, $A^\alpha$ is simple, and by Proposition \ref{AS}, $\alpha$ is saturated. Therefore, by Theorem \ref{CI}, $A\rtimes_\alpha G$ is simple.
		\end{proof}
	\end{theorem}
	
	\begin{corollary}\label{MESI}
		Let $A$ be a unital simple separable infinite-dimensional C*-algebra. Let $\alpha\colon G\to \mathrm{Aut}(A)$ be an action of a second-countable compact group with the weak tracial Rokhlin property with comparison. Then the crossed product $A\rtimes_\alpha G$ and the fixed point algebra $A^\alpha$ are strongly Morita equivalent and stably isomorphic.
		\begin{proof}
			Since $A\rtimes_\alpha G$ is simple, $A^\alpha$ is a full corner of $A\rtimes_\alpha G$ (see \cite[Corollary]{Ros1979}), hence they are strongly Morita equivalent. Furthermore, since both $A^\alpha$ and $A\rtimes_\alpha G$ are separable, they are stably isomorphic.
		\end{proof}
	\end{corollary}

	\section{Permanence properties}
	\label{sec8.3}
	In this section, we prove that for a unital simple separable infinite-dimensional C*-algebra under the action of a compact group with the weak tracial Rokhlin property with comparison, pure infiniteness, tracial $\mathcal{Z}$-stability, and the combination of amenability with $\mathcal{Z}$-stability can be transferred from the original algebra to the fixed point algebra and the crossed product algebra. We also prove that the radius of comparison of the fixed point algebra does not exceed that of the original algebra.
	
	The following proposition tells us that the compact group action with the weak tracial Rokhlin property of comparison is preserved under taking tensor products. However, for the technical handling of Cuntz comparison, it seems we can only require the action of one group to be trivial. Furthermore, since we need to use Kirchberg's Slice Lemma, we employ the minimal tensor product.
	
	\begin{proposition}
		Let $A,B$ be unital simple infinite-dimensional C*-algebras, let $G$ be a second-countable compact group, and let $\alpha\colon G\to \mathrm{Aut}(A)$ be a group action with the weak tracial Rokhlin property with comparison and $\beta\colon G\to \mathrm{Aut}(B)$ be a trivial group action. Then $\alpha\otimes\beta\colon G\to \mathrm{Aut}(A\otimes_{min}B)$ has the weak tracial Rokhlin property with comparison.
		\begin{proof}
			Let $F\subseteq A\otimes_{min}B$ and $S\subseteq C(G)$ be finite sets, let $\varepsilon>0$, $x\in (A\otimes_{min}B)_+$ with $\|x\|=1$, and $y\in (A\otimes_{min}B)^{\alpha\otimes\beta}_+\setminus\{0\}$. We need to find a positive contraction $d\in (A\otimes_{min}B)^{\alpha\otimes\beta}$ and a completely positive contractive map $\psi\colon C(G)\to \overline{d(A\otimes_{min}B)d}$ such that the following hold:
			\begin{enumerate}[{(}1{)}]
				\item $\psi$ is an $(S,F,\varepsilon)$-approximately equivariant central strongly order zero map.\label{p8.3.1}
				\item $\psi(1)=d$.\label{p8.3.2}
				\item $(1-d-\varepsilon)_+\precsim_{A\otimes_{min}B} x$.\label{p8.3.3}
				\item $(1-d-\varepsilon)_+\precsim_{(A\otimes_{min}B)^{\alpha\otimes\beta}} y$.\label{p8.3.4}
				\item $(1-d-\varepsilon)_+\precsim_{(A\otimes_{min}B)^{\alpha\otimes\beta}} d$.\label{p8.3.5}
				\item $\|dxd\|>1-\varepsilon$.\label{p8.3.6}
			\end{enumerate}
			We may assume there exist $c_1,c_2,\dots,c_m\in A$ and $e_1,e_2,\dots,e_m\in B$ such that $F=\{c_i\otimes e_i\colon i=1,2,\dots,m\}$ and $\|c_i\|\leq1,\|e_i\|\leq1$ for $i=1,2,\dots,m$. Choose $\delta$ such that $0<\delta<\frac{\varepsilon}{3}$ and $(1-5\delta-\delta^2)>1-\varepsilon$. Also, choose $\delta_1$ such that $\frac{1}{2}<\delta_1<1$ and $(1-5\delta-\delta^2)\delta_1^2>1-\varepsilon$. Set $z=(x-\delta_1)_+$. Since $\beta$ is trivial, $(A\otimes_{min}B)^{\alpha\otimes\beta}\cong A^\alpha\otimes_{min} B^\beta$. By Kirchberg's Slice Lemma (see Lemma \ref{KS}), there exist $x_1\in A_+$ and $x_2\in B_+$ such that $x_1\otimes x_2\precsim_{A\otimes_{min}B}z$ with $\|x_1\|=\|x_2\|=1$, and $y_1\in (A^\alpha)_+$ and $y_2\in (B^\beta)_+$ such that $y_1\otimes y_2\precsim_{(A\otimes_{min}B)^{\alpha\otimes\beta}}y$ with $\|y_1\|=\|y_2\|=1$. By \cite[Lemma 2.2, Lemma 2.4(i)]{KR2002}, there exists $w\in A\otimes_{min}B$ such that $\|wxw^*-x_1\otimes x_2\|<\delta^2$ and $\|w\|\leq\delta_1^{-\frac{1}{2}}<\delta_1^{-1}$. Thus there exists $v\in A\otimes_{min}B$, where $v=\sum_{i=1}^kv_i\otimes w_i$ for some $v_i\in A$ and $w_i\in B$, $i=1,2,\dots,k$ such that
			\begin{enumerate}[{(}1{)}]
				\setcounter{enumi}{6}
				\item $\|vxv^*-x_1\otimes x_2\|<\delta^2\ \mbox{and}\ \|v\|<\delta_1^{-1}$.\label{p8.3.7}
			\end{enumerate}
			Set $E=\{c_i\colon i=1,2,\dots,m\}\cup\{v_i\colon i=1,2,\dots,k\}$. Since $B$ and $B^\beta$ are simple and unital, by \cite[Proposition 2.7(v)]{KR2000}, there exist $n,m\in\mathbb{N}$ such that
			\begin{enumerate}[{(}1{)}]
				\setcounter{enumi}{7}
				\item $1_B\precsim_Bx_2\otimes1_n\ \mbox{and}\ 1_B\precsim_{B^\beta}y_2\otimes1_m$.\label{p8.3.8}
			\end{enumerate}
			Since $A$ and $A^\alpha$ are non-type I, by Lemma \ref{OCP}, there exist nonzero $x_0\in A_+$ and $y_0\in (A^\alpha)_+$ such that 
			\begin{enumerate}[{(}1{)}]
				\setcounter{enumi}{8}
				\item $x_0\otimes1_n\precsim_A x_1\ \mbox{and}\ y_0\otimes 1_m\precsim_{A^\alpha}y_1$.\label{p8.3.9}
			\end{enumerate}
			Set $M=1+\sum_{i=1}^k\|v_i\|+\sum_{i=1}^k\|w_i\|$. Choose $\eta>0$ such that $\eta<\frac{\delta}{Mk}$. Applying Lemma \ref{WTRD} to $\alpha$ for $\eta$, $E$, $S$, $x_0$, $x_1$, and $y_1$, we obtain a positive contraction $r\in A^\alpha$ and a completely positive contractive map $\varphi\colon C(G)\to \overline{rAr}$ such that the following hold:
			\begin{enumerate}[{(}1{)}]
				\setcounter{enumi}{9}
				\item $\varphi$ is an $(S,E,\eta)$-approximately equivariant central strongly order zero map.\label{p8.3.10}
				\item $\varphi(1)=r$.\label{p8.3.11}
				\item $(1-r-\eta)_+\precsim_A x_0$.\label{p8.3.12}
				\item $(1-r-\eta)_+\precsim_{A^\alpha} y_0$.\label{p8.3.13}
				\item $(1-r-\eta)_+\precsim_{A^\alpha} r$.\label{p8.3.14}
				\item $\|rx_1r\|>1-\eta$. \label{p8.3.15}
			\end{enumerate}
			Set $d=r\otimes1_B\in (A\otimes_{min}B)^{\alpha\otimes\beta}$ and define a completely positive contractive map $\psi\colon C(G)\to \overline{d(A\otimes_{min}B)d}$ by $\psi(f)=\varphi(f)\otimes1_B$ for all $f\in C(G)$. Clearly $\psi(1)=d$, which is (\ref{p8.3.2}).\\
			For $f\in S$, we have
			\begin{align*}
				&{}\sup_{g\in G}\|\psi(\mathrm{Lt}_g(f))-(\alpha_g\otimes\beta_g)(\psi(f))\|\\
				=&{}\sup_{g\in G}\|\varphi(\mathrm{Lt}_g(f))\otimes1_B-\alpha_g(\varphi(f))\otimes\beta_g(1_B)\|\\
				=&{}\sup_{g\in G}\|\varphi(\mathrm{Lt}_g(f))-\alpha_g(\varphi(f))\|\\
				<&{}\eta<\varepsilon.
			\end{align*}
			For $f\in S$ and $a=c_i\otimes e_i\in F$, we have
			\begin{align*}
				&{}\|\psi(f)a-a\psi(f)\|\\
				=&{}\|(\varphi(f)\otimes1_B)(c_i\otimes e_i)-(c_i\otimes e_i)(\varphi(f)\otimes1_B)\|\\
				=&{}\|(\varphi(f)c_i)\otimes e_i-(c_i\varphi(f))\otimes e_i\|\\
				\leq&{}\|\varphi(f)c_i-c_i\varphi(f)\|\\
				<&{}\eta<\varepsilon.
			\end{align*}
			For $f_1,f_2\in S$, we have
			\begin{align*}
				&{}\|\psi(1)\psi(f_1f_2)-\psi(f_1)\psi(f_2)\|\\
				=&{}|(\varphi(1)\otimes1_B)(\varphi(f_1f_2)\otimes1_B)-(\varphi(f_1)\otimes1_B)(\varphi(f_2)\otimes1_B)\|\\
				=&{}\|(\varphi(1)\varphi(f_1f_2)-\varphi(f_1)\varphi(f_2))\otimes1_B\|\\
				=&{}\|\varphi(1)\varphi(f_1f_2)-\varphi(f_1)\varphi(f_2)\|\\
				<&{}\eta<\varepsilon.
			\end{align*}
			Thus (\ref{p8.3.1}) holds.
			
			Since
			\[\|(1-d)-(1-r-\eta)_+\otimes1_B\|<\eta<\varepsilon,\]
			using part (4) of Lemma \ref{CCP} in the first step, (\ref{p8.3.12}) and (\ref{p8.3.8}) in the second step, (\ref{p8.3.9}) in the fourth step, we have
			\begin{align*}
				(1-d-\varepsilon)_+&{}\precsim_{A\otimes_{min}B}(1-r-\eta)_+\otimes1_B\\
				&{}\precsim_{A\otimes_{min}B}x_0\otimes(x_2\otimes1_n)\\
				&{}\precsim_{A\otimes_{min}B}(x_0\otimes1_n)\otimes x_2\\
				&{}\precsim_{A\otimes_{min}B}x_1\otimes x_2\\
				&{}\precsim_{A\otimes_{min}B}z\\
				&{}\precsim_{A\otimes_{min}B}x.
			\end{align*}
			Similarly, using part (4) of Lemma \ref{CCP} in the first step, (\ref{p8.3.13}) and (\ref{p8.3.8}) in the second step, (\ref{p8.3.9}) in the fourth step, we have
			\begin{align*}
				(1-d-\varepsilon)_+&{}\precsim_{(A\otimes_{min}B)^{\alpha\otimes\beta}}(1-r-\eta)_+\otimes1_B\\
				&{}\precsim_{(A\otimes_{min}B)^{\alpha\otimes\beta}}y_0\otimes(y_2\otimes1_m)\\
				&{}\precsim_{(A\otimes_{min}B)^{\alpha\otimes\beta}}(y_0\otimes1_m)\otimes y_2\\
				&{}\precsim_{(A\otimes_{min}B)^{\alpha\otimes\beta}}y_1\otimes y_2\\
				&{}\precsim_{(A\otimes_{min}B)^{\alpha\otimes\beta}}y.
			\end{align*}
			Moreover, we have
			\begin{align*}
				(1-d-\varepsilon)_+&{}\precsim_{(A\otimes_{min}B)^{\alpha\otimes\beta}}(1-r-\eta)_+\otimes1_B\\
				&{}\precsim_{(A\otimes_{min}B)^{\alpha\otimes\beta}}r\otimes1_B=d.
			\end{align*}
			To prove $\|dxd\|>1-\varepsilon$, first, by (\ref{p8.3.10}) and (\ref{p8.3.11}), we have
			\begin{align*}
				\|dv-vd\|&{}=\|\sum_{i=1}^k(rv_i)\otimes w_i-\sum_{i=1}^k(v_ir)\otimes w_i\|\\
				&{}\leq\sum_{i=1}^k\|(rv_i-v_ir)\otimes w_i\|\\
				&{}\leq Mk\eta<\delta.
			\end{align*}
			And by (\ref{p8.3.15}), we also have
			\begin{align*}
				\|d(x_1\otimes x_2)d\|&{}=\|(rx_1r)\otimes x_2\|\\
				&{}=\|rx_1r\|\|x_2\|\\
				&{}>1-\eta>1-\delta.
			\end{align*}
			Therefore, from the above two inequalities and (\ref{p8.3.7}), we get
			\begin{align*}
				1-\delta&{}<\|d(x_1\otimes x_2)d\|\\
				&{}\leq\|dvxv^*d\|+\|d(vxv^*-x_1\otimes x_2)d\|\\
				&{}<\|vdxv^*d\|+\|(dv-vd)xv^*d\|+\delta^2\\
				&{}\leq\|vdxdv^*\|+\|vdx(v^*d-dv^*)\|+\delta\|v\|+\delta^2\\
				&{}\leq\|v\|^2\|dxd\|+\delta\|v\|+\delta\|v\|+\delta^2\\
				&{}\leq\delta_1^{-2}\|dxd\|+2\delta\delta_1^{-1}+\delta^2\\
				&{}<\delta_1^{-2}\|dxd\|+4\delta+\delta^2.
			\end{align*} 
			Therefore, by the choice of $\delta_1$, we have
			\[\|dxd\|>(1-5\delta-\delta^2)\delta_1^2>1-\varepsilon.\]
		\end{proof}
	\end{proposition}
	
	\begin{corollary}\label{WTRT}
		Let $A$ be a unital simple infinite-dimensional C*-algebra. Let $\alpha\colon G\to \mathrm{Aut}(A)$ be an action of a second-countable compact group with the weak tracial Rokhlin property with comparison. Let $n\in\mathbb{N}$, then the action of $G$ on $M_n\otimes A$ given by $g\mapsto \mathrm{id}_{M_n}\otimes\alpha_g$ has the weak tracial Rokhlin property with comparison.
	\end{corollary}
	
	Now we turn to pure infiniteness. For convenience, we state the following lemma separately.
	
	\begin{lemma}\label{PI}
		Let $A$ be a unital simple C*-algebra. Then the following are equivalent:
		\begin{enumerate}[{(}1{)}]
			\item $A$ is purely infinite.
			\item For any $a\in A_+\setminus\{0\}$, we have $1\precsim_A a$.
		\end{enumerate}
		\begin{proof}
			Assuming (1), by \cite[Lemma 4.12]{MP2021}, there exists $b\in A$ such that $1=b^*ab\precsim_A a$.
			
			For the converse, assuming (2), by parts (1) and (8) of Lemma \ref{CCP}, we have $b\in A$ such that $1=b^*ab$. By \cite[Lemma 4.12]{MP2021}, $A$ is purely infinite.
		\end{proof}
	\end{lemma}
	
	\begin{theorem}\label{FAPI}
		Let $A$ be a unital simple separable infinite-dimensional purely infinite C*-algebra. Let $\alpha\colon G\to \mathrm{Aut}(A)$ be an action of a second-countable compact group with the weak tracial Rokhlin property with comparison. Then the fixed point algebra $A^\alpha$ is purely infinite.
		\begin{proof}
			By Lemma \ref{PI}, we need to show that for any $a\in (A^\alpha)_+\setminus\{0\}$, we have $1\precsim_{A^\alpha}a$. We may assume $\|a\|=1$. Since $A^\alpha$ is simple and non-type I, by Lemma \ref{OCP}, there exist two nonzero mutually orthogonal Cuntz equivalent positive contractions $a_0,a_1$ in $A^\alpha$ with $\|a_0\|=1$ and $\|a_1\|=1$ such that $a_0+a_1\precsim_{A^\alpha}a$.
			
			It suffices to show that for any $\varepsilon>0$, $1\sim_{A^\alpha}(1-\varepsilon)_+\precsim_{A^\alpha}a$.
			
			Set $F=\{a_1\}$. For any $n\in\mathbb{N}$, suppose $\{S_n\}$ is a sequence of finite sets in $A$ such that $\cup_{n\in\mathbb{N}}S_n$ is dense in $A$ and $F\subseteq S_n$ for all $n\in\mathbb{N}$. For $\frac{1}{n}$, $S_n$, $F_2$, $x$, and $y$, there exist a positive contraction $d_n\in A^\alpha$ and a completely positive contractive map $\varphi_n\colon A\to \overline{d_nA^\alpha d_n}$ such that the following hold:
			\begin{enumerate}[{(}1{)}]
				\item $\varphi_n(1)=d_n$.
				\item $\|\varphi_n(a)\varphi_n(b)-d_n\varphi_n(ab)\|<\frac{1}{n}$ for all $a,b\in S_n$.
				\item $\|d_na-ad_n\|<\frac{1}{n}$ for all $a\in S_n$.
				\item $\|\varphi_n(a)-d_n^\frac{1}{2}ad_n^\frac{1}{2}\|<\frac{1}{n}$ for all $a\in F$.
				\item $(1-d_n-\frac{1}{n})_+\precsim_{A^\alpha} a_0$.
			\end{enumerate}
			Let $\pi_{A^\alpha}\colon l^\infty(\mathbb{N},A^\alpha)\to (A^\alpha)_\infty$ be the quotient map. Denote $d=(d_n)_{n=1}^\infty\in l^\infty(\mathbb{N},A^\alpha)$ and $r=\pi_{A^\alpha}(d)$. Define a completely positive contractive order zero map $\varphi\colon A\to \overline{r(A^\alpha)_\infty r}$ by $\varphi(a)=\pi_{A^\alpha}(\{\varphi_1(a),\varphi_2(a),\dots\})$ for all $a\in A$. We have the following hold:
			\begin{enumerate}[{(}1{)}]
				\setcounter{enumi}{5}
				\item $\varphi(1)=r$.
				\item $ar=ra$ for all $a\in A$.
				\item $\varphi(a)=r^\frac{1}{2}ar^\frac{1}{2}$ for all $a\in F$.
				\item $1-r\precsim_{p.w.-A^\alpha}a_0$.
			\end{enumerate}
			Since $A$ is simple and purely infinite, by \cite[Lemma 4.12]{MP2021}, there exists $x\in A$ such that $x^*a_1x=1$. Therefore, by Theorem \ref{OZAH} and Corollary \ref{OZF}, we have
			\begin{align*}
				r=\varphi(1)=&{}\varphi(x^*a_1x)\\
				=&{}r\pi_\varphi(x^*a_1x)\\
				=&{}\pi_\varphi(x^*)r\pi_\varphi(a_1)\pi_\varphi(x)\\
				=&{}\pi_\varphi(x^*)\varphi(a_1)\pi_\varphi(x)\\
				=&{}\pi_\varphi(x^*)r^\frac{1}{2}a_1r^\frac{1}{2}\pi_\varphi(x)\\
				=&{}\sqrt{\varphi}(x^*)a_1\sqrt{\varphi}(x).
			\end{align*}
			We can lift $\sqrt{\varphi}(x)$ and $r$ to obtain $y\in A^\alpha$ and $d_{n_0}\in A^\alpha$ such that $(1-d_{n_0}-\frac{\varepsilon}{2})_+\precsim_{A^\alpha} a_0$ and
			\[\|d_{n_0}-y^*a_1y\|<\frac{\varepsilon}{2}.\]
			By part (4) of Lemma \ref{CCP}, we have
			\[(d_{n_0}-\frac{\varepsilon}{2})_+\precsim_{A^\alpha}y^*a_1y\precsim_{A^\alpha}a_1.\]
			Therefore, by Lemma \ref{CPA} and parts (5) and (6) of Lemma \ref{CCP}, we have
			\begin{align*}
				1\sim_{A^\alpha}(1-\varepsilon)_+&{}=((1-d_{n_0})+d_{n_0}-\varepsilon)_+\\
				&{}\precsim_{A^\alpha}(1-d_{n_0}-\frac{\varepsilon}{2})_++(d_{n_0}-\frac{\varepsilon}{2})_+\\
				&{}\precsim_{A^\alpha}a_0+a_1\\
				&{}\precsim_{A^\alpha}a.
			\end{align*}
		\end{proof}
	\end{theorem}
	
	\begin{corollary}
		Let $A$ be a unital simple separable infinite-dimensional purely infinite C*-algebra. Let $\alpha\colon G\to \mathrm{Aut}(A)$ be an action of a second-countable compact group with the weak tracial Rokhlin property with comparison. Then the crossed product $A\rtimes_\alpha G$ is purely infinite.
		\begin{proof}
			By Theorem \ref{FAPI}, the algebra $A^\alpha$ is purely infinite. Since $A^\alpha$ and $A\rtimes_\alpha G$ are stably isomorphic (see Corollary \ref{MESI}), by \cite[Theorem 4.23]{KR2000}, $A\rtimes_\alpha G$ is also purely infinite.
		\end{proof}
	\end{corollary}
	
	\begin{remark}
		It is obvious that if $A$ is (stably) finite, then $A^\alpha$ is (stably) finite, regardless of the action.
	\end{remark}
	
	Next, we consider $\mathcal{Z}$-stability and tracial $\mathcal{Z}$-stability. The following are preliminaries. We first recall the definition of tracial $\mathcal{Z}$-stability.
	
	\begin{definition}\cite[Definition 2.1]{HO2013}
		Let $A$ be a unital C*-algebra with $A\ncong\mathbb{C}$. Then $A$ is tracially $\mathcal{Z}$-absorbing (or tracially $\mathcal{Z}$-stable) if and only if for any finite set $F\subseteq A$, $\varepsilon>0$, $x\in A_+\setminus\{0\}$, and $n\in\mathbb{N}$, there exists a completely positive contractive order zero map $\psi\colon M_n\to A$ such that
		\begin{enumerate}[{(}1{)}]
			\item $1-\psi(1)\precsim_A x$.
			\item $\|\psi(z)a-a\psi(z)\|<\varepsilon$ for all $z\in M_n$ with $\|z\|=1$ and all $a\in F$.
		\end{enumerate}
	\end{definition}
	
	\begin{proposition}\cite[Proposition 4.16]{MP2021}\label{AMU}
		Let $n\in\mathbb{N}$ and let $\varepsilon>0$. Then there exists $\delta>0$ satisfying the following property. Let $D$ be a C*-algebra and for $j,k=1,2,\dots,n$, let $y_{j,k}\in D$, satisfying (in the unitization of $D$)
		\begin{enumerate}[{(}1{)}]
			\item $\|y_{j,k}y_{k,m}-y_{j,l}y_{l,m}\|<\delta$ for $j,k,l,m=1,2,\dots,n$.\label{p8.4.1}
			\item $\|y_{j,j}y_{k,k}\|<\delta$ for $j,k=1,2,\dots,n$ with $j\neq k$.\label{p8.4.2}
			\item $\|y_{j,k}-y_{k,j}^*\|<\delta$ for $j,k=1,2,\dots,n$.\label{p8.4.3}
			\item $\|y_{j,j}\|<1+\delta$ for $j=1,2,\dots,n$.\label{p8.4.4}
			\item $\|1-y_{j,j}\|<1+\delta$ for $j=1,2,\dots,n$.\label{p8.4.5}
		\end{enumerate}
		Then there exists a completely positive contractive order zero map $\rho\colon M_n\rightarrow D$ such that, with $M_n$ equipped with $(e_{j,k})_{j,k=1,2,\dots,n}$ as standard matrix units, for $j,k=1,2,\dots,n$, we have $\|\rho(e_{j,k})-y_{j,k}\|<\varepsilon$.
	\end{proposition}
	
	\begin{theorem}\label{FATZ}
		Let $A$ be a unital simple separable infinite-dimensional C*-algebra that is tracially $\mathcal{Z}$-stable. Let $\alpha\colon G\to \mathrm{Aut}(A)$ be an action of a second-countable compact group with the weak tracial Rokhlin property with comparison. Then the fixed point algebra $A^\alpha$ is tracially $\mathcal{Z}$-stable.
		\begin{proof}
			We need to show that for any $\varepsilon>0$, any finite set $F\subseteq A^\alpha$, any $x\in (A^\alpha)_+\setminus\{0\}$, and any $n\in\mathbb{N}$, there exists a completely positive contractive order zero map $\psi\colon M_n\to A^\alpha$ such that the following hold:
			\begin{enumerate}[{(}1{)}]
				\item $1-\psi(1)\precsim_{A^\alpha} x$.
				\item $\|\psi(z)a-a\psi(z)\|<\varepsilon$ for all $z\in M_n$ with $\|z\|=1$ and all $a\in F$.
			\end{enumerate}
			Without loss of generality, we may assume $\|a\|\leq1$ for all $a\in F$ and $\varepsilon<1$. Set $\varepsilon_0=\frac{\varepsilon}{16n^2}$. Applying Proposition \ref{AMU} to $\varepsilon_0, n$, we obtain $\delta>0$. Since $A^\alpha$ is simple and non-type I, by Lemma \ref{OCP}, there exist $x_1,x_2\in (A^\alpha)_+\setminus\{0\}$ such that $x_1x_2=0$, $x_1\sim_{A^\alpha}x_2$, and $x_1+x_2\in\overline{xA^\alpha x}$.
			
			Since $A$ is tracially $\mathcal{Z}$-stable, for given $n$ and $F$, for $\varepsilon_0$ and for $x_1$, we have a completely positive contractive order zero map $\psi_0\colon M_n\rightarrow A$ such that the following hold:
			\begin{enumerate}[{(}1{)}]
				\setcounter{enumi}{2}
				\item $1-\psi_0(1)\precsim_A x_1$.
				\item $\|\psi_0(z)a-a\psi_0(z)\|<\varepsilon_0$ for all $z\in M_n$ with $\|z\|=1$ and all $a\in F$.\label{t8.7.4}
			\end{enumerate}
			Set $\delta_0=\min(\delta,\frac{\varepsilon}{32})$. Let $(e_{j,k})_{j,k=1,2\dots,n}$ be the standard matrix units for $M_n$. Set
			\[F_0=F\cup\{\psi_0(e_{j,k})\colon j,k\in\{1,2,\dots,n\}\}\cup\{x_1\}.\]
			Suppose $\{S_n\}$ is a sequence of finite sets in $A$ such that $\cup_{n\in\mathbb{N}}S_n$ is dense in $A$ and $F_0\subseteq S_n$ for all $n\in\mathbb{N}$. For $\frac{1}{n}$, $S_n$, $F\cup\{x_1\}$, $x_2$, there exists a positive contraction $d_n\in A^\alpha$ and a completely positive contractive map $\varphi_n\colon A\to \overline{d_nA^\alpha d_n}$ such that the following hold:
			\begin{enumerate}[{(}1{)}]
				\setcounter{enumi}{4}
				\item $\varphi_n(1)=d_n$.
				\item $\|\varphi_n(a)\varphi_n(b)-d_n\varphi_n(ab)\|<\frac{1}{n}$ for all $a,b\in S_n$.
				\item $\|d_na-ad_n\|<\frac{1}{n}$ for all $a\in S_n$.
				\item $\|\varphi_n(a)-d_n^\frac{1}{2}ad_n^\frac{1}{2}\|<\frac{1}{n}$ for all $a\in F\cup\{x_1\}$.
				\item $(1-d_n-\frac{1}{n})_+\precsim_{A^\alpha} x_2$.
			\end{enumerate}
			Let $\pi_{A^\alpha}\colon l^\infty(\mathbb{N},A^\alpha)\to (A^\alpha)_\infty$ be the quotient map. Denote $d=(d_n)_{n=1}^\infty\in l^\infty(\mathbb{N},A^\alpha)$ and $r=\pi_{A^\alpha}(d)$. Define a completely positive contractive order zero map $\varphi\colon A\to \overline{r(A^\alpha)_\infty r}$ by $\varphi(a)=\pi_{A^\alpha}(\{\varphi_1(a),\varphi_2(a),\dots\})$ for all $a\in A$. We have the following hold:
			\begin{enumerate}[{(}1{)}]
				\setcounter{enumi}{9}
				\item $\varphi(1)=r$.
				\item $ar=ra$ for all $a\in A$.
				\item $\varphi(a)=r^\frac{1}{2}ar^\frac{1}{2}$ for all $a\in F\cup\{x_1\}$.
				\item $1-r\precsim_{p.w.-A^\alpha}x_2$.
			\end{enumerate}
			By Corollary \ref{OZCP}, we have $r-\varphi\circ\psi_0(1)=\varphi(1-\psi_0(1))\precsim_{(A^\alpha)_\infty}\varphi(x_1)=r^\frac{1}{2}x_1r^\frac{1}{2}\precsim_{(A^\alpha)_\infty}x_1$. By part (1) of Lemma \ref{CCP}, for $\frac{\varepsilon}{64}$, there exists $\delta_1>0$ such that $(\varphi(1-\psi_0(1))-\frac{\varepsilon}{64})_+\precsim_{(A^\alpha)_\infty}(x_1-\delta_1)_+$. Therefore, by part (8) of Lemma \ref{CCP}, there exists $v\in(A^\alpha)_\infty$ such that
			\[(\varphi(1-\psi_0(1))-\frac{\varepsilon}{64})_+=v^*x_1v.\]
			Let $a\in F$ and $z\in M_n$ with $\|z\|=1$. Then by Lemma \ref{OZAH}, we have
			\begin{align*}
				\varphi(a\psi_0(z))=&{}r\pi_\varphi(a\psi_0(z))\\
				=&{}r\pi_\varphi(a)\pi_\varphi(\psi_0(z))\\
				=&{}\varphi(a)\pi_\varphi(\psi_0(z))\\
				=&{}ar\pi_\varphi(\psi_0(z))\\
				=&{}a\varphi\circ\psi_0(z).
			\end{align*}
			Similarly, we have $\varphi(\psi_0(z)a)=\varphi\circ\psi_0(z)a$.
			Therefore, for all $a\in F$ and $z\in M_n$ with $\|z\|=1$, by (\ref{t8.7.4}) and (\ref{t8.7.eq1}), we have
			\begin{equation}\label{t8.7.eq1}
				\|a\varphi\circ\psi_0(z)-\varphi\circ\psi_0(z)a\|=\|\varphi(a\psi_0(z)-\psi_0(z)a)\|<\varepsilon_0<\frac{\varepsilon}{16}. 
			\end{equation}
			Thus we can choose $u=(u_n)_{n\in\mathbb{N}}\in l^\infty(\mathbb{N},A^\alpha)$ such that $\pi_{A^\alpha}(u)=v$ and a completely positive contractive map $\varphi_{n_0}\colon A\to \overline{d_{n_0}A^\alpha d_{n_0}}$ such that the following hold:
			\begin{enumerate}[{(}1{)}]
				\setcounter{enumi}{13}
				\item $\varphi_{n_0}(1)=d_{n_0}$.
				\item $\|\varphi_{n_0}(a)\varphi_{n_0}(b)-d_{n_0}\varphi_{n_0}(ab)\|<\delta_0$ for all $a,b\in S_{n_0}$.\label{t8.7.15}
				\item $\|d_{n_0}a-ad_{n_0}\|<\delta_0$ for all $a\in S_{n_0}$.
				\item $\|\varphi_{n_0}(a)-d_{n_0}^\frac{1}{2}ad_{n_0}^\frac{1}{2}\|<\delta_0$ for all $a\in F\cup\{x_1\}$.
				\item $(1-d_{n_0}-\frac{\varepsilon}{32})_+\precsim_{A^\alpha} x_2$.
				\item $\|(\varphi_{n_0}(1-\psi_0(1))-\frac{\varepsilon}{64})_+-u_{n_0}^*x_1u_{n_0}\|<\frac{\varepsilon}{64}$.\label{t8.7.19}
				\item $\|a\varphi_{n_0}\circ\psi_0(z)-\varphi_{n_0}\circ\psi_0(z)a\|<\frac{\varepsilon}{2}$ for all $a\in F$ and all $z\in M_n$ with $\|z\|=1$.\label{t8.7.20}
			\end{enumerate}
			For $j,k=1,2,\dots,n$, define $y_{j,k}=(\varphi_{n_0}\circ\psi_0)(e_{j,k})$. By \cite[Lemma 2.6(3)]{ABP2018}, for $j,k,l,m=1,2,\dots,n$, we have $\psi_0(e_{j,k})\psi_0(e_{k,m})=\psi_0(e_{j,l})\psi_0(e_{l,m})$. Then from (\ref{t8.7.15}) it is easy to see that for $j,k=1,2,\dots,n$, $y_{j,k}$ satisfy relations (\ref{p8.4.1}) and (\ref{p8.4.2}) in Proposition \ref{AMU}. Also, stronger versions of (\ref{p8.4.3}), (\ref{p8.4.4}), and (\ref{p8.4.5}) in Proposition \ref{AMU} hold immediately: $y_{j,k}=y_{k,j}^*$ for all $j,k=1,2,\dots,n$ and $\|y_{j,j}\|\leq 1$, $\|1-y_{j,j}\|\leq 1$ for all $j=1,2,\dots,n$. By the choice of $\delta$, there exists a completely positive contractive order zero map $\psi_1\colon M_n\to \overline{d_{n_0}A^\alpha d_{n_0}}$ such that for $j,k=1,2,\dots,n$,
			\[\|\psi_1(e_{j,k})-y_{j,k}\|<\varepsilon_0.\]
			It follows that for all $z\in M_n$,
			\begin{equation}\label{t8.7.eq2}
				\|\psi_1(z)-(\varphi_{n_0}\circ\psi_0)(z)\|< n^2\varepsilon_0\|z\|.
			\end{equation}
			From (\ref{t8.7.19}) and part (4) of Lemma \ref{CCP}, we have
			\[(d_{n_0}-\varphi_{n_0}\circ\psi_0(1)-\frac{\varepsilon}{32})_+=(\varphi_{n_0}(1-\psi_0(1))-\frac{\varepsilon}{32})_+\precsim_{A^\alpha}x_1.\]
			Therefore, by Lemma \ref{CPA} and parts (5) and (6) of Lemma \ref{CCP},
			\begin{align*}
				(1-\varphi_{n_0}\circ\psi_0(1)-\frac{\varepsilon}{16})_+&{}=((1-d_{n_0})+(d_{n_0}-\psi_{n_0}\circ\psi_0(1))-\frac{\varepsilon}{16})_+\\
				&{}\precsim_{A^\alpha}(1-d_{n_0}-\frac{\varepsilon}{32})_++(d_{n_0}-\varphi_{n_0}\circ\psi_0(1)-\frac{\varepsilon}{32})_+\\
				&{}\precsim_{A^\alpha}x_2+x_1\\
				&{}\precsim_{A^\alpha}x.
			\end{align*}
			Since
			\[\|(1-\psi_1(1))-(1-\varphi_{n_0}\circ\psi_0(1))\|<n^2\varepsilon_0=\frac{\varepsilon}{16},\]
			by part (4) of Lemma \ref{CCP}, we have
			\[(1-\psi_1(1)-\frac{\varepsilon}{8})_+\precsim_{A^\alpha}(1-\varphi_{n_0}\circ\psi_0(1)-\frac{\varepsilon}{16})_+\precsim_{A^\alpha}x.\]
			By \cite[Lemma 4.17]{MP2021}, we obtain a completely positive contractive order zero map $\psi\colon M_n\to \overline{d_{n_0}A^\alpha d_{n_0}}$ such that
			\begin{equation}\label{t8.7.eq3}
				\|\psi-\psi_1\|<\frac{\varepsilon}{8}
			\end{equation}
			and
			\[(1-\frac{\varepsilon}{8})(1-\psi(1))=(1-\psi_1(1)-\frac{\varepsilon}{8})_+.\]
			Then for all $z\in M_n$, by (\ref{t8.7.eq2}) and (\ref{t8.7.eq3}) we have
			\[\|\psi(z)-(\varphi_{n_0}\circ\psi_0)(z)\|<(\frac{\varepsilon}{8}+n^2\varepsilon_0)\|z\|<\frac{\varepsilon}{4}\|z\|.\]
			Let $a\in F$ and $z\in M_n$ with $\|z\|=1$, by (\ref{t8.7.20}) and (\ref{t8.7.eq1}), we have
			\begin{align*}
				&{}\|a\psi(z)-\psi(z)a\|\\
				\leq&{}\|a\psi(z)-a(\varphi_{n_0}\circ\psi_0)(z)\|+\|a(\varphi_{n_0}\circ\psi_0)(z)-(\varphi_{n_0}\circ\psi_0)(z)a\|\\
				&{}+\|(\varphi_{n_0}\circ\varphi_0)(z)a-\psi(z)a\|\\
				\leq&{}2\|a\|\|\psi(z)-(\varphi_{n_0}\circ\psi_0)(z)\|+\|a(\varphi_{n_0}\circ\psi_0)(z)-(\varphi_{n_0}\circ\psi_0)(z)a\|\\
				<&{}\frac{\varepsilon}{2}+\frac{\varepsilon}{2}\\
				=&{}\varepsilon.
			\end{align*}
			Similarly, we also have
			\[1-\psi(1)=(1-\frac{\varepsilon}{8})^{-1}(1-\psi_1(1)-\frac{\varepsilon}{8})_+\precsim_{A^\alpha}x.\]
		\end{proof}
	\end{theorem}
	
	\begin{corollary}
		Let $A$ be a unital simple separable infinite-dimensional C*-algebra that is tracially $\mathcal{Z}$-stable. Let $\alpha\colon G\to \mathrm{Aut}(A)$ be an action of a second-countable compact group with the weak tracial Rokhlin property with comparison. Then the crossed product $A\rtimes_\alpha G$ is tracially $\mathcal{Z}$-stable.
		\begin{proof}
			By Theorem \ref{FATZ}, the algebra $A^\alpha$ is tracially $\mathcal{Z}$-stable. Since $A^\alpha$ and $A\rtimes_\alpha G$ are stably isomorphic (see Corollary \ref{MESI}), by \cite[Corollary 4.12]{AGJP2021}, $A\rtimes_\alpha G$ is also tracially $\mathcal{Z}$-stable.
		\end{proof}
	\end{corollary}
	
	\begin{theorem}
		Let $A$ be a unital simple separable infinite-dimensional amenable C*-algebra that is $\mathcal{Z}$-stable. Let $\alpha\colon G\to \mathrm{Aut}(A)$ be an action of a second-countable compact group with the weak tracial Rokhlin property with comparison. Then both the fixed point algebra $A^\alpha$ and the crossed product $A\rtimes_\alpha G$ are $\mathcal{Z}$-stable.
		\begin{proof}
			The algebra $A$ is $\mathcal{Z}$-stable, hence by \cite[Proposition 2.2]{HO2013} it is tracially $\mathcal{Z}$-stable. Therefore, by Theorem \ref{FATZ}, $A^\alpha$ is tracially $\mathcal{Z}$-stable. $A$ is amenable, so $A^\alpha$ is also amenable; thus by \cite[Theorem 4.1]{HO2013}, $A^\alpha$ is $\mathcal{Z}$-stable. By \cite[Corollary 3.2]{TW2007}, $\mathcal{Z}$-stability is preserved under stable isomorphism, and $A^\alpha$ and $A\rtimes_\alpha G$ are stably isomorphic (see Corollary \ref{MESI}), so the crossed product $A\rtimes_\alpha G$ is also $\mathcal{Z}$-stable.
		\end{proof}
	\end{theorem}
	
	Finally, we consider the radius of comparison. Since it is only defined for unital C*-algebras, we only consider the fixed point algebra. We first introduce some preliminaries.
	\begin{proposition}\label{FAC}
		Let $A$ be a unital simple separable infinite-dimensional C*-algebra. Let $\alpha\colon G\to \mathrm{Aut}(A)$ be an action of a second-countable compact group with the weak tracial Rokhlin property with comparison. Let $n\in\mathbb{N}$, and let $a,b\in M_n(A)_+$. Suppose $0$ is a limit point of $\mathrm{sp}(b)$. Then $a\precsim_A b$ if and only if $a\precsim_{A^\alpha} b$.
		\begin{proof}
			We only need to show that if $a\precsim_A b$ then $a\precsim_{A^\alpha}b$. By Corollary \ref{WTRT}, we may assume $n=1$. Without loss of generality, we may assume $\|a\|\leq1$ and $\|b\|\leq1$.
			
			It suffices to show that for any $\varepsilon>0$, $(a-\varepsilon)_+\precsim_{A^\alpha}b$. By part (1) of Lemma \ref{CCP}, there exists $\delta>0$ such that $(a-\frac{\varepsilon}{4})_+\precsim_A(b-\delta)_+$. Set $a'=(a-\frac{\varepsilon}{4})_+$, $b'=(b-\delta)_+$.
			
			We first prove that if $B$ is a unital C*-algebra, $a,g\in B$ satisfy $0\leq a,g\leq 1$, and let $\varepsilon_1,\varepsilon_2>0$, then
			\begin{equation}\label{p8.5.1}
				(a-(\varepsilon_1+\varepsilon_2))_+\precsim_B (gag-\varepsilon_1)_+ \oplus (1-g^2-\varepsilon_2)_+.
			\end{equation}
			Indeed, using Lemma \ref{CPA} in the first step, part (10) of Lemma \ref{CCP} in the second step, and part (2) of Lemma \ref{CCP} in the third step, we get
			\begin{align*}
				(a-(\varepsilon_1+\varepsilon_2))_+\precsim_B &{}(a^\frac{1}{2}g^2a^\frac{1}{2}-\varepsilon_1)_+\oplus (a^\frac{1}{2}(1-g^2)a^\frac{1}{2}-\varepsilon_2)_+\\
				\sim_B &{}(gag-\varepsilon_1)_+\oplus((1-g^2)^\frac{1}{2}a(1-g^2)^\frac{1}{2}-\varepsilon_2)_+\\
				\precsim_B &{}(gag-\varepsilon_1)_+\oplus (1-g^2-\varepsilon_2)_+.
			\end{align*}
			Choose $\lambda\in \mathrm{sp}(b)\cap (0,\delta)$. Let $f\colon [0,\infty)\to [0,1]$ be a continuous function such that $f(\lambda)=1$ and $\mathrm{supp}(f)\subseteq (0,\delta)$. Then
			\begin{equation}\label{p8.5.eq1}
				\|f(b)\|=1,\ f(b)\perp b',\ \mbox{and} \ f(b)+b'\precsim_{A^\alpha} b.
			\end{equation}
			Set $F=\{a,b,a',b'\}$. Suppose $\{S_n\}$ is a sequence of finite sets in $A$ such that $\cup_{n\in\mathbb{N}}S_n$ is dense in $A$ and $F\subseteq S_n$ for all $n\in\mathbb{N}$. For $\frac{1}{n}$, $S_n$, $F$, $f(b)$, there exist a positive contraction $d_n\in A^\alpha$ and a completely positive contractive map $\varphi_n\colon A\to \overline{d_nA^\alpha d_n}$ such that the following hold:
			\begin{enumerate}[{(}1{)}]
				\item $\varphi_n(1)=d_n$.
				\item $\|\varphi_n(x)\varphi_n(y)-d_n\varphi_n(xy)\|<\frac{1}{n}$ for all $x,y\in S_n$.
				\item $\|d_nx-xd_n\|<\frac{1}{n}$ for all $x\in S_n$.
				\item $\|\varphi_n(a)-d_n^\frac{1}{2}xd_n^\frac{1}{2}\|<\frac{1}{n}$ for all $x\in F$.
				\item $(1-d_n-\frac{1}{n})_+\precsim_{A^\alpha} f(b)$.
			\end{enumerate}
			Let $\pi_{A^\alpha}\colon l^\infty(\mathbb{N},A^\alpha)\to (A^\alpha)_\infty$ be the quotient map. Denote $d=(d_n)_{n=1}^\infty\in l^\infty(\mathbb{N},A^\alpha)$ and $r=\pi_{A^\alpha}(d)$. Define a completely positive contractive order zero map $\varphi\colon A\to \overline{r(A^\alpha)_\infty r}$ by $\varphi(a)=\pi_{A^\alpha}(\{\varphi_1(a),\varphi_2(a),\dots\})$ for all $a\in A$. We have the following hold:
			\begin{enumerate}[{(}1{)}]
				\setcounter{enumi}{5}
				\item $\varphi(1)=r$.
				\item $xr=rx$ for all $x\in A$.
				\item $\varphi(x)=r^\frac{1}{2}xr^\frac{1}{2}$ for all $x\in F$.
				\item $1-r\precsim_{p.w.-A^\alpha}f(b)$.
			\end{enumerate}
			Therefore, by Lemma \ref{OZCP}, we have
			\[r^\frac{1}{2}a'r^\frac{1}{2}=\varphi(a')\precsim_{(A^\alpha)_\infty}\varphi(b')=r^\frac{1}{2}b'r^\frac{1}{2}\precsim_{(A^\alpha)_\infty}b'.\]
			Thus by part (1) of Lemma \ref{CCP}, for $\frac{\varepsilon}{4}$ there exists $\delta_1>0$ such that $(r^\frac{1}{2}a'r^\frac{1}{2}-\frac{\varepsilon}{4})_+\precsim_{(A^\alpha)_\infty}(b'-\delta_1)_+$. Therefore, by part (8) of Lemma \ref{CCP}, there exists $v\in(A^\alpha)_\infty$ such that
			\[(r^\frac{1}{2}a'r^\frac{1}{2}-\frac{\varepsilon}{4})_+=v^*b'v.\]
			So we can choose $u=(u_n)_{n\in\mathbb{N}}\in l^\infty(\mathbb{N},A^\alpha)$ such that $\pi_{A^\alpha}(u)=v$ and a completely positive contractive map $\varphi_{n_0}\colon A\to \overline{d_{n_0}A^\alpha d_{n_0}}$ such that the following hold:
			\begin{enumerate}[{(}1{)}]
				\setcounter{enumi}{9}
				\item $\varphi_{n_0}(1)=d_{n_0}$.
				\item $\|\varphi_{n_0}(x)\varphi_{n_0}(y)-d_{n_0}\varphi_{n_0}(xy)\|<\frac{\varepsilon}{4}$ for all $x,y\in S_{n_0}$. 
				\item $\|d_{n_0}x-xd_{n_0}\|<\frac{\varepsilon}{4}$ for all $x\in S_n$. 
				\item $\|\varphi_{n_0}(x)-d_{n_0}^\frac{1}{2}xd_{n_0}^\frac{1}{2}\|<\frac{\varepsilon}{4}$ for all $x\in F$. 
				\item $(1-d_{n_0}-\frac{\varepsilon}{4})_+\precsim_{A^\alpha} f(b)$.
				\item $\|(d_{n_0}^\frac{1}{2}a'd_{n_0}^\frac{1}{2}-\frac{\varepsilon}{4})_+-u_{n_0}^*b'u_{n_0}\|<\frac{\varepsilon}{4}$. 
			\end{enumerate}
			Therefore, applying part (4) of Lemma \ref{CCP} to (15), we have
			\begin{align*}
				(d_{n_0}^\frac{1}{2}a'd_{n_0}^\frac{1}{2}-\frac{\varepsilon}{2})_+\precsim_{A^\alpha}u_{n_0}^*b'u_{n_0}\precsim_{A^\alpha}b'.
			\end{align*}
			Thus, using part (3) of Lemma \ref{CCP} in the first step, (\ref{p8.5.1}) in the second step, part (6) of Lemma \ref{CCP} in the third step, and (\ref{p8.5.eq1}) in the fourth step, we get
			\begin{align*}
				(a-\varepsilon)_+&{}=(a'-\frac{3\varepsilon}{4})_+\\
				&{}\precsim_{A^\alpha} (d_{n_0}^\frac{1}{2}a'd_{n_0}^\frac{1}{2}-\frac{\varepsilon}{2})_+\oplus (1-d_{n_0}-\frac{\varepsilon}{4})_+\\
				&{}\precsim_{A^\alpha}b'+ f(b)\precsim_{A^\alpha}b.
			\end{align*}
		\end{proof}
	\end{proposition}
	
	The following proposition is a generalization of \cite[Proposition 4.22]{MP2021}. This result is weaker than its counterpart for compact group actions with the Rokhlin property; namely, we cannot obtain an embedding of the Cuntz semigroup of the fixed point algebra into the Cuntz semigroup of the original algebra, but only an embedding between the subsemigroups consisting of purely positive elements (i.e., elements not Cuntz equivalent to a projection).
		\begin{proposition}
		Let $A$ be a unital simple separable infinite-dimensional stably finite C*-algebra. Let $\alpha\colon G\to \mathrm{Aut}(A)$ be an action of a second-countable compact group with the weak tracial Rokhlin property with comparison. Let $\iota\colon A^\alpha\to A$ be the inclusion map. Then
		\begin{enumerate}
			\item the map $\mathrm{W}(\iota)\colon \mathrm{W}(A^\alpha)\to \mathrm{W}(A)$ induces an injective map of ordered semigroups from $\mathrm{W}(A^\alpha)_+\cup\{0\}$ to $\mathrm{W}(A)_+\cup\{0\}$;
			
			\item the map $\mathrm{Cu}(\iota)\colon \mathrm{Cu}(A^\alpha)\to \mathrm{Cu}(A)$ induces an injective map of ordered semigroups from $\mathrm{Cu}(A^\alpha)_+\cup\{0\}$ to $\mathrm{Cu}(A)_+\cup\{0\}$.
		\end{enumerate}
		\begin{proof}
			The proof is essentially the same as that of \cite[Proposition 4.22]{MP2021} when we use Proposition \ref{FAC} in place of \cite[Proposition 4.21]{MP2021}, so we omit it.
		\end{proof}
	\end{proposition}
	
	The following proposition establishes a correspondence between the 2-quasitraces on the fixed point algebra under the action of a second-countable compact group with the weak tracial Rokhlin property of comparison and the 2-quasitraces on the original algebra.
	\begin{proposition}
		Let $A$ be a unital simple separable infinite-dimensional C*-algebra. Let $\alpha\colon G\to \mathrm{Aut}(A)$ be an action of a second-countable compact group with the weak tracial Rokhlin property with comparison. Then the restriction map $\mathrm{QT}(\iota)\colon\mathrm{QT}(A)\to \mathrm{QT}({A^\alpha})$ is surjective.
		\begin{proof}
			Let $\tau'\in \mathrm{QT}(A)$. We will show that the map $\mathrm{QT}(\iota)\colon \mathrm{QT}(A)\to \mathrm{QT}(A^\alpha)$, sending $\tau'\mapsto \tau'\circ\iota$, is surjective. Let $z\in (A^\alpha)_+\setminus\{0\}$. Suppose $\{S_n\}$ is a sequence of finite sets in $A$ such that $\cup_{n\in\mathbb{N}}S_n$ is dense in $A$ and $F\subseteq S_n$ for all $n\in\mathbb{N}$. For $\frac{1}{n}$, $S_n$, $F$, $z$, there exist a positive contraction $d_n\in A^\alpha$ and a completely positive contractive map $\varphi_n\colon A\to \overline{d_nA^\alpha d_n}$ such that the following hold:
			\begin{enumerate}[{(}1{)}]
				\item $\varphi_n(1)=d_n$.
				\item $\|\varphi_n(a)\varphi_n(b)-d_n\varphi_n(ab)\|<\frac{1}{n}$ for all $a,b\in S_n$.
				\item $\|d_na-ad_n\|<\frac{1}{n}$ for all $a\in S_n$.
				\item $\|\varphi_n(a)-d_n^\frac{1}{2}ad_n^\frac{1}{2}\|<\frac{1}{n}$ for all $a\in F$.
				\item $(1-d_n-\frac{1}{n})_+\precsim_{A^\alpha} (z-\frac{1}{2})_+$.
			\end{enumerate}
			Let $\pi_{A^\alpha}\colon l^\infty(\mathbb{N},A^\alpha)\to (A^\alpha)_\infty$ be the quotient map. Denote $d=(d_n)_{n=1}^\infty\in l^\infty(\mathbb{N},A^\alpha)$ and $r=\pi_{A^\alpha}(d)$. Define a completely positive contractive order zero map $\varphi\colon A\to \overline{r(A^\alpha)_\infty r}$ by $\varphi(a)=\pi_{A^\alpha}(\{\varphi_1(a),\varphi_2(a),\dots\})$ for all $a\in A$. We have the following hold:
			\begin{enumerate}[{(}1{)}]
				\setcounter{enumi}{5}
				\item $\varphi(1)=r$.
				\item $ar=ra$ for all $a\in A$.
				\item $\varphi(a)=r^\frac{1}{2}ar^\frac{1}{2}$ for all $a\in F$.
				\item $1-r\precsim_{(A^\alpha)_\infty}z$.
			\end{enumerate}
			For any $n\in\mathbb{N}$, there exist a positive contraction $r_n\in (A^\alpha)_\infty\cap A'$ and a completely positive contractive order zero map $\psi_n \colon A \to \overline{r_n(A^\alpha)_\infty r_n}$, such that for each $v\in \mathrm{QT}((A^\alpha)_\infty)$, $d_{v}(1-r_n)<\frac{1}{n}$. Now, let $\tau\in\mathrm{QT}(A^\alpha)$, which naturally induces a quasitrace on $(A^\alpha)_\infty$, still denoted by $\tau$. We consider the weak-$*$ limit of the quasitraces $\tau_\infty\circ\psi_n$ in $\mathrm{QT}(A)$, denoted $w^*-\lim(\tau_\infty\circ\psi_n)$. Note that for $a\in A_+$ with $\|a\|\leq1$ and $v \in\mathrm{QT}((A^\alpha)_\infty)$, $v(a)\leq d_v(a)$. Then for each $a\in A^\alpha$
			\begin{align*}
				\mathrm{QT}([w^*-\lim(\tau_\infty\circ\psi_n)])(a)=&{}[w^*-\lim(\tau_\infty\circ\psi_n)](\iota(a))\\
				=&{}\lim_n\tau_\infty(\psi_n(\iota(a)))\\
				=&{}\lim_n\tau_\infty(r_n^\frac{1}{2}ar_n^\frac{1}{2})\\
				=&{}\tau_\infty(a)-\lim_n\tau_\infty(a(1-r_n))\\
				=&{}\tau(a)
			\end{align*}
		\end{proof}
	\end{proposition}
	
	\begin{theorem}\label{RC}
		Let $A$ be a unital simple separable infinite-dimensional stably finite C*-algebra. Let $\alpha\colon G\to \mathrm{Aut}(A)$ be an action of a second-countable compact group with the weak tracial Rokhlin property with comparison. Then $\mathrm{rc}(A^\alpha)\leq \mathrm{rc}(A)$.
		\begin{proof}
			The proof is essentially the same as that of \cite[Theorem 4.1]{AGP2021} when we use Proposition \ref{FAC} in place of \cite[Lemma 3.7]{AGP2021}, so we omit it.
		\end{proof}
	\end{theorem}
	
	\begin{corollary}
		Let $A$ be a unital simple separable infinite-dimensional stably finite C*-algebra with strict comparison. Let $\alpha\colon G\to \mathrm{Aut}(A)$ be an action of a second-countable compact group with the weak tracial Rokhlin property with comparison. Then both the fixed point algebra $A^\alpha$ and the crossed product $A\rtimes_\alpha G$ have strict comparison.
		\begin{proof}
			Applying Theorem \ref{RC} for $\mathrm{rc}(A)=0$, we get that $A^\alpha$ has strict comparison. Theorem \ref{MESI} implies that $A\rtimes_\alpha G$ is stably isomorphic to $A^\alpha$. Since stably isomorphic C*-algebras have canonically isomorphic Cuntz semigroups, $A\rtimes_\alpha G$ has strict comparison.
		\end{proof}
	\end{corollary}
	
	Due to Proposition \ref{FAC}, similar to the case of the tracial Rokhlin property with comparison, we have the following comparison properties.
	\begin{theorem}\label{WBC}
		Let $A$ be a unital simple separable infinite-dimensional stably finite C*-algebra with $\beta$-comparison. Let $\alpha\colon G\to \mathrm{Aut}(A)$ be an action of a second-countable compact group with the weak tracial Rokhlin property with comparison. Then the fixed point algebra $A^\alpha$ has $\beta$-comparison.
		\begin{proof}
			The proof is essentially the same as that of \cite[Theorem 4.4]{TF2025}, so we omit it.
		\end{proof}
	\end{theorem}
	
	\begin{corollary}\label{WBC1}
		Let $A$ be a unital simple separable infinite-dimensional stably finite C*-algebra with $\beta$-comparison. Let $\alpha\colon G\to \mathrm{Aut}(A)$ be an action of a second-countable compact group with the weak tracial Rokhlin property with comparison. Then the crossed product $A\rtimes_\alpha G$ has $\beta$-comparison.
		\begin{proof}
			By Theorem \ref{WBC}, the algebra $A^\alpha$ has $\beta$-comparison. Theorem \ref{MESI} implies that $A\rtimes_\alpha G$ is stably isomorphic to $A^\alpha$. Since stably isomorphic C*-algebras have canonically isomorphic Cuntz semigroups, $A\rtimes_\alpha G$ has $\beta$-comparison.
		\end{proof}
	\end{corollary}
	
	\begin{theorem}\label{WNC}
		Let $A$ be a unital simple separable infinite-dimensional stably finite C*-algebra with Winter's $n$-comparison. Let $\alpha\colon G\to \mathrm{Aut}(A)$ be an action of a second-countable compact group with the weak tracial Rokhlin property with comparison. Then the fixed point algebra $A^\alpha$ has Winter's $n$-comparison.
		\begin{proof}
			The proof is essentially the same as that of \cite[Theorem 4.6]{TF2025}, so we omit it.
		\end{proof}
	\end{theorem}
	
	\begin{corollary}\label{WNC1}
		Let $A$ be a unital simple separable infinite-dimensional stably finite C*-algebra with Winter's $n$-comparison. Let $\alpha\colon G\to \mathrm{Aut}(A)$ be an action of a second-countable compact group with the weak tracial Rokhlin property with comparison. Then the crossed product $A\rtimes_\alpha G$ has Winter's $n$-comparison.
		\begin{proof}
			By Theorem \ref{WNC}, the algebra $A^\alpha$ has Winter's $n$-comparison. Theorem \ref{MESI} implies that $A\rtimes_\alpha G$ is stably isomorphic to $A^\alpha$. Since stably isomorphic C*-algebras have canonically isomorphic Cuntz semigroups, $A\rtimes_\alpha G$ has Winter's $n$-comparison.
		\end{proof}
	\end{corollary}
	
	\section{The weak tracial Rokhlin property}
	\label{sec8.4}
	As discussed after Definition \ref{WTRC}, the weak tracial Rokhlin property with comparison defined in Definition \ref{WTRC} is stronger than the weak tracial Rokhlin property for finite groups (Definition \ref{WTRF}). In this section, we discuss the differences between them. Definition \ref{NWTR} is presented solely for discussion purposes. It is a straightforward translation of the definition of the weak tracial Rokhlin property for finite groups.
	
	\begin{definition}\label{NWTR}
		Let $G$ be a second-countable compact group, let $A$ be a unital simple infinite-dimensional C*-algebra, and let $\alpha\colon G\to \mathrm{Aut}(A)$ be an action. We say that the action $\alpha$ has the naive weak tracial Rokhlin property if for any $\varepsilon>0$, any finite set $F\subseteq A$, any finite set $S\subseteq C(G)$, and any $x\in A_+$ with $\|x\|=1$, there exist a positive contraction $d\in A^\alpha$ and a contractive completely positive map $\psi\colon C(G)\to \overline{dAd}$ such that the following hold:
		\begin{enumerate}[{(}1{)}]
			\item $\psi$ is an $(S,F,\varepsilon)$-approximately central equivariant strongly order zero map.
			\item $\psi(1)=d$.
			\item $(1-d-\varepsilon)_+\precsim_A x$.
			\item $\|dxd\|>1-\varepsilon$.
		\end{enumerate}
	\end{definition}
	
	\begin{remark}
		Similar to condition (4) in \cite[Definition 4.1]{Phi2014}, since we consider the transfer of properties from $A$ to $A^\alpha$, the Cuntz comparison condition in $A^\alpha$ is necessary to handle certain technical details. Thus, the condition $(1-d-\varepsilon)_+\precsim_{A^\alpha}y$ appears essential; without it, we may not be able to complete most of the previous proofs.
	\end{remark}
	
	We now show that condition (\ref{d8.3.4}) of Definition \ref{WTRC} is automatically satisfied for actions of finite groups.
	
	\begin{proposition}\label{CIFA}
		Let $A$ be a unital simple infinite-dimensional C*-algebra. Let $\alpha\colon G\to \mathrm{Aut}(A)$ be a finite group action with the weak tracial Rokhlin property. Then for any finite set $F\subseteq A$, any $\varepsilon>0$, any $x\in A_+$ with $\|x\|=1$, and any $y\in(A^\alpha)_+\setminus\{0\}$, there exist a positive contraction $d\in A^\alpha$ and positive contractions $(d_g)_{g\in G}$ such that the following hold:
		\begin{enumerate}[{(}1{)}]
			\item $d=\sum_{g\in G}d_g$.\label{p8.8.1}
			\item $\|d_ga-ad_g\|<\varepsilon$ for all $a\in F$ and all $g\in G$.\label{p8.8.2}
			\item $\|d_gd_h\|<\varepsilon$ for all $g,h\in G$ with $g\neq h$.\label{p8.8.3}
			\item $\|\alpha_g(d_h)-d_{gh}\|<\varepsilon$ for all $g,h\in G$.\label{p8.8.4}
			\item $(1-d-\varepsilon)_+\precsim_Ax$.\label{p8.8.5}
			\item $(1-d-\varepsilon)_+\precsim_{A^\alpha}y$.\label{p8.8.6}
			\item $\|dxd\|>1-\varepsilon$. \label{p8.8.7}
		\end{enumerate}
		\begin{proof}
			Let $\varepsilon>0$, let $F\subseteq A$ be a finite set, let $x\in A_+$ with $\|x\|=1$, and let $y\in(A^\alpha)_+\setminus\{0\}$. Since $A^\alpha$ is simple and non-type I, by \cite[Lemma 2.1]{Phi2014}, there exists a positive element $z\in\overline{yA^\alpha y}$ such that 0 is a limit point of $\mathrm{sp}(z)$. Define continuous functions $h,h_0\colon [0,1]\rightarrow[0,1]$ by
			\begin{equation*}
				h(\lambda)=
				\begin{cases}
					0,&\mbox{$\lambda\in[0,1-\frac{\varepsilon}{2}]$,}\\
					\frac{2}{\varepsilon}(\lambda-1)+1,&\mbox{$\lambda\in[1-\frac{\varepsilon}{2},1]$.}
				\end{cases}
			\end{equation*}
			and
			\begin{equation*}
				h_0(\lambda)=
				\begin{cases}
					(1-\frac{\varepsilon}{2})^{-1}\lambda,&\mbox{$\lambda\in[0,1-\frac{\varepsilon}{2}]$,}\\
					1,&\mbox{$\lambda\in[1-\frac{\varepsilon}{2},1]$.}
				\end{cases}
			\end{equation*}
			Then we have $\|x-h_0(x)\|\leq\frac{\varepsilon}{2}$. Moreover, since $\|x\|=1$, we have $h(x)\neq0$. By \cite[Lemma 2.6]{Phi2014}, we can choose a nonzero positive element $x_0\in\overline{h(x)Ah(x)}$ such that $x_0\precsim_A z$. We may assume $\|x_0\|=1$.
			
			Apply Lemma \ref{WTRFE} to $\alpha$ for given $x_0$, $\frac{\varepsilon}{2}$, and $F$. We obtain positive contractions $(d_g)_{g\in G}\in A$ and $d\in A^\alpha$ such that the following hold:
			\begin{enumerate}[{(}1{)}]
				\setcounter{enumi}{7}
				\item $d=\sum_{g\in G}d_g$.
				\item $\|d_ga-ad_g\|<\frac{\varepsilon}{2}$ for all $a\in F$ and all $g\in G$.
				\item $\|d_gd_h\|<\frac{\varepsilon}{2}$ for all $g,h\in G$ with $g\neq h$.
				\item $\|\alpha_g(d_h)-d_{gh}\|<\frac{\varepsilon}{2}$ for all $g,h\in G$.
				\item $(1-d-\frac{\varepsilon}{2})_+\precsim_Ax_0$.
				\item $\|dx_0d\|>1-\frac{\varepsilon}{2}$.
			\end{enumerate}
			Clearly, conditions (\ref{p8.8.1}), (\ref{p8.8.2}), (\ref{p8.8.3}), and (\ref{p8.8.4}) hold. Moreover, $(1-d-\varepsilon)_+\precsim_A(1-d-\frac{\varepsilon}{2})_+\precsim_Ax_0\precsim_Ax$, which is condition (\ref{p8.8.5}). Furthermore,
			\[(1-d-\varepsilon)_+\precsim_Ax_0\precsim_Az,\ (1-d-\varepsilon)_+,z\in A^\alpha,\ \mbox{and}\ 0\in\overline{{\textrm{sp}}(z)\setminus\{0\}}.\]
			By \cite[Lemma 3.7]{AGP2021}, $(1-d-\varepsilon)_+\precsim_{A^\alpha}z$. Since $z\in\overline{yA^\alpha y}$, we have $(1-d-\varepsilon)_+\precsim_{A^\alpha}y$, which is condition (\ref{p8.8.6}).\\
			It remains to prove condition (\ref{p8.8.7}). Since $h_0(x)h(x)=h(x)$, we have
			\[x_0=h_0(x)^\frac{1}{2}x_0h_0(x)\frac{1}{2}\leq h_0(x).\]
			Therefore, together with $\|dx_0d\|>1-\frac{\varepsilon}{2}$, we have
			\[\|dxd\|\geq\|dh_0(x)d\|-\|h_0(x)-x\|\geq\|dx_0d\|-\frac{\varepsilon}{2}>1-\frac{\varepsilon}{2}-\frac{\varepsilon}{2}=1-\varepsilon.\]		
		\end{proof}
	\end{proposition}
	
	We now present some situations where condition (\ref{d8.3.5}) of Definition \ref{WTRC} automatically holds for finite group actions.
	
	\begin{proposition}\label{RCL1}
		Let $A$ be a unital simple infinite-dimensional stably finite C*-algebra. Let $\alpha\colon G\to \mathrm{Aut}(A)$ be a finite group action with the weak tracial Rokhlin property. If $\mathrm{rc}(A)<1$, then $\alpha$ has the weak tracial Rokhlin property with comparison.
		
		In particular, if $A$ has strict comparison, then $\alpha$ has the weak tracial Rokhlin property with comparison.
		\begin{proof}
			Let $\varepsilon>0$, let $F\subseteq A$ be a finite set, let $x\in A_+$ with $\|x\|=1$, and let $y\in(A^\alpha)_+\setminus\{0\}$. Since $A^\alpha$ is simple and non-type I, by \cite[Lemma 2.1]{Phi2014}, there exists a positive element $z\in\overline{yA^\alpha y}$ such that 0 is a limit point of $\mathrm{sp}(z)$.
			
			Choose $n\in\mathbb{N}$ such that
			\begin{equation}\label{p8.9.eq1}
				n>\frac{2}{1-{\textrm{rc}}(A)}.
			\end{equation}
			Since $A^\alpha$ is simple and non-type I, by Lemma \ref{OCP}, there exist $n+1$ mutually orthogonal, Cuntz equivalent nonzero positive elements in $A^\alpha$. Let $e$ be one of them. Thus, by \cite[Lemma 2.1]{Phi2014}, there exists $r\in(\overline{eA^\alpha e})_+\setminus\{0\}$ such that 0 is a limit point of $\mathrm{sp}(r)$. Define a continuous function $h\colon [0,1]\to [0,1]$ by
			\begin{equation*}
				h(\lambda)=
				\begin{cases}
					0,&\mbox{$\lambda\in[0,1-\frac{\varepsilon}{2}]$,}\\
					\frac{2}{\varepsilon}(\lambda-1)+1,&\mbox{$\lambda\in[1-\frac{\varepsilon}{2},1]$.}
				\end{cases}
			\end{equation*}
			Apply \cite[Lemma 2.6]{Phi2014} to find $x_0\in\overline{h(x)Ah(x)}$ such that
			\[\|x_0\|=1,\ x_0\precsim_Ar,\ \mbox{and}\ x_0\precsim_Ay_0.\]
			Apply Lemma \ref{WTRFE} to $\alpha$ for given $x_0$, $\frac{\varepsilon}{2}$, and $F$. We obtain positive contractions $(d_g)_{g\in G}\in A$ and $d\in A^\alpha$ such that the following hold:
			\begin{enumerate}
				\item $d=\sum_{g\in G}d_g$.
				\item $\|d_ga-ad_g\|<\frac{\varepsilon}{2}$ for all $a\in F$ and all $g\in G$.
				\item $\|d_gd_h\|<\frac{\varepsilon}{2}$ for all $g,h\in G$ with $g\neq h$.
				\item $\|\alpha_g(d_h)-d_{gh}\|<\frac{\varepsilon}{2}$ for all $g,h\in G$.
				\item $(1-d-\frac{\varepsilon}{2})_+\precsim_Ax_0$.
				\item $\|dx_0d\|>1-\frac{\varepsilon}{2}$. 
			\end{enumerate}
			By Proposition \ref{CIFA}, we only need to show that $(1-d-\varepsilon)_+\precsim_{A^\alpha}d$. Since $(1-d-\varepsilon)_+\precsim_Ax_0\precsim_A r$, by \cite[Lemma 3.7]{AGP2021}, we have $(1-d-\varepsilon)_+\precsim_{A^\alpha}r$. Thus $(1-d-\varepsilon)_+\precsim_{A^\alpha}e$. Let $\tau\in \mathrm{QT}(A^\alpha)$. Then $d_\tau((1-d-\varepsilon)_+)\leq d_\tau(e)$, and by Lemma \ref{CPA},
			\[1\sim_{A^\alpha}(1-2\varepsilon)_+\precsim_{A^\alpha}(d-\varepsilon)_++(1-d-\varepsilon)_+,\]
			so we have
			\[d_\tau(d)\geq d_\tau((d-\varepsilon)_+)\geq1-d_\tau((1-d-\varepsilon)_+)\geq1-d_\tau(e).\] 
			The construction of $e$ ensures that $d_\tau(e)\leq\frac{1}{n+1}<\frac{1}{n}$,
			hence
			\begin{align}\label{p8.9.eq2}
				d_\tau((1-d-\varepsilon)_+)+1-\frac{2}{n}\leq &{}d_\tau(e)+1-\frac{2}{n}<1-\frac{1}{n}<1-d_\tau(e)\leq d_\tau(d). 
			\end{align}
			Using \cite[Theorem 4.1]{AGP2021} in the first step and (\ref{p8.9.eq1}) in the second step, we have $\mathrm{rc}(A^\alpha)\leq\mathrm{rc}(A)<1-\frac{2}{n}$. Since (\ref{p8.9.eq2}) holds for all $\tau\in \mathrm{QT}(A^\alpha)$, we have $(1-d-\varepsilon)_+\precsim_{A^\alpha}d$.
		\end{proof}
	\end{proposition}
	
	The condition $(1-d-\varepsilon)_+\precsim_{A^\alpha} d$ is only used to prove that the crossed product $A\rtimes_\alpha G$ and the fixed point algebra $A^\alpha$ are Morita equivalent (see Corollary \ref{MESI}). Without it, we cannot transfer structural properties from the original algebra to the crossed product. However, if we omit this condition, other conclusions still hold. In other words, without $(1-d-\varepsilon)_+\precsim_{A^\alpha} d$, we can still transfer structural properties from the original algebra to the fixed point algebra. Therefore, we consider the following version of the weak tracial Rokhlin property for actions of compact groups.
	
	\begin{definition}\label{WTR}
		Let $G$ be a second-countable compact group, let $A$ be a unital simple infinite-dimensional C*-algebra, and let $\alpha\colon G\to \mathrm{Aut}(A)$ be an action. We say that the action $\alpha$ has the weak tracial Rokhlin property if for any $\varepsilon>0$, any finite set $F\subseteq A$, any finite set $S\subseteq C(G)$, any $x\in A_+$ with $\|x\|=1$, and any $y\in(A^\alpha)_+\setminus\{0\}$, there exist a positive contraction $d\in A^\alpha$ and a contractive completely positive map $\psi\colon C(G)\to \overline{dAd}$ such that the following hold:
		\begin{enumerate}
			\item $\psi$ is an $(S,F,\varepsilon)$-approximately central equivariant strongly order zero map.
			\item $\psi(1)=d$.
			\item $(1-d-\varepsilon)_+\precsim_A x$.
			\item $(1-d-\varepsilon)_+\precsim_{A^\alpha} y$.
			\item $\|dxd\|>1-\varepsilon$.
		\end{enumerate}
	\end{definition}
	
	Since $(1-d-\varepsilon)_+\precsim_{A^\alpha}d$ is only used to prove that $\alpha$ is saturated, by Remark \ref{MTFR}, the following theorems still hold for actions with the weak tracial Rokhlin property.
	
	\begin{proposition}
		Let $A,B$ be unital simple infinite-dimensional C*-algebras, let $G$ be a second-countable compact group, and let $\alpha\colon G\to \mathrm{Aut}(A)$ be a group action with the weak tracial Rokhlin property and $\beta\colon G\to \mathrm{Aut}(B)$ be a trivial group action. Then $\alpha\otimes\beta\colon G\to \mathrm{Aut}(A\otimes_{min}B)$ has the weak tracial Rokhlin property.
	\end{proposition}
	
	\begin{theorem}
		Let $A$ be a unital simple separable infinite-dimensional purely infinite C*-algebra. Let $\alpha\colon G\to \mathrm{Aut}(A)$ be an action of a second-countable compact group with the weak tracial Rokhlin property. Then the fixed point algebra $A^\alpha$ is purely infinite.
	\end{theorem}
	
	\begin{theorem}
		Let $A$ be a unital simple separable infinite-dimensional C*-algebra that is tracially $\mathcal{Z}$-stable. Let $\alpha\colon G\to \mathrm{Aut}(A)$ be an action of a second-countable compact group with the weak tracial Rokhlin property. Then the fixed point algebra $A^\alpha$ is tracially $\mathcal{Z}$-stable.
	\end{theorem}
	
	\begin{theorem}
		Let $A$ be a unital simple separable infinite-dimensional amenable C*-algebra that is $\mathcal{Z}$-stable. Let $\alpha\colon G\to \mathrm{Aut}(A)$ be an action of a second-countable compact group with the weak tracial Rokhlin property. Then the fixed point algebra $A^\alpha$ is $\mathcal{Z}$-stable.
	\end{theorem}
	
	\begin{theorem}
		Let $A$ be a unital simple separable infinite-dimensional C*-algebra. Let $\alpha\colon G\to \mathrm{Aut}(A)$ be an action of a second-countable compact group with the weak tracial Rokhlin property. Let $n\in\mathbb{N}$, and let $a,b\in M_n(A)_+$. Suppose 0 is a limit point of $\mathrm{sp}(b)$. Then $a\precsim_A b$ if and only if $a\precsim_{A^\alpha} b$.
	\end{theorem}
	
	\begin{theorem}
		Let $A$ be a unital simple separable infinite-dimensional stably finite C*-algebra. Let $\alpha\colon G\to \mathrm{Aut}(A)$ be an action of a second-countable compact group with the weak tracial Rokhlin property. Then $\mathrm{rc}(A^\alpha)\leq \mathrm{rc}(A)$.
	\end{theorem}
	
	\begin{corollary}
		Let $A$ be a unital simple separable infinite-dimensional stably finite C*-algebra with strict comparison. Let $\alpha\colon G\to \mathrm{Aut}(A)$ be an action of a second-countable compact group with the weak tracial Rokhlin property. Then $A^\alpha$ has strict comparison.
	\end{corollary}
	
		\begin{proposition}
		Let $A$ be a unital simple separable infinite-dimensional stably finite C*-algebra. Let $\alpha\colon G\to \mathrm{Aut}(A)$ be an action of a second-countable compact group with the weak tracial Rokhlin property. Let $\iota\colon A^\alpha\to A$ be the inclusion map. Then
		\begin{enumerate}
			\item the map $\mathrm{W}(\iota)\colon \mathrm{W}(A^\alpha)\to \mathrm{W}(A)$ induces an injective map of ordered semigroups from $\mathrm{W}(A^\alpha)_+\cup\{0\}$ to $\mathrm{W}(A)_+\cup\{0\}$;
			
			\item the map $\mathrm{Cu}(\iota)\colon \mathrm{Cu}(A^\alpha)\to \mathrm{Cu}(A)$ induces an injective map of ordered semigroups from $\mathrm{Cu}(A^\alpha)_+\cup\{0\}$ to $\mathrm{Cu}(A)_+\cup\{0\}$.
		\end{enumerate}
	\end{proposition}
	
	\begin{proposition}
		Let $A$ be a unital, simple, separable, infinite-dimensional C*-algebra. Let $\alpha\colon G\to \mathrm{Aut}(A)$ be an action of a second-countable compact group with the weak tracial Rokhlin property. Then the restriction map $\mathrm{QT}(A)\to \mathrm{QT}({A^\alpha})$ is surjective.
	\end{proposition}
	
	\begin{theorem}
		Let $A$ be a unital simple separable infinite-dimensional stably finite C*-algebra with $\beta$-comparison. Let $\alpha\colon G\to \mathrm{Aut}(A)$ be an action of a second-countable compact group with the weak tracial Rokhlin property. Then the fixed point algebra $A^\alpha$ has $\beta$-comparison.
	\end{theorem}
	
	\begin{theorem}
		Let $A$ be a unital simple separable infinite-dimensional stably finite C*-algebra with Winter's $n$-comparison. Let $\alpha\colon G\to \mathrm{Aut}(A)$ be an action of a second-countable compact group with the weak tracial Rokhlin property. Then the fixed point algebra $A^\alpha$ has Winter's $n$-comparison.
	\end{theorem}
	
	Similar to the case of finite groups, condition (\ref{d8.3.5}) of Definition \ref{WTRC} is automatically satisfied for compact groups when $\mathrm{rc}(A)<1$.
	
	\begin{proposition}\label{RCL1C}
		Let $A$ be a unital simple separable infinite-dimensional stably finite C*-algebra. Let $\alpha\colon G\to \mathrm{Aut}(A)$ be an action of a second-countable compact group with the weak tracial Rokhlin property. If $\mathrm{rc}(A)<1$, then $\alpha$ has the weak tracial Rokhlin property with comparison.
		
		In particular, if $A$ has strict comparison, then $\alpha$ has the weak tracial Rokhlin property with comparison.
		\begin{proof}
			The proof is essentially the same as that of Proposition \ref{RCL1}, so we omit it.
		\end{proof}
	\end{proposition}
	
	Since tracial $\mathcal{Z}$-stability or $\mathcal{Z}$-stability implies strict comparison, we have the following corollaries.
	
	\begin{corollary}
		Let $A$ be a unital simple separable infinite-dimensional C*-algebra that is tracially $\mathcal{Z}$-stable. Let $\alpha\colon G\to \mathrm{Aut}(A)$ be an action of a second-countable compact group with the weak tracial Rokhlin property. Then the crossed product $A\rtimes_\alpha G$ is tracially $\mathcal{Z}$-stable.
	\end{corollary}
	
	\begin{corollary}
		Let $A$ be a unital simple separable infinite-dimensional amenable C*-algebra that is $\mathcal{Z}$-stable. Let $\alpha\colon G\to \mathrm{Aut}(A)$ be an action of a second-countable compact group with the weak tracial Rokhlin property. Then both the fixed point algebra $A^\alpha$ and the crossed product $A\rtimes_\alpha G$ are $\mathcal{Z}$-stable.
	\end{corollary}
	
	\begin{corollary}
		Let $A$ be a unital simple separable infinite-dimensional stably finite C*-algebra with strict comparison. Let $\alpha\colon G\to \mathrm{Aut}(A)$ be an action of a second-countable compact group with the weak tracial Rokhlin property. Then $A\rtimes_\alpha G$ has strict comparison.
	\end{corollary}
	
	We now prove that condition (\ref{d8.3.5}) of Definition \ref{WTRC} is automatically satisfied for compact groups when $A$ is purely infinite.
	
	\begin{proposition}\label{PIC}
		Let $A$ be a unital simple separable infinite-dimensional C*-algebra. Let $\alpha\colon G\to \mathrm{Aut}(A)$ be an action of a second-countable compact group with the weak tracial Rokhlin property. If $A$ is purely infinite, then $\alpha$ has the weak tracial Rokhlin property with comparison.
		\begin{proof}
			As in the proof of Theorem \ref{FAPI}, we can show that $A^\alpha$ is purely infinite, so $(1-d-\varepsilon)_+\precsim_{A^\alpha}d$ is automatically satisfied.
		\end{proof}
	\end{proposition}
	
	\begin{corollary}
		Let $A$ be a unital simple separable infinite-dimensional purely infinite C*-algebra. Let $\alpha\colon G\to \mathrm{Aut}(A)$ be an action of a second-countable compact group with the weak tracial Rokhlin property. Then $A\rtimes_\alpha G$ is purely infinite.
	\end{corollary}
	
	\section{A totally disconnected compact group action}
	\label{sec8.5}
	In this section, we construct an action of a totally disconnected infinite compact group on the Jiang-Su algebra $\mathcal{Z}$, which has the weak tracial Rokhlin property with comparison, but does not have the tracial Rokhlin property.
	
	Let $S_2$ denote the symmetric group of order $2$ (i.e., the permutation group on two elements). Define the group $G=\prod_{n=1}^\infty S_2$ (i.e., the infinite direct product of $S_2$). Its action is the infinite tensor product of the permutation actions on the twofold tensor products of $\mathcal{Z}$.
	\begin{notation}
		Let $B=\mathcal{Z}\otimes\mathcal{Z}$, $\gamma\colon S_2\to \mathrm{Aut}(B)$ be the permutation action.
		
		Define $G=\prod_{n=1}^\infty S_2$ and $A=\bigotimes_{n=1}^\infty B$. Let $\alpha\colon G\to \mathrm{Aut}(A)$ be the infinite tensor action, defined by
		\[\alpha_{(h_1,h_2,\dots)}(b_1\otimes b_2\otimes\cdots)=\gamma_{h_1}(b_1)\otimes\gamma_{h_2}(b_2)\otimes\cdots\]
		for $h_1,h_2,\dots\in S_2$ and $b_1,b_2,\dots\in B$ with $b_n=1$ for all but finitely many $n\in\mathbb{N}$.
	\end{notation}
	
	\begin{remark}
		By \cite[Theorem 1.10.16]{Lin2001}, $A\cong \mathcal{Z}$.
	\end{remark}
	
	\begin{notation}
		For $n\in\mathbb{N}$, let $B_n=B$, so $A=\bigotimes_{m=1}^\infty B_m$, and let $A_n=\bigotimes_{m=1}^n B_m$, so $A=\lim\limits_{\longrightarrow n}A_n$. Let $H_n=S_2$, so $G=\prod_{m=1}^\infty H_m$, and let $G_n=\prod_{m=1}^n H_m$, then $G=\lim\limits_{\longleftarrow n}G_n$. This gives
		\[C(G_n)=\bigotimes_{m=1}^nC(H_m),\quad \mbox{and}\quad C(G)=\lim_{\longrightarrow n}C(G_n)=\bigotimes_{m=1}^\infty C(H_m).\]
		We identify $C(G_n)$ with its image in $C(G)$. Denote the action of $G_n$ on $A_n$ by $\alpha_n$.
	\end{notation}
	
	\begin{lemma}\label{WTRZ}
		Let $A=\mathcal{Z}^{\otimes n}\cong\mathcal{Z}$, let $G=S_n$ be the permutation group on the set $\{1,2,\dots,n\}$. Let $\alpha\colon G\to \mathrm{Aut}(A)$ be the action defined on elementary tensors by
		\[\alpha_\sigma(z_1\otimes\dots\otimes z_n)=z_{\sigma^{-1}(1)}\otimes\dots\otimes z_{\sigma^{-1}(n)}.\]
		Then for every $\varepsilon>0$, every finite set $F\subseteq A$, and every $x\in A_+$ with $\|x\|=1$, there exist mutually orthogonal positive contractions $(d_g)_{g\in G}\in A$ satisfying $d=\sum_{g\in G}d_g\in A^\alpha$, such that the following hold:
		\begin{enumerate}[{(}1{)}]
			\item $\|d_ga-ad_g\|<\varepsilon$ for all $a\in F$ and all $g\in G$.
			
			\item $\alpha_g(d_h)=d_{gh}$ for all $g,h\in G$.
			
			\item $1-d\precsim_Ax$.
			
			\item $1-d\precsim_{A^\alpha}y$.
			
			\item $\|dxd\|>1-\varepsilon$.
		\end{enumerate}
		\begin{proof}
			This follows from \cite[Example 5.10]{HO2013} (see also \cite[Theorem 3.10]{AGJP2021}) and Proposition \ref{CIFA}.
		\end{proof}
	\end{lemma}
	
	\begin{lemma}
		Let $n\in\mathbb{N}$, let $A_1,A_2,\dots,A_n$ be unital C*-algebras, and for $m=1,2,\dots,n$ let $e_m\in A_m$ be a positive contraction and let $\tau_m$ be a tracial state on $A_m$. Let $A=A_1\otimes A_2\otimes\dots\otimes A_n$ (minimal tensor product), and let
		\[e=e_1\otimes e_2\otimes\dots\otimes e_n\in A\quad \mbox{and}\quad \tau=\tau_1\otimes\tau_2\otimes\dots\otimes\tau_n\in \mathrm{T}(A).\]
		Then
		\[d_\tau(e)=\prod_{m=1}^nd_{\tau_m}(e_m).\]
		\begin{proof}
			\begin{align*}
				d_\tau(e)=\lim_{k\to\infty}\tau(e^\frac{1}{k})=&{}\lim_{k\to\infty}\tau((e_1\otimes e_2\otimes\dots\otimes e_n)^\frac{1}{k})\\
				=&{}\lim_{k\to\infty}\tau(e_1^\frac{1}{k}\otimes e_2^\frac{1}{k}\otimes\dots\otimes e_n^\frac{1}{k})\\
				=&{}\lim_{k\to\infty}\tau_1(e_1^\frac{1}{k})\tau_2(e_2^\frac{1}{k})\cdots\tau_n(e_n^\frac{1}{k})\\
				=&{}\prod_{m=1}^nd_{\tau_m}(e_m).
			\end{align*}
		\end{proof}
	\end{lemma}
	
	\begin{lemma}\label{TPC}
		Let $n\in\mathbb{N}$, let $A_1,A_2,\dots,A_n$ be unital C*-algebras, and for $m=1,2,\dots,n$ let $e_m\in A_m$ be a positive contraction and let $\tau_m$ be a tracial state on $A_m$. Let $A=A_1\otimes A_2\otimes\dots\otimes A_n$ (minimal tensor product), and let
		\[e=e_1\otimes e_2\otimes\dots\otimes e_n\in A\quad \mbox{and}\quad \tau=\tau_1\otimes\tau_2\otimes\dots\otimes\tau_n\in \mathrm{T}(A).\]
		Then
		\[d_\tau(1-e)\leq\sum_{m=1}^nd_{\tau_m}(1-e_m).\]
		\begin{proof}
			We prove this by induction on $n$.
			
			The case $n=1$ is obvious.
			
			For $n=2$, since $1-e_1\otimes e_2\leq(1-e_1\otimes1)+(1-1\otimes e_2)$, we have
			\begin{align*}
				d_\tau(1-e_1\otimes e_2)\leq &{}d_\tau((1-e_1\otimes1)+(1-1\otimes e_2))\\
				\leq&{}d_\tau(1-e_1\otimes1)+d_\tau(1-1\otimes e_2)\\ =&{}d_{\tau_1}(1-e_1)d_{\tau_2}(1)+d_{\tau_1}(1)d_{\tau_2}(1-e_2)\\
				=&{}d_{\tau_1}(1-e_1)+d_{\tau_2}(1-e_2).
			\end{align*}
			So the case $n=2$ holds.
			
			Now suppose the statement holds for some $n\geq2$. Then for $n+1$ we have
			\begin{align*}
				&{}d_\tau(1-e_1\otimes e_2\otimes\dots\otimes e_n\otimes e_{n+1})\\
				\leq&{} d_\tau(1-e_1\otimes e_2\otimes\dots\otimes e_n)+d_{\tau_{n+1}}(1-e_{n+1})\\
				\leq&{}\sum_{m=1}^nd_{\tau_m}(1-e_m)+d_{\tau_{n+1}}(1-e_{n+1})\\
				=&{}\sum_{m=1}^{n+1}d_{\tau_m}(1-e_m).
			\end{align*}
		\end{proof}
	\end{lemma}
	
	\begin{theorem}
		The action $\alpha\colon G\to \mathrm{Aut}(A)$ has the weak tracial Rokhlin property with comparison.
		\begin{proof}
			Let $F\subseteq A$ and $S\subseteq C(G)$ be finite sets, let $\varepsilon>0$, let $x\in A_+$ with $\|x\|=1$, and let $y\in (A^\alpha)_+\setminus\{0\}$. Without loss of generality, we may assume $\|a\|\leq1$ for all $a\in F$, $\|f\|\leq1$ for all $f\in S$, and $\varepsilon<1$. We need to find a positive contraction $d\in A^\alpha$ and a completely positive contractive map $\psi\colon C(G)\to \overline{dAd}$ such that the following hold:
			\begin{enumerate}[{(}1{)}]
				\item $\psi$ is an $(S,F,\varepsilon)$-approximately central equivariant strong order zero map.\label{t8.19.1}
				\item $\psi(1)=d$.\label{t8.19.2}
				\item $(1-d-\varepsilon)_+\precsim_A x$.\label{t8.19.3}
				\item $(1-d-\varepsilon)_+\precsim_{A^\alpha} y$.\label{t8.19.4}
				\item $(1-d-\varepsilon)_+\precsim_{A^\alpha} d$.\label{t8.19.5}
				\item $\|dxd\|>1-\varepsilon$.\label{t8.19.6}
			\end{enumerate}
			
			Since $A\cong\mathcal{Z}$, we denote $\delta_1=d_\tau(x)>0$ where $\tau$ is the unique tracial state on $A$.
			
			By \cite[Example 5.10]{HO2013}, the action $\gamma$ has the weak tracial Rokhlin property. Therefore, by \cite[Lemma 3.3]{FG2020}, $B\rtimes_\gamma S_2$ is simple, and hence by Theorem \ref{CI}, $B^\gamma$ is simple. It is easy to verify that $A^\alpha$ can be identified with $\bigotimes_{m=1}^\infty B_m^\gamma$ (see \cite[Lemma 4.5]{MP2025}), hence it is simple. Meanwhile, since $B$ has strict comparison, $B^\gamma$ also has it, and thus $A^\alpha$ does as well (see \cite[Proposition 5.16.1(4)]{FHLRTVW2021}). Since $A^\alpha$ is unital, stably finite, and amenable, $\mathrm{T}(A^\alpha)$ is a nonempty compact set. We denote $\delta_2=\inf_{\tau\in \mathrm{T}(A^\alpha)}d_\tau(y)>0$.
			
			Since $\cup_{n=1}^\infty A_n$ is dense in $A$, there exists $N_1\in\mathbb{N}$ and a finite set $F_0\subseteq A_{N_1}\subseteq A$ such that for each $a\in F$, there exists $b\in F_0$ satisfying $\|a-b\|<\frac{\varepsilon}{8}$, and for all $b\in F_0$, $\|b\|\leq1$. We may assume there exist $c_{ij}\in B$ for $i=1,2,\dots,N_1$ and $j=1,2,\dots,k$ such that $F_0=\{c_{1j}\otimes c_{2j}\otimes\dots\otimes c_{N_1j}\colon 1\leq j\leq k\}$. Similarly, there exists $N_2\in\mathbb{N}$ and a finite set $S_0\subseteq C(G_{N_2})\subseteq C(G)$ such that for each $f\in F$, there exists $e\in S_0$ satisfying $\|f-e\|<\frac{\varepsilon}{4}$, and for all $e\in S_0$, $\|e\|\leq1$. 
			Let $N=\max(N_1,N_2)$. Let $M=\mathrm{card}(G_N)$, and let $\varepsilon_0=\min(\frac{\delta_1}{N},\frac{\delta_2}{N},\frac{1}{2N},\frac{\varepsilon}{4MN})$.
			
			Since the action $\gamma$ has the weak tracial Rokhlin property, by Lemma \ref{WTRZ}, for each $1\leq i\leq N$, set $F_i=\{c_{ij}\colon 1\leq j\leq k\}$, we obtain a positive contraction $d^{(i)}\in B^\gamma$ and mutually orthogonal positive contractions $(d_{g}^{(i)})_{g\in H_i}\in B$ such that the following hold:
			\begin{itemize}
				\item[(7)] $d^{(i)}=\sum_{g\in H_i}d_{g}^{(i)}$.
				
				\item[(8)] For all $a\in F_i$ and all $g\in H_i$, $\|d_{g}^{(i)}a-ad_{g}^{(i)}\|<\varepsilon_0$.
				
				\item[(9)] For all $g,h\in H_i$, $\gamma_g(d_{h}^{(i)})=d_{gh}^{(i)}$.
				
				\item[(10)] $d_\tau(1-d^{(i)})<\varepsilon_0$ where $\tau$ is the unique tracial state on $B$.
				
				\item[(11)] For all $\tau\in \mathrm{T}(B^\gamma)$, $d_\tau(1-d^{(i)})<\varepsilon_0$.
			\end{itemize}
			For $h=(h_1,h_2,\dots,h_N)\in G_N=\prod_{m=1}^N H_m$, let
			\[d_h=d_{h_1}^{(1)}\otimes d_{h_2}^{(2)}\otimes\dots\otimes d_{h_N}^{(N)}.\]
			It is easy to verify that for $g,h\in G_N$ with $g\neq h$, $d_gd_h=0$; for $g,h\in G_N$, $\alpha_{N_g}(d_h)=d_{gh}$; and for $a\in F_0$ and $g\in G_N$, $\|ad_g-d_ga\|<N\varepsilon_0$. Define $d=\sum_{h\in G_N}d_h\in A_N\subseteq A$. There exists an equivariant completely positive contractive order zero map $\psi_0\colon C(G_N)\to \overline{dA_Nd}\subseteq\overline{dAd}$ given by $\psi_0(f)=\sum_{h\in G_N}f(h)d_h$ for $f\in C(G_N)$, such that $d=\psi_0(1)$.\\
			For $f\in S_0$, we have
			\begin{align*}
				&{}\sup_{g\in G_N}\|\alpha_{N_g}(\psi_0(f))-\psi_0(\mathrm{Lt}_g(f))\|\\
				=&{}\sup_{g\in G_N}\|\alpha_{N_g}(\sum_{h\in G_N}f(h)d_h)-\sum_{h\in G_N}f(g^{-1}h)d_h\|\\
				=&{}\sup_{g\in G_N}\|\sum_{h\in G_N}f(h)\alpha_{N_g}(d_h)-\sum_{h\in G_N}f(h)d_{gh}\|\\
				=&{}0.
			\end{align*}
			For $f\in S_0$ and $a\in F_0$, we have
			\begin{align*}
				\|a\psi_0(f)-\psi_0(f)a\|=&{}\|a\sum_{h\in G_N}f(h)d_h-\sum_{h\in G_N}f(h)d_ha\|\\
				=&{}\|\sum_{h\in G_N}f(h)(ad_h-d_ha)\|\\
				\leq&{}MN\varepsilon_0<\frac{\varepsilon}{4}.
			\end{align*}
			For $f_1,f_2\in S_0$, we have
			\begin{align*}
				&{}\|\psi_0(1)\psi_0(f_1f_2)-\psi_0(f_1)\psi_0(f_2)\|\\
				=&{}\|\sum_{g\in G_N}d_g\sum_{h\in G_N}f_1(h)f_2(h)d_h-\sum_{g\in G_N}f_1(g)d_g\sum_{h\in G_N}f_2(h)d_h\|\\
				\leq&{}\|\sum_{g\in G_N}d_g\sum_{h\in G_N}f_1(h)f_2(h)d_h-\sum_{h\in G_N}f_1(h)f_2(h)d_h^2\|\\
				+&{}\|\sum_{h\in G_N}f_1(h)f_2(h)d_h^2-\sum_{g\in G_N}f_1(g)d_g\sum_{h\in G_N}f_2(h)d_h\|\\
				=&{}0.
			\end{align*}
			Let $K=\prod_{m=N+1}^\infty H_m$, then $G=G_N\times K$. Let $\mu$ be the normalized Haar measure on $K$. Then there exists a conditional expectation $P\colon C(G)\to C(G_N)$, given by $P(f)(h)=\int_Kf(h,g)d_\mu(g)$ for $f\in C(G)$ and $h\in G_N$. By identifying $C(G_N)$ as a subalgebra of $C(G)$, which is invariant under the left translation action $g\mapsto \mathrm{Lt}_g$ of the group $G$. Under this identification, the map $P$ is equivariant. It is easy to verify that $d\in A^\alpha$. Therefore, $\psi=\psi_0\circ P\colon C(G)\to \overline{dAd}$ is an equivariant completely positive contractive map with $\psi(1)=d$, i.e., condition (\ref{t8.19.2}) holds.
			
			Let $f_1,f_2\in S$ and $a\in F$, choose $b\in F_0$ and $e_1,e_2\in S_0$ such that
			\[\|b-a\|<\frac{\varepsilon}{8},\ \|e_1-f_1\|<\frac{\varepsilon}{4},\ \text{and} \ \|e_2-f_2\|<\frac{\varepsilon}{4}.\]
			Then we have
			\begin{align*}
				&{}\|\psi(f_1)a-a\psi(f_1)\|=\|\psi_0\circ P(f_1)a-a\psi_0\circ P(f_1)\|\\
				\leq&{}\|\psi_0\circ P(f_1)a-\psi_0\circ P(e_1)a\|+\|\psi_0\circ P(e_1)a-\psi_0\circ P(e_1)b\|\\
				+&{}\|\psi_0\circ P(e_1)b-b\psi_0\circ P(e_1)\|+\|b\psi_0\circ P(e_1)-b\psi_0\circ P(f_1)\|\\
				+&{}\|b\psi_0\circ P(f_1)-a\psi_0\circ P(f_1)\|\\
				<&{}\frac{\varepsilon}{4}+\frac{\varepsilon}{8}+\frac{\varepsilon}{4}+\frac{\varepsilon}{4}+\frac{\varepsilon}{8}=\varepsilon.
			\end{align*}
			And
			\begin{align*}
				&{}\|\psi(1)\psi(f_1f_2)-\psi(f_1)\psi(f_2)\|=\|\psi_0\circ P(1)\psi_0\circ P(f_1f_2)-\psi_0\circ P(f_1)\psi_0\circ P(f_2)\|\\
				\leq&{}\|\psi_0\circ P(1)\psi_0\circ P(f_1f_2)-\psi_0\circ P(1)\psi_0\circ P(e_1e_2)\|\\
				+&{}\|\psi_0\circ P(1)\psi_0\circ P(e_1e_2)-\psi_0\circ P(e_1)\psi_0\circ P(e_2)\|\\
				+&{}\|\psi_0\circ P(e_1)\psi_0\circ P(e_2)-\psi_0\circ P(f_1)\psi_0\circ P(f_2)\|\\
				<&{}\frac{\varepsilon}{2}+\frac{\varepsilon}{2}=\varepsilon.
			\end{align*}
			Hence condition (\ref{t8.19.1}) holds.
			
			For condition (\ref{t8.19.3}), since $d=d^{(1)}\otimes d^{(2)}\otimes\dots\otimes d^{(N)}$, let $\tau$ be the unique tracial state on $A$, by Lemma \ref{TPC}, we have
			\[d_\tau(1-d)\leq\sum_{m=1}^Nd_{\tau_m}(1-d^{(m)})<N\varepsilon_0<\delta_1=d_\tau(x).\]
			Since $A\cong \mathcal{Z}$ has strict comparison, we obtain $(1-d-\varepsilon)_+\precsim_A 1-d\precsim_A x$.
			
			For condition (\ref{t8.19.4}), similarly, let $\tau\in \mathrm{T}(A^\alpha)$, for all $\tau\in \mathrm{T}(A^\alpha)$, since the set of tracial states on the minimal tensor product of C*-algebras is the closed convex hull of the tensor products of tracial states on the component C*-algebras, we may assume $\tau$ is an elementary tensor product. Therefore, we have
			\[d_\tau(1-d)\leq\sum_{m=1}^Nd_{\tau_m}(1-d^{(m)})<N\varepsilon_0<\delta_2=\inf_{\tau\in \mathrm{T}(A^\alpha)}d_\tau(y).\]
			Since $A^\alpha$ has strict comparison, we obtain $(1-d-\varepsilon)_+\precsim_{A^\alpha}1-d\precsim_{A^\alpha}y$.
			
			Moreover, we also have
			\[d_\tau(1-d)\leq\sum_{m=1}^Nd_{\tau_m}(1-d^{(m)})<N\varepsilon_0<\frac{1}{2}.\]
			Since $d_\tau(d)\geq1-d_\tau(1-d)>\frac{1}{2}>d_\tau(1-d)$, and $A^\alpha$ has strict comparison, we obtain $(1-d-\varepsilon)_+\precsim_{A^\alpha}1-d\precsim_{A^\alpha}d$.
			Hence condition (\ref{t8.19.5}) holds.
			
			Since $A\cong\mathcal{Z}$ is stably finite, condition (\ref{t8.19.6}) ($\|dxd\|>1-\varepsilon$) is redundant.
		\end{proof}
	\end{theorem}
	
	\begin{remark}
		In fact, we can use the above technique to construct examples of permutation actions on unital simple tracially $\mathcal{Z}$-stable C*-algebras with few projections. For example, the algebras constructed in \cite{FL2022} are all unital simple tracially $\mathcal{Z}$-stable with few projections, have strict comparison, and stable rank one (hence stably finite). The verification of conditions (\ref{t8.19.3}) - (\ref{t8.19.5}) can also be performed using the approach in the proof of \cite[Theorem 4.8]{MP2025}.
	\end{remark}
	
	
	\subsection*{Funding}
	No funding was received.
	
	
	\bibliographystyle{plain}
	\bibliography{EX5_3}

@article {AGJP2021,
	AUTHOR = {Amini, Massoud and Golestani, Nasser and Jamali, Saeid and
	Phillips, N. Christopher},
	TITLE = {Finite group and integer actions on simple tracially $\mathcal{Z}$-absorbing {C}*-algebras},
	JOURNAL = {J. Operator Theory},
	FJOURNAL = {Journal of Operator Theory},
	VOLUME = {92},
	YEAR = {2024},
	NUMBER = {2},
	PAGES = {505--547},
	ISSN = {0379-4024,1841-7744},
	MRCLASS = {46L55 (46L05 46L40)},
	MRNUMBER = {4849105},
}

@incollection {APT2009,
	AUTHOR = {Ara, Pere and Perera, Francesc and Toms, Andrew S.},
	TITLE = {K-theory for operator algebras. {C}lassification of
	{C}*-algebras},
	BOOKTITLE = {Aspects of operator algebras and applications},
	SERIES = {Contemp. Math.},
	VOLUME = {534},
	PAGES = {1--71},
	PUBLISHER = {Amer. Math. Soc., Providence, RI},
	YEAR = {2011},
	ISBN = {978-0-8218-4905-7},
	MRCLASS = {46L35 (19K14 46L05 46L80)},
	MRNUMBER = {2767222},
	MRREVIEWER = {Vladimir\ Manuilov},
	DOI = {10.1090/conm/534/10521},
	URL = {https://doi.org/10.1090/conm/534/10521},
}

@article {ABP2018,
	AUTHOR = {Archey, Dawn and Buck, Julian and Phillips, N. Christopher},
	TITLE = {Centrally large subalgebras and tracial $\mathcal{Z}$-absorption},
	JOURNAL = {Int. Math. Res. Not. IMRN},
	FJOURNAL = {International Mathematics Research Notices. IMRN},
	YEAR = {2018},
	NUMBER = {6},
	PAGES = {1857--1877},
	ISSN = {1073-7928,1687-0247},
	MRCLASS = {46L05 (46L10 46L35)},
	MRNUMBER = {3801476},
	MRREVIEWER = {Bachir\ Bekka},
	DOI = {10.1093/imrn/rnw292},
	URL = {https://doi.org/10.1093/imrn/rnw292},
}

@article {AP2020,
	AUTHOR = {Archey, Dawn E. and Phillips, N. Christopher},
	TITLE = {Permanence of stable rank one for centrally large subalgebras
	and crossed products by minimal homeomorphisms},
	JOURNAL = {J. Operator Theory},
	FJOURNAL = {Journal of Operator Theory},
	VOLUME = {83},
	YEAR = {2020},
	NUMBER = {2},
	PAGES = {353--389},
	ISSN = {0379-4024,1841-7744},
	MRCLASS = {46L05 (46L55)},
	MRNUMBER = {4078704},
	MRREVIEWER = {Zhuang\ Niu},
	DOI = {10.7900/jot},
	URL = {https://doi.org/10.7900/jot},
}

@article {Asa2024,
	AUTHOR = {Asadi-Vasfi, M. Ali},
	TITLE = {Weakly tracially approximately representable actions},
	JOURNAL = {J. Operator Theory},
	FJOURNAL = {Journal of Operator Theory},
	VOLUME = {91},
	YEAR = {2024},
	NUMBER = {1},
	PAGES = {3--25},
	ISSN = {0379-4024,1841-7744},
	MRCLASS = {46L55 (19K14 46L80)},
	MRNUMBER = {4718795},
	MRREVIEWER = {Xiao\ Chun\ Fang},
}

@article {AGP2021,
	AUTHOR = {Asadi-Vasfi, M. Ali and Golestani, Nasser and Phillips, N.
	Christopher},
	TITLE = {The {C}untz semigroup and the radius of comparison of the
	crossed product by a finite group},
	JOURNAL = {Ergodic Theory Dynam. Systems},
	FJOURNAL = {Ergodic Theory and Dynamical Systems},
	VOLUME = {41},
	YEAR = {2021},
	NUMBER = {12},
	PAGES = {3541--3592},
	ISSN = {0143-3857,1469-4417},
	MRCLASS = {46L55 (19K14 46L80)},
	MRNUMBER = {4336489},
	MRREVIEWER = {Xiao\ Chun\ Fang},
	DOI = {10.1017/etds.2020.121},
	URL = {https://doi.org/10.1017/etds.2020.121},
}

@article {APT2018,
	AUTHOR = {Antoine, Ramon and Perera, Francesc and Thiel, Hannes},
	TITLE = {Tensor products and regularity properties of {C}untz
	semigroups},
	JOURNAL = {Mem. Amer. Math. Soc.},
	FJOURNAL = {Memoirs of the American Mathematical Society},
	VOLUME = {251},
	YEAR = {2018},
	NUMBER = {1199},
	PAGES = {viii+191},
	ISSN = {0065-9266,1947-6221},
	ISBN = {978-1-4704-2797-9; 978-1-4704-4282-8},
	MRCLASS = {06B35 (06F05 15A69 16W80 18B35 19K14 46L05 54F05)},
	MRNUMBER = {3756921},
	DOI = {10.1090/memo/1199},
	URL = {https://doi.org/10.1090/memo/1199},
}

@article {BH1982,
	AUTHOR = {Blackadar, Bruce and Handelman, David},
	TITLE = {Dimension functions and traces on {C}*-algebras},
	JOURNAL = {J. Funct. Anal.},
	FJOURNAL = {Journal of Functional Analysis},
	VOLUME = {45},
	YEAR = {1982},
	NUMBER = {3},
	PAGES = {297--340},
	ISSN = {0022-1236},
	MRCLASS = {46L05},
	MRNUMBER = {650185},
	MRREVIEWER = {J.\ W.\ Bunce},
	DOI = {10.1016/0022-1236(82)90009-X},
	URL = {https://doi.org/10.1016/0022-1236(82)90009-X},
}

@article {Con1975,
	AUTHOR = {Connes, Alain},
	TITLE = {Outer conjugacy classes of automorphisms of factors},
	JOURNAL = {Ann. Sci. \'Ecole Norm. Sup. (4)},
	FJOURNAL = {Annales Scientifiques de l'\'Ecole Normale Sup\'erieure.
	Quatri\`eme S\'erie},
	VOLUME = {8},
	YEAR = {1975},
	NUMBER = {3},
	PAGES = {383--419},
	ISSN = {0012-9593},
	MRCLASS = {46L10},
	MRNUMBER = {394228},
	MRREVIEWER = {Hisashi\ Choda},
	URL = {http://www.numdam.org/item?id=ASENS_1975_4_8_3_383_0},
}

@article {Cun1978,
	AUTHOR = {Cuntz, Joachim},
	TITLE = {Dimension functions on simple {C}*-algebras},
	JOURNAL = {Math. Ann.},
	FJOURNAL = {Mathematische Annalen},
	VOLUME = {233},
	YEAR = {1978},
	NUMBER = {2},
	PAGES = {145--153},
	ISSN = {0025-5831,1432-1807},
	MRCLASS = {46L05},
	MRNUMBER = {467332},
	MRREVIEWER = {J.\ W.\ Bunce},
	DOI = {10.1007/BF01421922},
	URL = {https://doi.org/10.1007/BF01421922},
}

@article {ERS2011,
	AUTHOR = {Elliott, George A. and Robert, Leonel and Santiago, Luis},
	TITLE = {The cone of lower semicontinuous traces on a {C}*-algebra},
	JOURNAL = {Amer. J. Math.},
	FJOURNAL = {American Journal of Mathematics},
	VOLUME = {133},
	YEAR = {2011},
	NUMBER = {4},
	PAGES = {969--1005},
	ISSN = {0002-9327,1080-6377},
	MRCLASS = {46L30},
	MRNUMBER = {2823868},
	DOI = {10.1353/ajm.2011.0027},
	URL = {https://doi.org/10.1353/ajm.2011.0027},
}

@article {FHLRTVW2021,
	AUTHOR = {Farah, Ilijas and Hart, Bradd and Lupini, Martino and Robert,
	Leonel and Tikuisis, Aaron and Vignati, Alessandro and Winter,
	Wilhelm},
	TITLE = {Model theory of {C}*-algebras},
	JOURNAL = {Mem. Amer. Math. Soc.},
	FJOURNAL = {Memoirs of the American Mathematical Society},
	VOLUME = {271},
	YEAR = {2021},
	NUMBER = {1324},
	PAGES = {viii+127},
	ISSN = {0065-9266,1947-6221},
	ISBN = {978-1-4704-4757-1; 978-1-4704-6626-8},
	MRCLASS = {03C98 (03C20 03C25 03E15 46L05 46L35)},
	MRNUMBER = {4279915},
	MRREVIEWER = {Christopher\ J.\ Eagle},
	DOI = {10.1090/memo/1324},
	URL = {https://doi.org/10.1090/memo/1324},
}

@article {FG2020,
	AUTHOR = {Forough, Marzieh and Golestani, Nasser},
	TITLE = {The weak tracial {R}okhlin property for finite group actions
	on simple {C}*-algebras},
	JOURNAL = {Doc. Math.},
	FJOURNAL = {Documenta Mathematica},
	VOLUME = {25},
	YEAR = {2020},
	PAGES = {2507--2552},
	ISSN = {1431-0635,1431-0643},
	MRCLASS = {46L55 (46L05 46L40)},
	MRNUMBER = {4216445},
}

@article {FL2020,
	AUTHOR = {Fu, Xuanlong and Lin, Huaxin},
	TITLE = {Tracial approximation in simple {C}*-algebras},
	JOURNAL = {Canad. J. Math.},
	FJOURNAL = {Canadian Journal of Mathematics. Journal Canadien de
	Math\'ematiques},
	VOLUME = {74},
	YEAR = {2022},
	NUMBER = {4},
	PAGES = {942--1004},
	ISSN = {0008-414X,1496-4279},
	MRCLASS = {46L05 (46L35)},
	MRNUMBER = {4464578},
	DOI = {10.4153/S0008414X21000158},
	URL = {https://doi.org/10.4153/S0008414X21000158},
}

@article {FL2022,
	AUTHOR = {Fu, Xuanlong and Lin, Huaxin},
	TITLE = {Nonamenable simple {C}*-algebras with tracial
	approximation},
	JOURNAL = {Forum Math. Sigma},
	FJOURNAL = {Forum of Mathematics. Sigma},
	VOLUME = {10},
	YEAR = {2022},
	PAGES = {Paper No. e14, 50},
	ISSN = {2050-5094},
	MRCLASS = {46L35 (46L05)},
	MRNUMBER = {4387777},
	DOI = {10.1017/fms.2021.79},
	URL = {https://doi.org/10.1017/fms.2021.79},
}

@article {Gar2014,
	AUTHOR = {Gardella, Eusebio},
	TITLE = {Regularity properties and {R}okhlin dimension for compact
	group actions},
	JOURNAL = {Houston J. Math.},
	FJOURNAL = {Houston Journal of Mathematics},
	VOLUME = {43},
	YEAR = {2017},
	NUMBER = {3},
	PAGES = {861--889},
	ISSN = {0362-1588},
	MRCLASS = {46L35 (46L55)},
	MRNUMBER = {3739037},
	MRREVIEWER = {G\'abor\ Szab\'o},
}

@article {Gar2017,
	AUTHOR = {Gardella, Eusebio},
	TITLE = {Crossed products by compact group actions with the {R}okhlin
	property},
	JOURNAL = {J. Noncommut. Geom.},
	FJOURNAL = {Journal of Noncommutative Geometry},
	VOLUME = {11},
	YEAR = {2017},
	NUMBER = {4},
	PAGES = {1593--1626},
	ISSN = {1661-6952,1661-6960},
	MRCLASS = {46L55 (46L35)},
	MRNUMBER = {3743232},
	MRREVIEWER = {Xiao\ Chun\ Fang},
	DOI = {10.4171/JNCG/11-4-11},
	URL = {https://doi.org/10.4171/JNCG/11-4-11},
}

@article {Gar2019,
	AUTHOR = {Gardella, Eusebio},
	TITLE = {Compact group actions with the {R}okhlin property},
	JOURNAL = {Trans. Amer. Math. Soc.},
	FJOURNAL = {Transactions of the American Mathematical Society},
	VOLUME = {371},
	YEAR = {2019},
	NUMBER = {4},
	PAGES = {2837--2874},
	ISSN = {0002-9947,1088-6850},
	MRCLASS = {46L55 (19K99 37A55 46L35 46L80)},
	MRNUMBER = {3896099},
	MRREVIEWER = {Valentin\ Deaconu},
	DOI = {10.1090/tran/7523},
	URL = {https://doi.org/10.1090/tran/7523},
}

@article {GHS2021,
	AUTHOR = {Gardella, Eusebio and Hirshberg, Ilan and Santiago, Luis},
	TITLE = {Rokhlin dimension: duality, tracial properties, and crossed
	products},
	JOURNAL = {Ergodic Theory Dynam. Systems},
	FJOURNAL = {Ergodic Theory and Dynamical Systems},
	VOLUME = {41},
	YEAR = {2021},
	NUMBER = {2},
	PAGES = {408--460},
	ISSN = {0143-3857,1469-4417},
	MRCLASS = {46L55 (22D25 37A55 46L35 46L80)},
	MRNUMBER = {4177290},
	MRREVIEWER = {Yongle\ Jiang},
	DOI = {10.1017/etds.2019.68},
	URL = {https://doi.org/10.1017/etds.2019.68},
}

@article{GP2024,
	title={The modern theory of {C}untz semigroups of {C}*-algebras},
	author={Gardella, Eusebio and Perera, Francesc},
	journal={EMS Surv. Math. Sci.},
	year={2024},
	DOI ={10.4171/EMSS/84},
	URL=  {https://ems.press/journals/emss/articles/14298154},
}

@article {GS2016,
	AUTHOR = {Gardella, Eusebio and Santiago, Luis},
	TITLE = {Equivariant {$*$}-homomorphisms, {R}okhlin constraints and
	equivariant {UHF}-absorption},
	JOURNAL = {J. Funct. Anal.},
	FJOURNAL = {Journal of Functional Analysis},
	VOLUME = {270},
	YEAR = {2016},
	NUMBER = {7},
	PAGES = {2543--2590},
	ISSN = {0022-1236,1096-0783},
	MRCLASS = {46L55 (46L05)},
	MRNUMBER = {3464050},
	MRREVIEWER = {G\'abor\ Szab\'o},
	DOI = {10.1016/j.jfa.2016.02.010},
	URL = {https://doi.org/10.1016/j.jfa.2016.02.010},
}

@article {GT2022,
	AUTHOR = {Gardella, Eusebio and Thiel, Hannes},
	TITLE = {Weighted homomorphisms between {$\rm C^*$}-algebras},
	JOURNAL = {Doc. Math.},
	FJOURNAL = {Documenta Mathematica},
	VOLUME = {30},
	YEAR = {2025},
	NUMBER = {3},
	PAGES = {587--610},
	ISSN = {1431-0635,1431-0643},
	MRCLASS = {47B65 (47B49)},
	MRNUMBER = {4916104},
	DOI = {10.4171/dm/1008},
	URL = {https://doi.org/10.4171/dm/1008},
}

@article {GLP1994,
	AUTHOR = {Gootman, Elliot C. and Lazar, Aldo J. and Peligrad, Costel},
	TITLE = {Spectra for compact group actions},
	JOURNAL = {J. Operator Theory},
	FJOURNAL = {Journal of Operator Theory},
	VOLUME = {31},
	YEAR = {1994},
	NUMBER = {2},
	PAGES = {381--399},
	ISSN = {0379-4024},
	MRCLASS = {46L55 (22D25)},
	MRNUMBER = {1331784},
	MRREVIEWER = {Erik\ B\'edos},
}

@article {Han1981,
	AUTHOR = {Handelman, David},
	TITLE = {Homomorphisms of {C}*-algebras to finite
	{AW}*-algebras},
	JOURNAL = {Michigan Math. J.},
	FJOURNAL = {Michigan Mathematical Journal},
	VOLUME = {28},
	YEAR = {1981},
	NUMBER = {2},
	PAGES = {229--240},
	ISSN = {0026-2285,1945-2365},
	MRCLASS = {46L05},
	MRNUMBER = {616272},
	MRREVIEWER = {J.\ W.\ Bunce},
	URL = {http://projecteuclid.org/euclid.mmj/1029002512},
}

@article {HJ1982,
	AUTHOR = {Herman, Richard H. and Jones, Vaughan F. R.},
	TITLE = {Period two automorphisms of {UHF} {C}*-algebras},
	JOURNAL = {J. Funct. Anal.},
	FJOURNAL = {Journal of Functional Analysis},
	VOLUME = {45},
	YEAR = {1982},
	NUMBER = {2},
	PAGES = {169--176},
	ISSN = {0022-1236},
	MRCLASS = {46L40},
	MRNUMBER = {647069},
	MRREVIEWER = {G.\ A.\ Elliott},
	DOI = {10.1016/0022-1236(82)90016-7},
	URL = {https://doi.org/10.1016/0022-1236(82)90016-7},
}

@article {HJ1983,
	AUTHOR = {Herman, Richard H. and Jones, Vaughan F. R.},
	TITLE = {Models of finite group actions},
	JOURNAL = {Math. Scand.},
	FJOURNAL = {Mathematica Scandinavica},
	VOLUME = {52},
	YEAR = {1983},
	NUMBER = {2},
	PAGES = {312--320},
	ISSN = {0025-5521,1903-1807},
	MRCLASS = {46L55},
	MRNUMBER = {702960},
	MRREVIEWER = {Colin\ E.\ Sutherland},
	DOI = {10.7146/math.scand.a-12008},
	URL = {https://doi.org/10.7146/math.scand.a-12008},
}

@article {HO2013,
	AUTHOR = {Hirshberg, Ilan and Orovitz, Joav},
	TITLE = {Tracially $\mathcal{Z}$-absorbing {C}*-algebras},
	JOURNAL = {J. Funct. Anal.},
	FJOURNAL = {Journal of Functional Analysis},
	VOLUME = {265},
	YEAR = {2013},
	NUMBER = {5},
	PAGES = {765--785},
	ISSN = {0022-1236,1096-0783},
	MRCLASS = {46L05},
	MRNUMBER = {3063095},
	MRREVIEWER = {Stuart\ A.\ White},
	DOI = {10.1016/j.jfa.2013.05.005},
	URL = {https://doi.org/10.1016/j.jfa.2013.05.005},
}

@article {HP2015,
	AUTHOR = {Hirshberg, Ilan and Phillips, N. Christopher},
	TITLE = {Rokhlin dimension: obstructions and permanence properties},
	JOURNAL = {Doc. Math.},
	FJOURNAL = {Documenta Mathematica},
	VOLUME = {20},
	YEAR = {2015},
	PAGES = {199--236},
	ISSN = {1431-0635,1431-0643},
	MRCLASS = {46L55},
	MRNUMBER = {3398712},
	MRREVIEWER = {G\'abor\ Szab\'o},
}

@article {HW2007,
	AUTHOR = {Hirshberg, Ilan and Winter, Wilhelm},
	TITLE = {Rokhlin actions and self-absorbing {C}*-algebras},
	JOURNAL = {Pacific J. Math.},
	FJOURNAL = {Pacific Journal of Mathematics},
	VOLUME = {233},
	YEAR = {2007},
	NUMBER = {1},
	PAGES = {125--143},
	ISSN = {0030-8730,1945-5844},
	MRCLASS = {46L05 (46L55)},
	MRNUMBER = {2366371},
	MRREVIEWER = {Efren\ Ruiz},
	DOI = {10.2140/pjm.2007.233.125},
	URL = {https://doi.org/10.2140/pjm.2007.233.125},
}

@article {KR2000,
	AUTHOR = {Kirchberg, Eberhard and {R\o rdam}, Mikael},
	TITLE = {Non-simple purely infinite {C}*-algebras},
	JOURNAL = {Amer. J. Math.},
	FJOURNAL = {American Journal of Mathematics},
	VOLUME = {122},
	YEAR = {2000},
	NUMBER = {3},
	PAGES = {637--666},
	ISSN = {0002-9327,1080-6377},
	MRCLASS = {46L05 (46L35)},
	MRNUMBER = {1759891},
	MRREVIEWER = {Judith\ A.\ Packer},
	URL =
	{http://muse.jhu.edu/journals/american_journal_of_mathematics/v122/122.3kirchberg.pdf},
}

@article {KR2002,
	AUTHOR = {Kirchberg, Eberhard and {R\o rdam}, Mikael},
	TITLE = {Infinite non-simple {C}*-algebras: absorbing the {C}untz
	algebras $\mathcal{O}_\infty$},
	JOURNAL = {Adv. Math.},
	FJOURNAL = {Advances in Mathematics},
	VOLUME = {167},
	YEAR = {2002},
	NUMBER = {2},
	PAGES = {195--264},
	ISSN = {0001-8708,1090-2082},
	MRCLASS = {46L05 (46L35)},
	MRNUMBER = {1906257},
	MRREVIEWER = {Judith\ A.\ Packer},
	DOI = {10.1006/aima.2001.2041},
	URL = {https://doi.org/10.1006/aima.2001.2041},
}

@article {Kis1981,
	AUTHOR = {Kishimoto, Akitaka},
	TITLE = {Outer automorphisms and reduced crossed products of simple
	{C}*-algebras},
	JOURNAL = {Comm. Math. Phys.},
	FJOURNAL = {Communications in Mathematical Physics},
	VOLUME = {81},
	YEAR = {1981},
	NUMBER = {3},
	PAGES = {429--435},
	ISSN = {0010-3616,1432-0916},
	MRCLASS = {46L40 (46L05)},
	MRNUMBER = {634163},
	MRREVIEWER = {Dorte\ Olesen},
	URL = {http://projecteuclid.org/euclid.cmp/1103920327},
}

@book {Lin2001,
	AUTHOR = {Lin, Huaxin},
	TITLE = {An introduction to the classification of amenable
	{C}*-algebras},
	PUBLISHER = {World Scientific Publishing Co., Inc., River Edge, NJ},
	YEAR = {2001},
	PAGES = {xii+320},
	ISBN = {981-02-4680-3},
	MRCLASS = {46Lxx (19K35)},
	MRNUMBER = {1884366},
	MRREVIEWER = {Xiao\ Chun\ Fang},
	DOI = {10.1142/9789812799883},
	URL = {https://doi.org/10.1142/9789812799883},
}

@article {Lin2023,
	AUTHOR = {Lin, Huaxin},
	TITLE = {Hereditary uniform property {$\Gamma$}},
	JOURNAL = {Sci. China Math.},
	FJOURNAL = {Science China. Mathematics},
	VOLUME = {66},
	YEAR = {2023},
	NUMBER = {8},
	PAGES = {1813--1830},
	ISSN = {1674-7283,1869-1862},
	MRCLASS = {46L35 (46L05)},
	MRNUMBER = {4621130},
	MRREVIEWER = {Jorge\ Castillejos},
	DOI = {10.1007/s11425-022-2005-x},
	URL = {https://doi.org/10.1007/s11425-022-2005-x},
}

@article {MS2012b,
	AUTHOR = {Matui, Hiroki and Sato, Yasuhiko},
	TITLE = {$\mathcal{Z}$-stability of crossed products by strongly outer
	actions},
	JOURNAL = {Comm. Math. Phys.},
	FJOURNAL = {Communications in Mathematical Physics},
	VOLUME = {314},
	YEAR = {2012},
	NUMBER = {1},
	PAGES = {193--228},
	ISSN = {0010-3616,1432-0916},
	MRCLASS = {46L55},
	MRNUMBER = {2954514},
	MRREVIEWER = {Sriwulan\ Adji},
	DOI = {10.1007/s00220-011-1392-9},
	URL = {https://doi.org/10.1007/s00220-011-1392-9},
}

@ARTICLE{MP2021,
	author = {{Mohammadkarimi}, Javad and {Phillips}, N. Christopher},
	title = "{Compact Group Actions with the Tracial Rokhlin Property}",
	journal = {arXiv e-prints},
	year = 2021,
	month = oct,
	eid = {arXiv:2110.12135},
	pages = {arXiv:2110.12135},
	doi = {10.48550/arXiv.2110.12135},
	archivePrefix = {arXiv},
	eprint = {2110.12135},
	primaryClass = {math.OA},
	adsurl = {https://ui.adsabs.harvard.edu/abs/2021arXiv211012135M}
}

@ARTICLE{MP2025,
	author = {{Mohammadkarimi}, Javad and {Phillips}, N. Christopher},
	title = "{Compact Group Actions with the Tracial Rokhlin Property II: Examples and Nonexistence Theorems}",
	journal = {arXiv e-prints},
	year = 2025,
	month = may,
	eid = {arXiv:2505.04661},
	pages = {arXiv:2505.04661},
	doi = {10.48550/arXiv.2505.04661},
	archivePrefix = {arXiv},
	eprint = {2505.04661},
	primaryClass = {math.OA},
	adsurl = {https://ui.adsabs.harvard.edu/abs/2025arXiv250504661M}
}

@article {OP2012,
	AUTHOR = {Osaka, Hiroyuki and Phillips, N. Christopher},
	TITLE = {Crossed products by finite group actions with the {R}okhlin
	property},
	JOURNAL = {Math. Z.},
	FJOURNAL = {Mathematische Zeitschrift},
	VOLUME = {270},
	YEAR = {2012},
	NUMBER = {1-2},
	PAGES = {19--42},
	ISSN = {0025-5874,1432-1823},
	MRCLASS = {46L55 (46L35)},
	MRNUMBER = {2875821},
	MRREVIEWER = {Ilan\ Hirshberg},
	DOI = {10.1007/s00209-010-0784-4},
	URL = {https://doi.org/10.1007/s00209-010-0784-4},
}

@book {Phi2006,
	AUTHOR = {Phillips, N. Christopher},
	TITLE = {Equivariant {K}-theory and freeness of group actions on
	{C}*-algebras},
	SERIES = {Lecture Notes in Mathematics},
	VOLUME = {1274},
	PUBLISHER = {Springer-Verlag, Berlin},
	YEAR = {1987},
	PAGES = {viii+371},
	ISBN = {3-540-18277-2},
	MRCLASS = {46L80 (19L47 22D25 46L55 46M20 58G12)},
	MRNUMBER = {911880},
	MRREVIEWER = {Judith\ A.\ Packer},
	DOI = {10.1007/BFb0078657},
	URL = {https://doi.org/10.1007/BFb0078657},
}

@article {Phi2011,
	AUTHOR = {Phillips, N. Christopher},
	TITLE = {The tracial {R}okhlin property for actions of finite groups on
	{C}*-algebras},
	JOURNAL = {Amer. J. Math.},
	FJOURNAL = {American Journal of Mathematics},
	VOLUME = {133},
	YEAR = {2011},
	NUMBER = {3},
	PAGES = {581--636},
	ISSN = {0002-9327,1080-6377},
	MRCLASS = {46L55 (46L40)},
	MRNUMBER = {2808327},
	MRREVIEWER = {Hiroyuki\ Osaka},
	DOI = {10.1353/ajm.2011.0016},
	URL = {https://doi.org/10.1353/ajm.2011.0016},
}

@ARTICLE{Phi2012,
	author = {{Phillips}, N. Christopher},
	title = "{The tracial Rokhlin property is generic}",
	journal = {arXiv e-prints},
	year = 2012,
	month = sep,
	eid = {arXiv:1209.3859},
	pages = {arXiv:1209.3859},
	doi = {10.48550/arXiv.1209.3859},
	archivePrefix = {arXiv},
	eprint = {1209.3859},
	primaryClass = {math.OA},
	adsurl = {https://ui.adsabs.harvard.edu/abs/2012arXiv1209.3859P}
}

@ARTICLE{Phi2014,
	author = {{Phillips}, N. Christopher},
	title = "{Large subalgebras}",
	journal = {arXiv e-prints},
	year = 2014,
	month = aug,
	eid = {arXiv:1408.5546},
	pages = {arXiv:1408.5546},
	doi = {10.48550/arXiv.1408.5546},
	archivePrefix = {arXiv},
	eprint = {1408.5546},
	primaryClass = {math.OA},
	adsurl = {https://ui.adsabs.harvard.edu/abs/2014arXiv1408.5546P}
}

@article {Ros1979,
	AUTHOR = {Rosenberg, Jonathan},
	TITLE = {Appendix to: ``{C}rossed products of {UHF} algebras by product
	type actions''\ [{D}uke {M}ath. {J}. {\bf 46} (1979), no. 1,
	1--23;\ {MR} 82a:46063 above]\ by {O}. {B}ratteli},
	JOURNAL = {Duke Math. J.},
	FJOURNAL = {Duke Mathematical Journal},
	VOLUME = {46},
	YEAR = {1979},
	NUMBER = {1},
	PAGES = {25--26},
	ISSN = {0012-7094,1547-7398},
	MRCLASS = {46L05},
	MRNUMBER = {523599},
	MRREVIEWER = {\c S.\ Str\u atil\u a},
	URL = {http://projecteuclid.org/euclid.dmj/1077313252},
}

@article {Ror1992,
	AUTHOR = {{R\o rdam}, Mikael},
	TITLE = {On the structure of simple {C}*-algebras tensored with a
	{UHF}-algebra. {II}},
	JOURNAL = {J. Funct. Anal.},
	FJOURNAL = {Journal of Functional Analysis},
	VOLUME = {107},
	YEAR = {1992},
	NUMBER = {2},
	PAGES = {255--269},
	ISSN = {0022-1236,1096-0783},
	MRCLASS = {46L05 (46L85)},
	MRNUMBER = {1172023},
	MRREVIEWER = {Mahmood\ Khoshkam},
	DOI = {10.1016/0022-1236(92)90106-S},
	URL = {https://doi.org/10.1016/0022-1236(92)90106-S},
}

@book {RS2002,
	AUTHOR = {{R\o rdam}, M. and {St\o rmer}, E.},
	TITLE = {Classification of nuclear {C}*-algebras. {E}ntropy in
	operator algebras},
	SERIES = {Encyclopaedia of Mathematical Sciences},
	VOLUME = {126},
	NOTE = {Operator Algebras and Non-commutative Geometry, 7},
	PUBLISHER = {Springer-Verlag, Berlin},
	YEAR = {2002},
	PAGES = {x+198},
	ISBN = {3-540-42305-X},
	MRCLASS = {46Lxx (46-06)},
	MRNUMBER = {1878881},
	DOI = {10.1007/978-3-662-04825-2},
	URL = {https://doi.org/10.1007/978-3-662-04825-2},
}

@article {San2015,
	AUTHOR = {Santiago, Luis},
	TITLE = {Crossed products by actions of finite groups with the
	{R}okhlin property},
	JOURNAL = {Internat. J. Math.},
	FJOURNAL = {International Journal of Mathematics},
	VOLUME = {26},
	YEAR = {2015},
	NUMBER = {7},
	PAGES = {1550042, 31},
	ISSN = {0129-167X,1793-6519},
	MRCLASS = {46L55 (46L05)},
	MRNUMBER = {3357031},
	MRREVIEWER = {Jia-jie\ Hua},
	DOI = {10.1142/S0129167X15500421},
	URL = {https://doi.org/10.1142/S0129167X15500421},
}

@article {TF2025,
	AUTHOR = {Tian, Haotian and Fang, Xiaochun},
	TITLE = {Some PERMANENCE PROPERTIES FOR CROSSED PRODUCTS
	BY COMPACT GROUP ACTIONS WITH THE TRACIAL
	{R}OKHLIN PROPERTY},
	JOURNAL = {Rocky Mountain J. Math.},
	FJOURNAL = {The Rocky Mountain Journal of Mathematics},
	VOLUME = {55},
	YEAR = {2025},
	NUMBER = {3},
	PAGES = {847--866},
	ISSN = {0035-7596,1945-3795},
	MRCLASS = {99-06},
	MRNUMBER = {4927094},
	DOI = {10.1216/rmj.2025.55.847},
	URL = {https://doi.org/10.1216/rmj.2025.55.847},
}

@article {TT2015,
	AUTHOR = {Tikuisis, Aaron Peter and Toms, Andrew},
	TITLE = {On the structure of {C}untz semigroups in (possibly) nonunital
	{C}*-algebras},
	JOURNAL = {Canad. Math. Bull.},
	FJOURNAL = {Canadian Mathematical Bulletin. Bulletin Canadien de
	Math\'ematiques},
	VOLUME = {58},
	YEAR = {2015},
	NUMBER = {2},
	PAGES = {402--414},
	ISSN = {0008-4395,1496-4287},
	MRCLASS = {46L35 (46L80 47L40)},
	MRNUMBER = {3334936},
	MRREVIEWER = {Daniel\ Belti\c t\u a},
	DOI = {10.4153/CMB-2014-040-5},
	URL = {https://doi.org/10.4153/CMB-2014-040-5},
}

@article {Tom2006,
	AUTHOR = {Toms, Andrew S.},
	TITLE = {Flat dimension growth for {C}*-algebras},
	JOURNAL = {J. Funct. Anal.},
	FJOURNAL = {Journal of Functional Analysis},
	VOLUME = {238},
	YEAR = {2006},
	NUMBER = {2},
	PAGES = {678--708},
	ISSN = {0022-1236,1096-0783},
	MRCLASS = {46L05},
	MRNUMBER = {2253738},
	MRREVIEWER = {Nadia\ S.\ Larsen},
	DOI = {10.1016/j.jfa.2006.01.010},
	URL = {https://doi.org/10.1016/j.jfa.2006.01.010},
}

@article {TW2007,
	AUTHOR = {Toms, Andrew S. and Winter, Wilhelm},
	TITLE = {Strongly self-absorbing {C}*-algebras},
	JOURNAL = {Trans. Amer. Math. Soc.},
	FJOURNAL = {Transactions of the American Mathematical Society},
	VOLUME = {359},
	YEAR = {2007},
	NUMBER = {8},
	PAGES = {3999--4029},
	ISSN = {0002-9947,1088-6850},
	MRCLASS = {46L05 (46L35)},
	MRNUMBER = {2302521},
	MRREVIEWER = {Hua\ Xin\ Lin},
	DOI = {10.1090/S0002-9947-07-04173-6},
	URL = {https://doi.org/10.1090/S0002-9947-07-04173-6},
}

@book {Wan2013,
	AUTHOR = {Wang, Qingyun},
	TITLE = {Tracial {R}okhlin {P}roperty and {N}on-{C}ommutative
	{D}imensions},
	NOTE = {Thesis (Ph.D.)--Washington University in St. Louis},
	PUBLISHER = {ProQuest LLC, Ann Arbor, MI},
	YEAR = {2013},
	PAGES = {87},
	ISBN = {978-1303-02778-9},
	MRCLASS = {99-05},
	MRNUMBER = {3153204},
	URL =
	{http://gateway.proquest.com/openurl?url_ver=Z39.88-2004&rft_val_fmt=info:ofi/fmt:kev:mtx:dissertation&res_dat=xri:pqm&rft_dat=xri:pqdiss:3558358},
}

@article {Wan2018,
	AUTHOR = {Wang, Qingyun},
	TITLE = {The tracial {R}okhlin property for actions of amenable groups
	on {C}*-algebras},
	JOURNAL = {Rocky Mountain J. Math.},
	FJOURNAL = {The Rocky Mountain Journal of Mathematics},
	VOLUME = {48},
	YEAR = {2018},
	NUMBER = {4},
	PAGES = {1307--1344},
	ISSN = {0035-7596,1945-3795},
	MRCLASS = {46L55},
	MRNUMBER = {3859760},
	MRREVIEWER = {Roberto\ Quezada},
	DOI = {10.1216/RMJ-2018-48-4-1307},
	URL = {https://doi.org/10.1216/RMJ-2018-48-4-1307},
}

@article {Win2010,
	AUTHOR = {Winter, Wilhelm},
	TITLE = {Decomposition rank and $\mathcal{Z}$-stability},
	JOURNAL = {Invent. Math.},
	FJOURNAL = {Inventiones Mathematicae},
	VOLUME = {179},
	YEAR = {2010},
	NUMBER = {2},
	PAGES = {229--301},
	ISSN = {0020-9910,1432-1297},
	MRCLASS = {46L35 (46L05)},
	MRNUMBER = {2570118},
	MRREVIEWER = {Hua\ Xin\ Lin},
	DOI = {10.1007/s00222-009-0216-4},
	URL = {https://doi.org/10.1007/s00222-009-0216-4},
}

@article {WZ2009,
	AUTHOR = {Winter, Wilhelm and Zacharias, Joachim},
	TITLE = {Completely positive maps of order zero},
	JOURNAL = {M\"unster J. Math.},
	FJOURNAL = {M\"unster Journal of Mathematics},
	VOLUME = {2},
	YEAR = {2009},
	PAGES = {311--324},
	ISSN = {1867-5778,1867-5786},
	MRCLASS = {46L05},
	MRNUMBER = {2545617},
	MRREVIEWER = {William\ Paschke},
}
\end{document}